\documentclass[a4paper,onecolumn,11pt]{article}
\usepackage{geometry}
\usepackage{amsmath}
\usepackage{indentfirst}
\usepackage{ntheorem}
\usepackage{authblk}
\usepackage{bbm}
\usepackage{graphicx}
\usepackage{multirow}
\numberwithin{equation}{section}
\newtheorem{thm}{Theorem}[section]
\newtheorem{prop}[thm]{Proposition}
\newtheorem{lem}[thm]{Lemma}

\newcommand{\RNum}[1]{\uppercase\expandafter{\romannumeral #1\relax}}
\usepackage{bm}
\usepackage{amssymb,amsmath,mathrsfs,txfonts}

\usepackage[colorlinks=true]{hyperref}
\hypersetup{urlcolor=blue, citecolor=red}

\title{\bf Asymptotic stability of planar stationary solution to a 2D hyperbolic-elliptic coupled system of the radiating gas in half space}
\author[1]{Minyi Zhang}
\author[2]{Changjiang Zhu\thanks{Corresponding author. Email: machjzhu@scut.edu.cn}}
\affil[1,2]{School of Mathematics, South China University of Technology, Guangzhou 510641, P.R. China}
\date{} 

\begin{document}
\maketitle

\begin{abstract}
  \indent This paper is concerned with the asymptotic stability of planar stationary solution to an initial-boundary value problem for a two-dimensional hyperbolic-elliptic coupled system of the radiating gas in half space. We show that the solution to the problem converges to the corresponding planar stationary solution as time tends to infinity under small initial perturbation. These results are proved by the standard $L^2$-energy method. Moreover, we prove that the solution $(u,q)$ converges to the corresponding planar stationary solution at the rate $t^{-\alpha/2-1/4}$ for non-degenerate case, and $t^{-1/4}$ for degenerate case. The proof is based on the time and space weighted energy method.

  \vspace{6mm}


\vspace{3mm}


\end{abstract}

\hypersetup{%
  linkcolor  = black
}
\hypersetup{%
  citecolor  = black
}

\tableofcontents

\section{Introduction}
In this paper, we consider the stability of a planar stationary solution to the following initial boundary value problem for a two-dimensional hyperbolic-elliptic coupled system of the radiating gas in half space $\mathbbm{R}_+\times\mathbbm{R}$:

\begin{equation}\label{yfc}
	\begin{cases}
		u_t+f(u)_x+g(u)_y+{\rm{div}} q=0,\\
		-\nabla {\rm{div}} q+q+\nabla u=0,
	\end{cases}
\end{equation}
where $\mathbbm{R}_+:=(0,\infty)$, $u(x,y,t)$ is a scalar function of position $(x,y)\in \mathbbm{R}_+\times\mathbbm{R}$ and time $t\ge0$, $q=(q_1,q_2)$ is a vector function, $f$ and $g$ are smooth functions. We assume that $f$ is strictly convex, i.e. there exists a positive constant $\alpha$ such that for any $u\in \mathbbm{R}$
\begin{equation}\label{eq-f-tj}
f''(u)\ge \alpha>0.
\end{equation}
We also assume that
\begin{equation*}
f(0)=f'(0)=0.
\end{equation*}
The initial and boundary conditions for the equation $\eqref{yfc}$ are prescribed as
\begin{equation}\label{yfctj}
u(x,y,0)=u_0(x,y), \ \ \ \ u(0,y,t)=u_-,
\end{equation}
where $u_-$ is a constant. The initial function $u_0(x,y)$ satisfies
\begin{equation}\label{yfcu0}
	u_0(x,y)\rightarrow u_+ \ \text{ as } \ x\rightarrow +\infty.
\end{equation}
Throughout this paper, we assume that $u_-<u_+\le0$. The inequality $u_-<0$ implies that the fluid blows out through the boundary $\{x=0\}$, which yields that the characteristics of the hyperbolic system $\eqref{yfc}$ around the boundary are negative. Therefore, the boundary condition $\eqref{yfctj}_2$ is necessary and sufficient for the well-posedness. This initial-boundary value problem is called the outflow problem.\\
\indent The system $\eqref{yfc}$ is a simplified version of a radiating gas model in two-dimension. The model $\eqref{yfc}$ gives a good approximation to the equations of radiation hydrodynamics. The derivation refers to \cite{DiFrancesco, Gao2008, Hamer, Vincenti1965}.

As $\eqref{yfc}_2$ is a linear elliptic, Liu and Kawashima in \cite{Liu2011} showed that $q$ can be expressed in format of $\nabla u$, that is $q=-(1-\Delta)^{-1}\nabla{u}$. Therefore, the system $\eqref{yfc}$ is formally equivalent to the Burgers equation with non-linear term
\begin{equation}\label{yfc-dengjia}
	u_t+f(u)_x+g(u)_y-(1-\Delta)^{-1}\Delta u=0.
\end{equation}
Since the system $\eqref{yfc-dengjia}$ is similar to a viscous conservation laws in format, many methods of energy estimate for viscous conservation laws can be used for reference in system $\eqref{yfc-dengjia}$ which is equivalent to the system $\eqref{yfc}$.
And we know that the stability of nonlinear waves to the Cauchy problem of the following one-dimensional viscous conservation laws $\eqref{eq-viscousconlaw}$ with $n=1$ has been studied thoroughly, after the pioneering work in 1960 by Il'in and Oleinik in $\cite{Ilin1960}$:
\begin{equation}\label{eq-viscousconlaw}
	u_t+\sum_{j=1}^nf(u)_x=\Delta u, \quad \text{for } \ (x,t)\in \mathbbm{R}^n\times(0,\infty).
\end{equation}
In the cases of one-dimensional whole space, Kawashima, Matsumura and Nishihara studied the asymptotic stability of traveling waves, cf. \cite{Kawashima1985} \cite{Nishihara1985}; Harabetian, Hattori, Matsumura and Nishihara proved the asymptotic stability of rarefaction waves, cf. \cite{Harabetian1988,Hattori1991,Matsumura1986,Matsumura1992}. In the cases of two-dimensional whole space, Xin in \cite{Xin1990} first showed the stability of the planar rarefaction wave. Ito in \cite{Ito1996} obtained the convergence rate toward the planar rarefaction wave. After that, Nishikawa \cite{Nishikawa2000} improved their results with smallness condition. Furthermore, Kawashima, Nishibata and Nishikawa extended the $L^2$ energy method to general $L^p$ space for multi-dimensional viscous conversation laws. Their techniques have been applied to many other models such as the compressible Navier-Stokes equations and hyperbolic-elliptic coupled system.

 The stability of nonlinear waves to the Cauchy problem of the hyperbolic-elliptic coupled system for one-dimensional radiating gas has been investigated thoroughly. We refer to \cite{Kawashima1985,kawashima1999,kawashima1998,kawashima99,Nishibata2000,Nishikawa2007,Lattanzio2007,Lattanzio2009,Nguyen2010} for the work of shock waves, and we refer to \cite{Kawashima2004,Ruan2007} for rarefaction waves and \cite{Tanaka1995,Iguchi2002} for diffusion waves. In the multi-dimensional case, Di Francesco in \cite{DiFrancesco} showed the global-in-time existence and the uniqueness of weak solutions to the system $\eqref{yfc}$ and analyzed the relaxation limits. In \cite{Wang2009}, Wang and Wang obtained the global existence and pointwise estimates of the solution in multi-dimensions by employing the method of the Green's function combined with some energy estimates. Gao and Zhu in \cite{Gao2008} studied the stability of solutions to the Cauchy problem $\eqref{yfc}$ toward the planar rarefaction waves in two-dimension with the far field states condition $u_-<u_+$, and they obtained the asymptotic time-decay rates. And in \cite{Gao}, Gao, Ruan and Zhu further extended to $3\le n\le5$ in the spatial dimension. In \cite{Ruan2010}, Ruan and Zhu proved the existence and uniqueness of the global solutions and obtained the $L^p$-convergence rates toward the diffusion waves for $1\le n<8$. In \cite{Liu2011}, Liu and Kawashima obtained the global existence and optimal time-decay rates of solutions to the linear diffusion waves, which is given in terms of the heat kernel, for $n\ge2$ with small initial data. More recently, Wang and Wang in \cite{Wang2013} investigated the blow up of classical solutions to the nonlinear Cauchy problem with large initial data by employing the characteristic method. And further they proved that under some restrict conditions on the initial data $u_0\in H^2$, the global existence theory is established.\\
\indent In the cases of one-dimensional half space, Liu, Matsumura and Nishihara in \cite{Liu1997,Liu1998} first considered the asymptotic stability to the initial-boundary value problem for viscous conversation laws. They showed that the stability of several types of nonlinear waves in one-dimensional half space. In \cite{Liu1998}, the asymptotic states of the solutions of the system $\eqref{eq-viscousconlaw}$ with $n=1$ were classified into the following three cases depending on the signatures of $f'(u_\pm)$:
\begin{equation*}
	\begin{aligned}
		(a)\ \ f'(u_-)<f'(u_+)\le0;\\
		(b)\ \  f'(u_-)<0<f'(u_+);\\
		(c)\ \  0\le f'(u_-)<f'(u_+).
	\end{aligned}
\end{equation*}
 In case $(a)$, the solution of $\eqref{eq-viscousconlaw}$ converges to stationary solution. In case $(b)$, the solution behaviors as the superposition of stationary solution and rarefaction wave. In case $(c)$, the asymptotic state is rarefaction wave. From $\eqref{eq-f-tj}$, case $(a)$ is equivalent to $u_-<u_+\le 0$, case $(b)$ is equivalent to $u_-<0<u_+$ and case $(c)$ is equivalent to $0\le u_-<u_+$. After the work of \cite{Liu1998}, a series of related studies have been done. Xin in \cite{Xin1990} first showed the asymptotic stability of planar rarefaction waves for viscous conversation laws $\eqref{eq-viscousconlaw}$ in two dimensions, and then Ito \cite{Ito1996} and Nishikawa \cite{Nishikawa2000} gained the decay rate. Recently, Kawashima \cite{Kawashima2003} first showed the asymptotic stability of planar stationary solution for viscous conversation laws in two-dimension half space and obtained the convergence rate. More recently, Ruan and Zhu \cite{Ruan2008} studied the initial-boundary value problem for a hyperbolic-elliptic system in one-dimensional half space to a rarefaction wave for the case $0=u_-<u_+$. Ji, Zhang and Zhu \cite{Ji} further investigated the asymptotic stability of the stationary solution, rarefaction wave, and the superposition of the rarefaction wave and the stationary solution for the case $u_-<u_+\le 0$, $0\le u_-<u_+$ and $u_-<0<u_+$, respectively. For the initial-boundary value problem, there are also many results for other physical models such as the compressible Navier-Stokes equations. Kawashima, Nishibata and Zhu in \cite{Kawashima2003NS} studied the existence and the asymptotic stability of a stationary solution for the compressible Navier-Stokes equations in one-dimensional half space and then the convergence rate was obtained in \cite{Nakamura2007}. Kagei and Kawashima in \cite{Kagei2006} investigated the stability of planar stationary solutions for the barotropic compressible Navier-Stokes equations in the multi-dimensional half space and the decay rate is given in \cite{Nakamura2009}.\\
 \indent Motivated by these preceding results, especially \cite{Liu1998}, \cite{Kawashima2003} and \cite{Ji}, we will study the asymptotic stability of planar stationary solution for system $\eqref{yfc}$-$\eqref{yfcu0}$. We attempt to extend the results in \cite{Kawashima2003} that the vicious conservation laws converges to the planar rarefaction wave (including asymptotic stability and decay rate) to the equations $\eqref{yfc}$. However, the one-dimensional $H^2$ energy estimate in \cite{Ji} cannot be adapted directly to the two-dimension case for $H^3$ Sobolev space caused by the \emph {a priori} assumption $\eqref{eq-pri-assump}$. Fortunately, we can overcome these difficulties and get the stability result by using the elementary energy methods and conducting suitable treatment on the boundary ${x=0}$.

 The main difficulties in the present paper are as follows.
 \begin{itemize}
 \item[(1)] In order to close the \emph{a priori} assumption $\eqref{eq-pri-assump}$ by applying Gagliardo-Nirenberg inequality, we need to prove \emph{a priori} estimates in $H^3$-Sobolev space. To gain the $H^3$-estimates, we have to deal with the boundary effect of higher order partial derivatives of $v$ and $p$ at the boundary $x=0$, which is a difficult problem. Thanks to $u_-<0$, in proving $H^3$-estimates of $v$, some additional boundary terms are provided:
     \begin{equation}\label{1.7}
     \int_0^t\int_{\mathbbm{R}}\int_{\mathbbm{R}_+} (Kv(0,y,t))^2 \,\mathrm{d}x \mathrm{d}y \mathrm{d}\tau,
     \end{equation}
     where the operator $K=\nabla, \Delta, \nabla\Delta$. However, there is no boundary good term of $p$ in $\eqref{1.7}$. To overcome it, we convert from the boundary term containing $p$ to the boundary term of $v$ by making full use of the structure of perturbation equations $\eqref{eq-vp-rd}$. For example, in proving $H^1$-estimates of $v$, by utilizing the formation $v\times\eqref{eq-vp-rd}_1+p\cdot \eqref{eq-vp-rd}_2+\nabla{v}\cdot \nabla \eqref{eq-vp-rd}_1- \nabla{\rm{div}}{p}\cdot \eqref{eq-vp-rd}_2$, we make the terms containing $p$ into the form of $\eqref{eq-vp-rd}_2$ to eliminate the terms containing $p$ and get a good term $|\nabla v|^2$ (can see $\eqref{eq-vnabv}$). Moreover, the boundary term ${\rm{div}}{p}(0,y,t)$ will appear in estimating the second derivative of $v$ (see $\eqref{eq-dlnabv-jf}$). We can get the relation $\eqref{eq-vx0-bj}$ from $\eqref{eq-vp-rd}_1$, and thereby convert from ${\rm{div}}{p}(0,y,t)$ to $v_x(0,y,t)$.
  \item[(2)]  We know that \cite{Ji, Kawashima2003} require only to prove the \emph{a priori} estimates in $H^2$ Sololev space to close the energy estimates, but $H^3$-estimates is required in the present paper. From $\eqref{eq-vp-rd}_1$ and boundary condition $v(0,y, t)=0$, we obtain the boundary information $\mathrm{div}p(0,y,t)$, but not $\nabla p(0,y,t)$. Therefore, to get the estimates of $v$, we apply operator $\nabla$, ${\rm{div}}\nabla$, $\nabla{\rm{div}}\nabla$ to $\eqref{eq-vp-rd}_1$ successively, and apply operator ${\rm{div}}$, $\nabla{\rm{div}}$, $\nabla{\rm{div}}$ to $\eqref{eq-vp-rd}_2$ successively, which yields the $L^2$-estimates of $v$, $\Delta v$ and $\nabla\Delta v$. Since $\|\Delta v\|_{L^2}^2=\|\nabla^2 v\|_{L^2}^2$, the $H^2$-estimates of $v$ can be obtained automatically. However, we still do not have $H^3$-estimates of $v$ owing to $\|\nabla\Delta v\|_{L^2}^2\not=\|\nabla^3 v\|_{L^2}^2$. This is the difference and new difficulty between this paper and \cite{Ji, Kawashima2003}. To show the estimate of the third order derivatives of $v$, from $v_{yy}(0,y,t)=0$, we can get further a new $L^2$-estimate of $\nabla v_{yy}$. By utilizing the equivalence relation between $\|  \nabla^3{v} \|_{L^2}$ and $\|  \nabla\Delta v \|_{L^2}+ \| \nabla v_{yy} \|_{L^2} $, the $L^2$-estimate of $\nabla^3 v$ can be obtained.
  \item[(3)]  Since we are dealing with the initial boundary value problem, we can not derive the $H^2$-estimates of $p$ directly by applying elliptic estimates $\| p \|_{H^2}\le C \| v \|_{L^2} $ from the elliptic equation $\eqref{eq-vp-rd}_2$ like Cauchy problem. Therefore, we require to find new ways to estimate $p$ and its low order derivatives. To overcome this difficulty, the $H^2$-estimates of ${\rm{div}}{p}$ can be obtained under the \emph{a priori} assumption $\eqref{eq-vt-prioriass}$ that requires small $\|v_t\|_{H^2}$. According to perturbation equation $\eqref{eq-vp-rd}_1$, we get the $H^2$-estimates of $\mathrm{div}p$, and then combining $\eqref{eq-vp-rd}_2$ we further obtain the $H^1$-estimates of $p$.
  \item[(4)] The fourth difficulty is how to estimate the higher order derivatives of $p$ and ${\rm{div}}{p}$. First of all, in order to give estimate of $\|\nabla^3{\rm{div}}{p}\|_{L^2} $, we notice that $\|\nabla^3{\rm{div}}{p}\|_{L^2}$ is equivalent to $\| \nabla\Delta{\rm{div}}{p} \|_{L^2} +\| \nabla{\rm{div}}{p}_{xx} \|_{L^2}$ by the definition of $\|\nabla^3{\rm{div}}{p}\|_{L^2} $. Since the estimate of $\| \nabla\Delta{\rm{div}}{p} \|_{L^2}$ is not difficult to obtain, we just need to show the estimate of $\| \nabla{\rm{div}}{p}_{xx} \|_{L^2} $ to complete the proof of $\|\nabla^3{\rm{div}}{p}\|_{L^2} $. It is inevitable that we need to deal with a boundary term $p_{1xx}(0,y,t)$. To overcome it, rewriting the perturbation equation $\eqref{eq-vp-rd}_2$ into $\eqref{eq-rd2-dengjia}$ and then applying operator $\nabla \partial_x$ to $\eqref{eq-rd2-dengjia}_1$, we can get some good terms ($\| \nabla{\rm{div}}{p}_{xx} \|_{L^2}^2, \| \nabla p_{1x} \|_{L^2}^2,  \| p_{1xxx} \|_{L^2}^2 $). These good terms can be used to control the boundary effect produced by the boundary term $p_{1xx}(0,y,t)$ (the details can see $\eqref{eq-divp-H3-3}$-$\eqref{eq-divp-H3-5}$). Secondly, in order to get the estimate of $\| \nabla^3{p} \|_{L^2} $, by making full use of $\eqref{eq-vp-rd}_2$ and $\eqref{eq-p1p2}$, we derive the relation between the third partial derivatives of $p_1$ and $p_2$: $\eqref{eq-p1xxx}$ and $\eqref{eq-p2xxx}$. By utilizing the $L^2$-estimate of second order derivatives of ${\rm{div}}{p}$, proving the $L^2$-estimate of $\nabla^3 p$ is converted to showing the $L^2$-estimate of $p_{1xxx}$ and $p_{2xyy}$, which reduces the complexity of the problem and helps us obtain the estimate of $\| \nabla^3{p} \|_{L^2}$.
 \end{itemize}

 The contents of this paper are as follows. In Section $\ref{sec-2}$, we summarize the spatial asymptotic property and monotonicity of the stationary solution which is given by \cite{Ji}. Then we give the main theorem of this paper. In Section $\ref{sec-3}$, we begin detailed discussion with a reformulation of the problem $\eqref{yfc}$ and $\eqref{yfctj}$ to that for the perturbation from the corresponding planar stationary solution. Then we prove the \emph{a priori} estimates of the perturbation. Moreover, we prove the large time behavior of the planar stationary solution by applying the Gagliardo-Nirenberg inequality. In Section $\ref{sec-5}$, using the time and space weighted energy method, we show the convergence rate toward the stationary solution under the assumption that the initial perturbation $v_0$ decays algebraically in $x$-direction with exponent $\alpha/2$. Using the interpolation inequality, we obtain the $L^\infty$ convergence rate $t^{-\alpha/2}$ of $(v,p)$ for non-degenerate case $(u_+<0)$. Finally, we show the time weighted estimate of derivatives in $y$-direction, which yields the convergence rate $t^{-1/4}$ in $L^\infty$ norm of $(v,p)$ for both of non-degenerate and degenerate cases. Combining the results of convergence rate, we get the $L^\infty$ decay rate $t^{-\alpha/2-1/4}$ of $(v,p)$ for non-degenerate case.

\section{Preliminaries and Main Theorems} \label{sec-2}
\subsection{Preliminaries}
Our main purpose in this paper is to show the asymptotic stability towards the planar stationary solution to the initial boundary value problem \eqref{yfc}-\eqref{yfcu0}, when $u_-<u_+\le0$. The corresponding planar stationary solution $(\bar{u},\bar{q})$ is defined as
\begin{equation}\label{eq-fhb}
\begin{cases}
 f(\bar{u})_x+\bar{q}_x=0,\\
  -\bar{q}_{xx}+\bar{q}+\bar{u}_x=0,\\
 \bar{u}(0)=u_-,\\
 \bar{u}(+\infty)=u_+, \ \ \ \ \bar{q}(+\infty)=0.
 \end{cases}
\end{equation}
The stationary solution $(\bar{u},\bar{q})$ is independent of $y$ and $t$. Ji, Zhang and Zhu in \cite{Ji} showed the existence and the precise asymptotic behavior of the stationary solution by applying the singular phase plane analysis method. Here, we summarize some properties of stationary solution $(\bar{u},\bar{q})$ in Proposition $\ref{prop-fhb}$. The proof of Proposition $\ref{prop-fhb}$ can be found in \cite{Ji}. For convenience of calculation, in the following, we assume that $f(u)=\frac{1}{2} u^2$ and the proof of general convex function $f(u)$ can be slightly modified.
\begin{prop}\label{prop-fhb}
	Suppose $u_-<u_+\le 0$ and let $\delta:=|u_--u_+|$. Then there exists a monotone increasing solution $\bar{u}(x)$ $(\bar{u}_x>0)$ to the stationary problem \eqref{eq-fhb}, such that the following estimates hold.
	\\[1mm] \indent
    $\mathrm{(i)}$ Non-degenerate case $\textbf{[ND]}$: Suppose that $u_+< 0$. Then, the solution $\bar{u}(x)$ satisfies
    $$\left|\partial_x^k\left\{\bar{u}(x)-u_+\right\}\right|\le C\delta  \mathrm{e}^{-\lambda x}, \ \ \ \ k=0,1,2,3,4,$$
    for some positive constants C and $\lambda$.
    \\[1mm] \indent
    $\mathrm{(ii)}$ Degenerate case $\textbf{[D]}$: Suppose that $u_+= 0$. Then, the solution $\bar{u}(x)$ satisfies
    $$|\partial_x^k\bar{u}(x)|  \le C \frac{\delta^{k+1}}{(1+\delta x)^{k+1}}, \ \ \ \ k=0,1,2,3,4,$$
    for some positive constants C.
\end{prop}

 Before we state our main theorem, we introduce some very useful inequalities in Lemmas $\ref{lem-prioriLinf}$-$\ref{lem-GN}$. The inequalities $\eqref{eq-Lyinf}$ and $\eqref{eq-Lxinf}$ are often used to deal with nonlinear higher derivatives in the \emph{a prior} estimate. The inequalities in Lemma $\ref{lem-decayLinf}$ and in Lemma $\ref{lem-GN}$ are used in calculating the convergence rate of $L^\infty$ norm and stating the large time behavior of planar stationary solution, respectively.
\begin{lem}\label{lem-prioriLinf} Let $f(x)\in H^1(\mathbbm{R})$, $h(x)\in H^1(\mathbbm{R}_+)$, it holds that
\begin{align}
	\| f \|_{L^\infty(\mathbbm{R})} \le \sqrt{2}\| f \|_{L^2(\mathbbm{R})}^\frac{1}{2}  \| f_x \|_{L^2(\mathbbm{R})}^\frac{1}{2},\label{eq-f-R}\\
	\| h \|_{L^\infty(\mathbbm{R}_+)} \le \sqrt{2}\| h \|_{L^2(\mathbbm{R}_+)}^\frac{1}{2} \| h_x \|_{L^2(\mathbbm{R}_+)}^\frac{1}{2}\label{eq-h-Rzheng}.
\end{align}
As an application of the inequalities $\eqref{eq-f-R}$, $\eqref{eq-h-Rzheng}$, we give two inequalities for the $L_x^2(\mathbbm{R}_+)L_y^\infty(\mathbbm{R})$ norm and $L_x^\infty(\mathbbm{R}_+)L_y^2(\mathbbm{R})$ norm, which will be often utilized later. For any function $f(x,y)\in H^1(\mathbbm{R}_+\times\mathbbm{R})$,
\begin{align}
	\| f \|_{L_x^2(L_y^\infty)}\le C \| f \|_{L^2}^\frac{1}{2} \| f_y \|_{L^2}^\frac{1}{2},\label{eq-Lyinf} \\
	\| f \|_{L_x^\infty(L_y^2)} \le C \| f \|_{L^2}^\frac{1}{2} \| f_x \|_{L^2}^\frac{1}{2}.\label{eq-Lxinf}
\end{align}
\end{lem}

\begin{lem}\label{lem-decayLinf} Suppose that $f(x,y)\in H^3(\mathbbm{R}_+\times\mathbbm{R})$. Using the Sobolev inequalities $\eqref{eq-h-Rzheng}$, $\eqref{eq-Lyinf}$, we have the following interpolation inequalities for the $L^\infty(\mathbbm{R}_+\times\mathbbm{R})$ norm:
\begin{align}
	&\|f\|_{L^\infty}\le C\|f\|_{L^2}^\frac{1}{4} \|f_x\|_{L^2}^\frac{1}{4}\|f_y\|_{L^2}^\frac{1}{4}\|f_{xy}\|_{L^2}^\frac{1}{4},\label{eq-v-Linfty} \\
	&\|\nabla{f}\|_{L^\infty}\le C\|\nabla{f}\|_{L^2}^\frac{1}{4} \|\nabla^2 f\|_{L^2}^\frac{1}{4}\| \nabla{f}_y\|_{L^2}^\frac{1}{4}\| \nabla^2 f_y\|_{L^2}^\frac{1}{4},\label{eq-nabv-Linfty}\\
	&\|f_y\|_{L^\infty}\le C\|f_y\|_{L^2}^\frac{1}{4} \|\nabla{f}_y\|_{L^2}^\frac{1}{4}\| f_{yy}\|_{L^2}^\frac{1}{4}\|\nabla{f}_{yy}\|_{L^2}^\frac{1}{4}.\label{eq-nabvy-inf}
\end{align}
\end{lem}

\begin{lem}(Gagliardo-Nirenberg's inequality)\label{lem-GN} Let $j$ and $m$ be integers with $0\le j<m$ and let $2\le p\le \infty$. Then any $f\in H^m(\mathbbm{R}^N)$ with $N\ge 1$ satisfies
\begin{equation}\label{eq-parj-f}
	\| \partial^j f \|_{L^p(\mathbbm{R}^N)} \le C \| f \|_{L^2(\mathbbm{R}^N)}^{1-\theta}\| \partial^m f \|_{L^2(\mathbbm{R}^N)}^\theta,
\end{equation}
where
$$\frac{1}{p}-\frac{j}{N}=(1-\theta)\frac{1}{2} +\theta\left(\frac{1}{2} -\frac{m}{N} \right)$$
for $\theta$ satisfing $j/m\le \theta<1$. Here $\partial^j$ denotes the totality of all the $j$-th order derivatives with respect to $x\in \mathbbm{R}^N$.
\end{lem}

\subsection{Main Theorems}
The main theorem of the present paper is stated as follows.
\begin{thm}\label{thm-main}
	Suppose that $u_-<u_+\le 0$ and that $\delta=|u_+-u_-|$ is sufficiently small. If $u_0- \bar{u}\in H^3(\mathbbm{R}_+\times \mathbbm{R})$, and there exists a sufficiently small constant $\delta'>0$ such that $\|u_0- \bar{u}\|_{H^3}^2+\delta\le \delta'$, then the initial boundary value problem \eqref{yfc}-\eqref{yfctj} has a unique solution $u\in C([0,\infty);H^3(\mathbbm{R}_+\times \mathbbm{R}))$, which satisfies the asymptotic behavior
    \begin{equation}\label{eq-asymptotic-behavior}
    \begin{aligned}
    &\sup \limits_{(x,y)\in\mathbbm{R}_+\times \mathbbm{R}}|\nabla^k (u(x,y,t)- \bar{u})|\rightarrow 0 \ \ \text{as} \ \ t \rightarrow \infty, \quad k=0,1,\\
    &\sup \limits_{(x,y)\in\mathbbm{R}_+\times \mathbbm{R}}|\nabla^k ({q}(x,y,t)-(\bar{q},0)^t)|\rightarrow 0 \ \ \text{as} \ \ t \rightarrow \infty, \quad k=0,1,\\
    &\sup \limits_{(x,y)\in\mathbbm{R}_+\times \mathbbm{R}}|\nabla ({\rm{div}}{q}(x,y,t)- \bar{q}_x)|\rightarrow 0 \ \ \text{as} \ \ t \rightarrow \infty.
    \end{aligned}	
    \end{equation}
    Here, $(\bar{q},0)^t$ represents the transpose of the vector $(\bar{q},0)$. 

    Moreover, for $\textbf{[ND]}$, if $u_0- \bar{u}\in L_{\alpha,2}^2(\mathbbm{R}_+\times \mathbbm{R})$ for $\alpha\ge0$, then the solution $u$ satisfies the decay estimates
	\begin{equation}\label{eq-convergence-rate-thm}
	\begin{aligned}
	 	&\sup \limits_{(x,y)\in\mathbbm{R}_+\times \mathbbm{R}}|u(x,y,t)- \bar{u}|\le C (\|u_0- \bar{u}\|_{H^3}+|u_0- \bar{u}|_{\alpha,2}) (1+t)^{-\frac{\alpha}{2} -\frac{1}{4}},\\
	 	&\sup \limits_{(x,y)\in\mathbbm{R}_+\times \mathbbm{R}}|(u(x,y,t)- \bar{u})_x|\le C (\|u_0- \bar{u}\|_{H^3}+|u_0- \bar{u}|_{\alpha,2}) (1+t)^{-\frac{\alpha}{2} -\frac{1}{4}},\\
	 	&\sup \limits_{(x,y)\in\mathbbm{R}_+\times \mathbbm{R}}|(u(x,y,t)- \bar{u})_y|\le C (\|u_0- \bar{u}\|_{H^3}+|u_0- \bar{u}|_{\alpha,2}) (1+t)^{-\frac{\alpha}{2} -\frac{3}{4}},		
	\end{aligned}
	 \end{equation}
	 and the solution $q$ satisfies the decay estimates
	 \begin{equation}\label{eq-q-sjsjgj}
	 	\begin{aligned}
	 	    &\sup \limits_{(x,y)\in\mathbbm{R}_+\times \mathbbm{R}}|\partial_x^k(q_1(x,y,t)- \bar{q}(x))|\le C (\|u_0- \bar{u}\|_{H^3}+|u_0- \bar{u}|_{\alpha,2}) (1+t)^{-\frac{\alpha}{2}-\frac{1}{4}}, \quad k=0,1,\\
	 	    &\sup \limits_{(x,y)\in\mathbbm{R}_+\times \mathbbm{R}}|(q_1(x,y,t)- \bar{q}(x))_y|\le C (\|u_0- \bar{u}\|_{H^3}+|u_0- \bar{u}|_{\alpha,2}) (1+t)^{-\frac{\alpha}{2}-\frac{1}{2}},\\
	 		&\sup \limits_{(x,y)\in\mathbbm{R}_+\times \mathbbm{R}}|\nabla^k q_2(x,y,t)|\le C (\|u_0- \bar{u}\|_{H^3}+|u_0- \bar{u}|_{\alpha,2}) (1+t)^{-\frac{\alpha}{2}-\frac{1}{2}}, \quad k=0,1,\\
	 		&\sup \limits_{(x,y)\in\mathbbm{R}_+\times \mathbbm{R}}|({\rm{div}}{q}(x,y,t)- \bar{q}_x(x))_x|\le C (\|u_0- \bar{u}\|_{H^3}+|u_0- \bar{u}|_{\alpha,2}) (1+t)^{-\frac{\alpha}{2}-\frac{1}{4}},\\
            &\sup \limits_{(x,y)\in\mathbbm{R}_+\times \mathbbm{R}}|({\rm{div}}{q}(x,y,t)- \bar{q}_x(x))_y|\le C (\|u_0- \bar{u}\|_{H^3}+|u_0- \bar{u}|_{\alpha,2}) (1+t)^{-\frac{\alpha}{2}-\frac{1}{2}},\\
	 	\end{aligned}
	 \end{equation}
	 for $t>0$ and some positive constant $C$, which is independent of $t$. For $\textbf{[D]}$, the solution $(u,q)$ satisfies the decay rate $\eqref{eq-convergence-rate-thm}$ and $\eqref{eq-q-sjsjgj}$ with $\alpha=0$.
\end{thm}
\textbf{Notations.} Throughout this paper, without any ambiguity, we denote a generic positive constant $C$ which may vary from line to line. For two functions $f$ and $h$, $f\sim h$ means
$$C^{-1}f(x)\le h(x)\le Cf(x).$$
 For any nonnegative constant $p$ $(1\le p\le \infty)$, $L^p=L^p(\mathbbm{R}_+\times\mathbbm{R})$ denotes usual Lebesgue space over $\mathbbm{R}_+\times\mathbbm{R}$, equipped with the norm $\|\cdot\|_{L^p}$. For any $l\ge0$, $H^l=H^l(\mathbbm{R}_+\times\mathbbm{R})$ denotes the usual Sobolev space over $\mathbbm{R}_+\times\mathbbm{R}$ with norm $\|\cdot\|_{H^l}$. We use the notation $\nabla^k f$ as in the meaning
$$\nabla^k f=(\partial_x^k f,\partial_x^{(k-1)}\partial_y f,...,\partial_x\partial_y^{(k-1)} f,\partial_y^k f)$$
where $f=f(x,y,t)$ and $\nabla^0 f=f$. And the notation $\Delta:=\partial_x^2+\partial_y^2$ denotes the Laplacian. For $\alpha\in \mathbbm{R}$, $L_\alpha^2=L_\alpha^2(\mathbbm{R}_+\times\mathbbm{R})$ denotes the space of measurable function $f$ satisfying $(1+x)^{\frac{\alpha}{2} }f\in L^2$ with the norm $|f|_\alpha:=(\int_{\mathbbm{R}}\int_{\mathbbm{R}_+} (1+x)^\alpha|f(x,y)|^2 \,\mathrm{d}x \mathrm{d}y )^{1/2} $. We also define $L_{\alpha,k}^2:=\{\nabla^j f\in L_\alpha^2; j=0,1,2,\cdots,k\}$ with the norm $|f|_{\alpha,k}:=(\sum_{j=0}^k|\nabla^j f|_\alpha^2)^{1/2}$. The Gaussian bracket $[x]$ denotes the greatest integer which does not exceed $x$.

\section{Asymptotic stability of planar stationary solution}\label{sec-3}
\subsection{\emph{A priori} estimates}
In this section, we derive the \emph{a priori} estimate in the $H^3$ Sobolev space. To do this, we introduce a perturbation $(v,p)$ by
\begin{equation*}
	\begin{aligned}
		&v(x,y,t)=u(x,y,t)-\bar{u}(x),\\
        &p(x,y,t)=q(x,y,t)-
          \left(                 
           \begin{array}{ccc}   
              \bar{q}(x)\\  
               0\\  
          \end{array}
          \right),
	\end{aligned}
\end{equation*}
and rewrite the problem \eqref{yfc}-\eqref{yfcu0} as
\begin{equation}\label{eq-vp-rd}
	\begin{cases}
		v_t+vv_x+(\bar{u}v)_x+g(v+\bar{u})_y+{\rm{div}}{p}=0,\\
		-\nabla{\rm{div}}{p}+p+\nabla{v}=0,
	\end{cases}
\end{equation}
with boundary condition and initial data
\begin{equation}\label{eq-vp-initial}
	v(0,y,t)=0, \ \ \ \ v(x,y,0)=v_0(x,y)=u_0(x,y)-\bar{u}(x).
\end{equation}
Notice that the equation $\eqref{eq-vp-rd}_2$ can be rewritten as
\begin{equation}\label{eq-rd2-dengjia}
	\begin{cases}
		-{\rm{div}}{p}_x+p_1+{v}_x=0,\\
		-{\rm{div}}{p}_y+p_2+{v}_y=0.
	\end{cases}
\end{equation}
Differentiating $\eqref{eq-rd2-dengjia}_1$ and $\eqref{eq-rd2-dengjia}_2$ with respect to $y$ and $x$, respectively. Then subtracting the two resulting equations, we can deduce that
\begin{equation}\label{eq-p1p2}
	p_{1y}=p_{2x},
\end{equation}
which will be frequently used in estimating the perturbation $p$.\\

The global existence follow from the combination of the local existence (Proposition $\ref{prop-local}$) and the a \emph{priori} estimates (Proposition $\ref{prop-priori}$). The local existence is proved by a standard iteration method with utilizing the heat kernel and its proof is omitted. In this section, we will denote ourselves to establish a \emph{priori} estimates under the \emph{a priori} assumption
\begin{equation}\label{eq-pri-assump}
	\|\nabla v(t)\|_{L^\infty}\le C \varepsilon_0,
\end{equation}
where $0<\varepsilon_0\ll 1$.

To illustrate the results of a \emph{priori} estimates and convergence rate, we put
\begin{equation}\label{eq-M0-Mal}
	M_0^2:=\|v_0\|_{H^3}^2, \ \ \ \ M_\alpha^2:=\|v_0\|_{H^3}^2+|v_0|_{\alpha,2}^2.
\end{equation}

For any $0<T\le \infty$, we seek the solution of the initial boundary value problem
$\eqref{eq-vp-rd}$-$\eqref{eq-vp-initial}$ in the set of functions $X(0,T)$ defined by
\begin{equation*}
X(0,T) = \left\{
\begin{tabular}{c|c}
   \multirow{2}{*}{({\it v, p})}    & $v\in C^0([0,T);H^3), \ \ \nabla{v}\in L^2(0,T;H^2)$   \\[1mm]
             & $p\in C^0([0,T);H^3)\cap L^2(0,T;H^3), \ \ {\rm{div}}{p}\in C^0([0,T);H^3)\cap L^2(0,T;H^3)$
\end{tabular}
       \right\}.
\end{equation*}
\begin{prop}[Local existence]\label{prop-local}
Suppose the boundary condition satisfies $u_-<u_+\le 0$ and the initial data satisfies $v_0 \in H^3(\mathbbm{R}_+\times\mathbbm{R})$. Also suppose that the initial data $M_0$ and $\delta$ are both small enough. Then there are two positive constants $C$ and $T_0$ such that
the problem $\eqref{eq-vp-rd}$-$\eqref{eq-vp-initial}$ has a unique solution $(v,p)\in X(0,T_0)$, which satisfies
\begin{equation*}
\|v(t)\|_{H^3}^2+\|p(t)\|_{H^3}^2+\| {\rm{div}}{p}(t) \|_{H^3}^2 +\int_0^t \  ( \| \nabla{v}(\tau) \|_{H^2}^2+\| p(\tau) \|_{H^3}^2+\| {\rm{div}}{p}(\tau) \|_{H^3}^2)  \,d{\tau} \le CM_0^2, \quad \forall t\in [0,T_0].
\end{equation*}
\end{prop}

\begin{prop}[A priori estimates] \label{prop-priori}
Let $T$ be a positive constant. Suppose that the problem $\eqref{eq-vp-rd}$-$\eqref{eq-vp-initial}$ has a unique solution $(v,p)\in X(0,T)$. Then there exist positive constants $C$ and $\delta_1$ such that if $M_0+\delta\le \delta_1$ $(0<\delta_1\ll 1)$, then we get the estimate
\begin{equation*}
\|v(t)\|_{H^3}^2+\|p(t)\|_{H^3}^2+\| {\rm{div}}{p}(t) \|_{H^3}^2+\int_0^t \ ( \| \nabla v(\tau)\|_{H^2}^2+\| p(\tau) \|_{H^3}^2+\| {\rm{div}}{p}(\tau) \|_{H^3}^2)  \,d{\tau}
\le C M_0^2, \quad \forall t\in [0,T].
\end{equation*}
\end{prop}

In order to prove the a \emph{priori} estimates, we firstly introduce Lemma $\ref{lem-DND}$ which plays an important role in estimating nonlinear terms. The proof of Lemma $\ref{lem-DND}$ is proved by the similar computation as in \cite{Nikkuni1999,Kawashima2003}, we omit the proof.
\begin{lem}\label{lem-DND}
For \textbf{[ND]} and \textbf{[D]}, there exists a positive constants $C$, such that it holds that
	\begin{align}
		&\int_{\mathbbm{R}}\int_{\mathbbm{R}_+} |\partial_x^j\bar{u}|^2|v|^2 \,\mathrm{d}x \mathrm{d}y \le C \delta\int_{\mathbbm{R}}\int_{\mathbbm{R}_+} |v_x|^2 \,\mathrm{d}x \mathrm{d}y, \ \ \ \ j=1,2,3,4, \label{eq-DND-v}\\
		&\int_{\mathbbm{R}}\int_{\mathbbm{R}_+} |\partial_x^j\bar{u}|^2|\partial_y^iv|^2 \,\mathrm{d}x \mathrm{d}y \le C \delta\int_{\mathbbm{R}}\int_{\mathbbm{R}_+} |\partial_y^iv_{x}|^2 \,\mathrm{d}x \mathrm{d}y,\ \ \ \ i=1,2.\label{eq-DND-vy}		
	\end{align}
For \textbf{[ND]}, it also holds that for any $\beta\in[0,\alpha]$, there exist positive constants $C$ and $C_\alpha$ satisfying
\begin{equation}\label{eq-ND-beta}
	\int_{\mathbbm{R}}\int_{\mathbbm{R}_+} (1+x)^\beta |\partial_x^j\bar{u}||v|^2 \,\mathrm{d}x \mathrm{d}y \le (C+C_\alpha \beta) \delta\int_{\mathbbm{R}}\int_{\mathbbm{R}_+} |v_x|^2 \,\mathrm{d}x \mathrm{d}y, \ \ \ \ j=1,2,3,4.
\end{equation}
\end{lem}

Lemma $\ref{lem-nab2v-dengjia}$ given below indicates that the $L^2$-estimates of $|\nabla^2 v|$ can convert to the $L^2$-estimates of $|\Delta v|+|\nabla v_y|$, and $|\nabla^3 v|$ can convert to $|\nabla\Delta v|+|\nabla v_{yy}|$ in $L^2$ norm. The former equivalence relation is used to calculate the $H^2$ convergence rate of $v$, while the latter equivalence relation is important in estimating the third derivatives of $v$ in the \emph{a priori} estimate.
\begin{lem}\label{lem-nab2v-dengjia}
	There exist positive constants $C_1$ and $C_2$ such that
	\begin{align}
		&C_1(|\Delta v|^2+|\nabla v_y|^2)\le |\nabla^2 v|^2\le C_2(|\Delta v|^2+|\nabla v_y|^2),\label{dengjia-nab2}\\
		&C_1(|\nabla\Delta v|^2+|\nabla v_{yy}|^2)\le |\nabla^3 v|^2\le C_2(|\nabla\Delta v|^2+|\nabla v_{yy}|^2)\label{dengjia-nab3}.
	\end{align}
	Namely, $|\nabla^2 v|$ is equivalent to $|\Delta v|+|\nabla v_y|$, and $|\nabla^3 v|$ is equivalent to $|\nabla\Delta v|+|\nabla v_{yy}|$.
\end{lem}

We divide the proof of the \emph{a priori} estimates into the following lemmas. Without the boundary information given to $p$, it is very difficult to obtain the estimates of $p$ and $v$ at the same time. So we try to get the estimate of $v$ first, and then get the estimate of $p$ through the existing estimates of $v$.
\begin{lem}\label{lem-v} Under the same assumptions as Proposition $\ref{prop-priori}$, then the following estimates are hold:
\begin{equation}\label{eq-v-jg}
\begin{aligned}
	\|v(t)\|_{H^1}^2+\int_0^t (\|\sqrt{\bar{u}_x}v(\tau)\|_{L^2}^2+\|\sqrt{\bar{u}_x}\nabla v(\tau)\|_{L^2}^2) \,\mathrm{d}\tau +\int_0^t \|\nabla v(\tau)\|_{L^2}^2 \,\mathrm{d}\tau +\int_0^t\int_{\mathbbm{R}} \left|\nabla v(0,y,t)\right|^2 \,\mathrm{d}y \mathrm{d}\tau \le C  M_0^2,\\
\end{aligned}	
\end{equation}
	for $t\in[0,T]$.
\end{lem}
{\it\bfseries Proof.}
We can obtain from $v\times\eqref{eq-vp-rd}_1+p\cdot \eqref{eq-vp-rd}_2$ that
\begin{equation}\label{eq-vnabv-djf1}
	\begin{aligned}[b]
		\left\{\frac{1}{2}v^2\right\}_t&+\frac{1}{2} \bar{u}_xv^2+{\rm{div}}\{pv\}-\nabla{\rm{div}}{p}\cdot p+|p|^2\\
		&+\left\{ \frac{1}{3}v^3+\frac{1}{2} \bar{u}v^2 \right\}_x+\left\{ g(v+\bar{u})v- \int_{\bar{u}}^{\bar{u}+v}g(s)  \,\mathrm{d}s \right\}_y=0.
	\end{aligned}
\end{equation}
After integration in $(x,y)\in \mathbbm{R}_+\times\mathbbm{R}$, the terms $-\nabla{\rm{div}}{p}\cdot p+|p|^2$ produce two good terms $\| {\rm{div}}{p}(t) \|_{L^2}^2 $ and $\| p(t) \|_{L^2}^2 $, but they also produce a boundary term $\int_{\mathbbm{R}} {\rm{div}}{p}(0,y,t)p_1(0,y,t) \, \mathrm{d}y $ that we can't deal with owing to no information for $p$ at the boundary $x=0$. The equation $\eqref{eq-vp-rd}_2$ gives
\begin{equation}\label{eq-v-gj}
\begin{aligned}[b]
	|\nabla{v}|^2&=|\nabla{\rm{div}}{p}|^2-2 \nabla{\rm{div}}{p}\cdot p+|p|^2\\
	             &=|\nabla{\rm{div}}{p}|^2+2({\rm{div}}{p})^2+|p|^2-2{\rm{div}}\{{\rm{div}}{p} p\},
\end{aligned}
\end{equation}
which indicates that we can replace the good term of $p$ with the good term of $v$ by the equation $\eqref{eq-v-gj}$. For the form of $\eqref{eq-v-gj}$, we need to calculate the higher order derivative of $v$. We can get from $\nabla{v}\cdot \nabla \eqref{eq-vp-rd}_1- \nabla{\rm{div}}{p}\cdot \eqref{eq-vp-rd}_2$ that
\begin{equation}\label{eq-nabvnab2v-djf1}
	\begin{aligned}[b]
		\left\{\frac{1}{2}|\nabla{v}|^2\right\}_t&+\frac{1}{2} \bar{u}_x|\nabla{v}|^2+\bar{u}_xv_x^2+|\nabla{\rm{div}}{p}|^2-\nabla{\rm{div}}{p}\cdot p\\
		&+\left\{ \frac{1}{2} v|\nabla{v}|^2+\frac{1}{2} \bar{u}|\nabla{v}|^2 \right\}_x
		+\left\{ \frac{1}{2} g'(v+\bar{u})|\nabla{v}|^2  \right\}_y \\
		&=-\frac{1}{2} v_x|\nabla{v}|^2- \bar{u}_{xx}vv_x- \frac{1}{2} g''(v+\bar{u})v_y|\nabla{v}|^2-g''(v+\bar{u})\bar{u}_xv_xv_y.
	\end{aligned}
\end{equation}
Adding $\eqref{eq-vnabv-djf1}$ and $\eqref{eq-nabvnab2v-djf1}$, we have
\begin{equation}\label{eq-vnabv}
	\begin{aligned}[b]
		 &\ \ \ \ \ \ \ \ \ \ \ \ \ \ \ \   \left\{\frac{1}{2}v^2+\frac{1}{2}|\nabla{v}|^2\right\}_t+\frac{1}{2} \bar{u}_x(v^2+|\nabla{v}|^2)+\bar{u}_xv_x^2+|\nabla{v}|^2+{\rm{div}}\{pv\}\\
		&+\left\{ \frac{1}{3}v^3+\frac{1}{2} \bar{u}v^2+\frac{1}{2} v|\nabla{v}|^2+\frac{1}{2} \bar{u}|\nabla{v}|^2 \right\}_x
		+\left\{ g(v+\bar{u})v- \int_{\bar{u}}^{\bar{u}+v}g(s)  \,\mathrm{d}s+\frac{1}{2} g'(v+\bar{u})|\nabla{v}|^2  \right\}_y \\
		&\ \ \ \ \ \ \ \ \ \ \ \ \ \ \ \ \ \ \ \ \ \ \ \  \ \ \ \ \ \ \ \ =-\frac{1}{2} v_x|\nabla{v}|^2- \bar{u}_{xx}vv_x- \frac{1}{2} g''(v+\bar{u})v_y|\nabla{v}|^2-g''(v+\bar{u})\bar{u}_xv_xv_y.
	\end{aligned}
\end{equation}

Note that the term $\{\cdots\}_y$ disappear after integration in $y\in \mathbbm{R}$, here and hereafer. For $v(0,y,t)=0$, the first three term in $\{\cdots\}_x$ disappear after integration in $x\in \mathbbm{R}_+$. And the last term in $\{\cdots\}_x$ is a good term after integration in $x\in \mathbbm{R}_+$ owing to $u_-<0$. Threrefore, integrating $\eqref{eq-vnabv}$ over $\mathbbm{R}_+\times \mathbbm{R}\times [0,t]$, we have
\begin{equation}\label{eq-v-jf}
	\begin{aligned}
		\|v(t)\|_{L^2}^2+&\|\nabla{v}(t)\|_{L^2}^2+\int_0^t (\|\sqrt{\bar{u}_x}v(\tau)\|_{L^2}^2+\|\sqrt{\bar{u}_x}\nabla v(\tau))\|_{L^2}^2) \,\mathrm{d}\tau +\int_0^t \|\nabla v(\tau)\|_{L^2}^2 \,\mathrm{d}\tau +\int_0^t\int_{\mathbbm{R}} \left|\nabla v(0,y,t)\right|^2 \,\mathrm{d}y \mathrm{d}\tau \\
		&\le CM_0^2+C \int_0^t\int_{\mathbbm{R}}\int_{\mathbbm{R}_+} (|v_x||\nabla{v}|^2+|\bar{u}_{xx}vv_x|+|v_y||\nabla{v}|^2+\bar{u}_x|v_xv_y|) \,\mathrm{d}x \mathrm{d}y \mathrm{d}\tau.
	\end{aligned}
\end{equation}
By Proposition $\ref{prop-fhb}$ and Lemma $\ref{lem-DND}$, using the Cauchy inequality, under the a \emph{priori} assumption $\eqref{eq-pri-assump}$, the right-hand side of $\eqref{eq-v-jf}$ can be estimated as follows:
\begin{equation}\label{eq-v-gj1}
	\int_0^t\int_{\mathbbm{R}}\int_{\mathbbm{R}_+} (|v_x||\nabla{v}|^2+|v_y||\nabla{v}|^2) \,\mathrm{d}x \mathrm{d}y \mathrm{d}\tau
	\le \varepsilon_0 \int_0^t \|\nabla{v}(\tau)\|_{L^2}^2 \,\mathrm{d}\tau,
\end{equation}
\begin{equation}\label{3.19}
\begin{aligned}[b]
	\int_0^t\int_{\mathbbm{R}}\int_{\mathbbm{R}_+} |\bar{u}_{xx}vv_x| \,\mathrm{d}x \mathrm{d}y \mathrm{d}\tau
	&\le \frac{1}{4} \int_0^t \|\nabla{v}(\tau)\|_{L^2}^2 \,\mathrm{d}\tau +C \int_0^t\int_{\mathbbm{R}}\int_{\mathbbm{R}_+} |\bar{u}_{xx}|^2|v|^2 \,\mathrm{d}x \mathrm{d}y \mathrm{d}\tau  \\
	&\le \frac{1}{4} \int_0^t \|\nabla{v}(\tau)\|_{L^2}^2 \,\mathrm{d}\tau +C \delta\int_0^t \| v_x(\tau))\|_{L^2}^2 \,\mathrm{d}\tau,	
\end{aligned}
\end{equation}
and
\begin{equation}\label{eq-v-gj3}
	\int_0^t\int_{\mathbbm{R}}\int_{\mathbbm{R}_+} \bar{u}_x|v_xv_y| \,\mathrm{d}x \mathrm{d}y \mathrm{d}\tau
	\le C \| \bar{u}_x \|_{L^\infty} \int_0^t \|\nabla{v}(\tau)\|_{L^2}^2 \,\mathrm{d}\tau
	\le C \delta \int_0^t \|\nabla{v}(\tau)\|_{L^2}^2 \,\mathrm{d}\tau.
\end{equation}
Substituting $\eqref{eq-v-gj1}$-$\eqref{eq-v-gj3}$ into $\eqref{eq-v-jf}$, for some small $\varepsilon_0$ and $\delta$, we get $\eqref{eq-v-jg}$.
$\hfill\Box$

\begin{lem}\label{lem-dlnabv}Under the same assumptions as Proposition $\ref{prop-priori}$, we get the higher order derivative estimation:
\begin{equation}\label{eq-dlnabv}
\begin{aligned}[b]
	\|\nabla v(t)\|_{L^2}^2+\|\nabla^2 v(t)\|_{L^2}^2+\int_0^t (\|\sqrt{\bar{u}_x}\nabla v(\tau)\|_{L^2}^2+ \|\sqrt{\bar{u}_x}\Delta v(\tau)\|_{L^2}^2) \,\mathrm{d}\tau+\int_0^t \|\nabla^2 v(\tau)\|_{L^2}^2 \,\mathrm{d}\tau& \\	
	+\int_0^t\int_{\mathbbm{R}} (\Delta v(0,y,t))^2 \,\mathrm{d}y \mathrm{d}\tau &\le C M_0^2,	
\end{aligned}
\end{equation}
 	for $t\in[0,T]$.
\end{lem}
{\it\bfseries Proof.}
We can get from $\nabla{v}\cdot \nabla \eqref{eq-vp-rd}_1+ {\rm{div}}{p}\times {\rm{div}}\eqref{eq-vp-rd}_2+ \Delta v\times \Delta \eqref{eq-vp-rd}_1- \Delta {\rm{div}}{p}\times {\rm{div}}\eqref{eq-vp-rd}_2$ that
\begin{equation}\label{eq-nabdelv-wf}
\begin{aligned}[b]
	& \ \ \ \ \ \  \left\{\frac{1}{2} |\nabla{v}|^2+\frac{1}{2} (\Delta v)^2\right\}_t+\frac{1}{2} \bar{u}_x(|\nabla{v}|^2+ (\Delta v)^2)+\bar{u}_x v_x^2+|\nabla^2 v|^2+{\rm{div}}\{{\rm{div}}{p} \nabla{v}\}\\
	&+\left\{\frac{1}{2} v|\nabla{v}|^2\!+\!\frac{1}{2} \bar{u}|\nabla{v}|^2\!+\!\frac{1}{2} v(\Delta v)^2\!+\!\frac{1}{2} \bar{u}(\Delta v)^2\!-\!2v_{xy}v_y\right\}_x\!+\!\left\{\frac{1}{2} g'(v\!+\!\bar{u})|\nabla{v}|^2\!+\!2v_{xx}v_y\right\}_y\\
	&=-\frac{1}{2} v_x|\nabla{v}|^2 - \bar{u}_{xx}vv_x- \frac{1}{2} g''(v+\bar{u})v_y|\nabla{v}|^2-g''(v+\bar{u})\bar{u}_xv_xv_y-\frac{1}{2} v_x(\Delta v)^2-2 \nabla v\cdot \nabla v_x \Delta v \\
	& \ \ \ \ - \bar{u}_{xxx}v  \Delta v-3 \bar{u}_{xx}v_x \Delta v-2 \bar{u}_xv_{xx}\Delta v-\Delta (g(v+\bar{u}))\Delta v,
\end{aligned}
\end{equation}
where
\begin{equation}\label{eq-delg-delv}
\begin{aligned}[b]
\Delta (g(v+\bar{u}))\Delta v
=&\frac{1}{2} g''(v+\bar{u})v_y(\Delta v)^2 +g^{(3)}(v+\bar{u})(v_x+\bar{u}_x)^2v_y \Delta v+g''(v+\bar{u})\bar{u}_{xx}v_y \Delta v\\ &+g^{(3)}(v+\bar{u})v_y^3 \Delta v+2g''(v+\bar{u})(\nabla{v}\cdot \nabla{v}_y+\bar{u}_xv_{xy})\Delta v+\left\{\frac{1}{2} g'(v\!+\!\bar{u})(\Delta v)^2\right\}_y.
\end{aligned}
\end{equation}
Here we have used the fact from $\eqref{eq-vp-rd}_2$ that
\begin{align}
	 &(\Delta v)^2= (\Delta {\rm{div}}{p})^2+ ({\rm{div}}{p})^2-2(\Delta {\rm{div}}{p}){\rm{div}}{p},\\
	 &(\Delta v)^2=|\nabla^2 v|^2+2\{v_{xx}v_y\}_y-2\{v_{xy}v_y\}_x,\label{eq-delv-nab2v}\\
	 &\Delta (g'(v+\bar{u}))=g^{(3)}(v+\bar{u})(v_x+\bar{u}_x)^2+g''(v+\bar{u})\bar{u}_{xx}+g^{(3)}(v+\bar{u})v_y^2+g''(v+\bar{u})\Delta v,\label{eq-del-g}
\end{align}
and for any functions $f(x,y)$ and $h(x,y)$,
\begin{equation}\label{eq-del-fh}
	\Delta (fh)=(\Delta f) h+f (\Delta h)+2 \nabla{f}\cdot \nabla h.\\
\end{equation}
Notice that the first four terms on the right-hand side of $\eqref{eq-nabdelv-wf}$ is bounded by $CM_0^2$ after integration in $(x,y,t)\in\mathbbm{R}_+\times\mathbbm{R}\times[0,t]$ owing to $\eqref{eq-v-gj1}$-$\eqref{eq-v-gj3}$. From $\eqref{eq-delv-nab2v}$ and $v_y(0,y,t)=0$, we can conclude that
$$\int_{\mathbbm{R}}\int_{\mathbbm{R}_+} (\Delta v)^2 \,\mathrm{d}x \mathrm{d}y=\int_{\mathbbm{R}}\int_{\mathbbm{R}_+} |\nabla^2 v|^2 \,\mathrm{d}x  \mathrm{d}y.$$
Thus, integrating $\eqref{eq-nabdelv-wf}$ over $\mathbbm{R}_+\times\mathbbm{R}\times[0,t]$, we have
\begin{equation}\label{eq-dlnabv-jf}
	\begin{aligned}[b]
		&\| \nabla{v}(t) \|_{L^2}^2+ \|\nabla^2 v(t)\|_{L^2}^2+\int_0^t (\|\sqrt{\bar{u}_x} \nabla{v}(\tau)\|_{L^2}^2+\|\sqrt{\bar{u}_x}\Delta v(\tau)\|_{L^2}^2) \,\mathrm{d}\tau +\int_0^t \|\nabla^2 v(\tau)\|_{L^2}^2 \,\mathrm{d}\tau \\
		&\ \ \ \ \ \ \ \ \ \ \ \  \ \ \ \ \ \ \ \ \ \ \ 	+\int_0^t\int_{\mathbbm{R}} (\Delta v(0,y,t))^2 \,\mathrm{d}y \mathrm{d}\tau \\
		\le& C M_0^2\!+\!\int_0^t\int_{\mathbbm{R}} {\rm{div}}{p}(0,y,t)v_x(0,y,t) \,\mathrm{d}y \mathrm{d}\tau\!+\!C \int_0^t\int_{\mathbbm{R}}\int_{\mathbbm{R}_+} (|\nabla{v}|(\Delta v)^2\!+\!|\bar{u}_{xxx}||v||\Delta v|\!+\!\bar{u}_x|v_{xx}||\Delta v|) \,\mathrm{d}x \mathrm{d}y \mathrm{d}\tau\\
		&+C \int_0^t\int_{\mathbbm{R}}\int_{\mathbbm{R}_+} (|\bar{u}_{xx}||\nabla{v}||\Delta v|+|\nabla{v}|^3|\Delta v|+\bar{u}_x^2|\nabla{v}||\Delta v|+|\nabla{v}||\nabla{v}_x|| \Delta v|+(|\nabla{v}|+\bar{u}_x)|\nabla{v}_y|| \Delta v|) \,\mathrm{d}x \mathrm{d}y \mathrm{d}\tau.
	\end{aligned}
\end{equation}
Since $\eqref{eq-vp-rd}_1$ gives
\begin{equation}
	{\rm{div}}{p}(0,y,t)=-u_-v_x(0,y,t),
\end{equation}\label{eq-vx0-bj}
it follows that the second term on right-hand side of $\eqref{eq-dlnabv-jf}$ can be estimated as
\begin{equation}\label{eq-dlnabv-bj}
	\int_0^t\int_{\mathbbm{R}} {\rm{div}}{p}(0,y,t)v_x(0,y,t) \,\mathrm{d}y \mathrm{d}\tau\le C \int_{\mathbbm{R}}\int_{\mathbbm{R}_+} |\nabla{v}(0,y,t)|^2 \,\mathrm{d}y \mathrm{d}y.
\end{equation}
Under the a \emph{priori} assumption $\eqref{eq-pri-assump}$, we get
\begin{equation}\label{eq-dlnabv-1}
	\int_0^t\int_{\mathbbm{R}}\int_{\mathbbm{R}_+} |\nabla{v}|(\Delta v)^2 \,\mathrm{d}x \mathrm{d}y \mathrm{d}\tau
	\le C \varepsilon_0 \int_0^t \|\nabla^2 v(\tau)\|_{L^2}^2 \,\mathrm{d}\tau.
\end{equation}
By Proposition $\ref{prop-fhb}$ and Lemma $\ref{lem-DND}$, using Cauchy inequality, the remaining term on the right-hand side of $\eqref{eq-dlnabv-jf}$ is bounded by
\begin{equation}
\begin{aligned}[b]
	\int_0^t\int_{\mathbbm{R}}\int_{\mathbbm{R}_+} |\bar{u}_{xxx}||v||\Delta v| \,\mathrm{d}x \mathrm{d}y \mathrm{d}\tau
	&\le \frac{1}{4} \int_0^t \|\nabla^2 v(\tau)\|_{L^2}^2 \,\mathrm{d}\tau +C \int_0^t\int_{\mathbbm{R}}\int_{\mathbbm{R}_+} |\bar{u}_{xxx}|^2|v|^2 \,\mathrm{d}x \mathrm{d}y \mathrm{d}\tau\\
	&\le \frac{1}{4} \int_0^t \|\nabla^2 v(\tau)\|_{L^2}^2 \,\mathrm{d}\tau + C \int_0^t\int_{\mathbbm{R}}\int_{\mathbbm{R}_+} |\nabla{v}(\tau)|^2 \,\mathrm{d}x \mathrm{d}y \mathrm{d}\tau, 	
\end{aligned}
\end{equation}
\begin{equation}\label{eq-dlnabv-3}
\begin{aligned}[b]
	\int_0^t\int_{\mathbbm{R}}\int_{\mathbbm{R}_+} &(|\bar{u}_{xx}| + |\nabla{v}|^2 + \bar{u}_x^2)|\nabla{v}||\Delta v| \,\mathrm{d}x \mathrm{d}y \mathrm{d}\tau\\
	&\le \frac{1}{4} \int_0^t \|\nabla^2 v(\tau)\|^2 \,\mathrm{d}\tau  + \sup\limits_{0\le \tau \le t }\{\|\bar{u}_{xx}\|_{L^\infty}^2 + \|\nabla{v}\|_{L^\infty}^2 + \|\bar{u}_x\|_{L^\infty}^2\}\int_0^t \|\nabla{v}(\tau)\|^2 \,\mathrm{d}\tau,	
\end{aligned}
\end{equation}
and
\begin{equation}\label{eq-dlnabv-4}
	\begin{aligned}
		 \int_0^t\int_{\mathbbm{R}}\int_{\mathbbm{R}_+} ((\bar{u}_x|v_{xx}|+|\nabla{v}||\nabla{v}_x|)| \Delta v|+(|\nabla{v}|+\bar{u}_x)|\nabla{v}_y|| \Delta v|) \,\mathrm{d}x \mathrm{d}y \mathrm{d}\tau
		 \le C (\varepsilon_0+\delta) \int_0^t \| \nabla^2 v(\tau)\|_{L^2}^2  \,\mathrm{d}\tau.
	\end{aligned}
\end{equation}
Substituting $\eqref{eq-dlnabv-bj}$-$\eqref{eq-dlnabv-4}$ into $\eqref{eq-dlnabv-jf}$, combining the results of Lemma $\ref{lem-v}$, for some small $\varepsilon_0$ and $\delta$, we get $\eqref{eq-dlnabv}$.
 $\hfill\Box $
\\

 Based on the $H^2$-estimate of $v$, we can get the corresponding estimates of $p$ from the elliptic equation $\eqref{eq-vp-rd}_2$.
\begin{lem}\label{lem-p-1-4}Under the same assumptions as Proposition $\ref{prop-priori}$, the perturbation $p$ can be estimated as follows:
\begin{equation}\label{eq-p-h2-gj}
	\int_0^t \left(\|\nabla{\rm{div}}p(\tau)\|_{L^2}^2+\|{\rm{div}}p(\tau)\|_{L^2}^2+\|p(\tau)\|_{L^2}^2\right) \, \mathrm{d}\tau \le C M_0^2,
\end{equation}
	for $t\in[0,T]$.
\end{lem}
 {\it\bfseries Proof.}
By $\eqref{eq-v-gj}$, we get further
\begin{equation}\label{eq-p-cjfgj}
\begin{aligned}[b]
	&\int_0^t\int_{\mathbbm{R}}\int_{\mathbbm{R}_+} (|\nabla{\rm{div}}{p}|^2+2({\rm{div}}{p})^2+|p|^2) \,\mathrm{d}x \mathrm{d}y \mathrm{d}\tau \\
	=&\int_0^t \|\nabla{v}(\tau)\|_{L^2}^2 \,\mathrm{d}\tau +2 \int_0^t\int_{\mathbbm{R}}\int_{\mathbbm{R}_+} {\rm{div}}\{({\rm{div}}{p})p\} \,\mathrm{d}x \mathrm{d}y \mathrm{d}\tau \\
	=&\int_0^t \|\nabla{v}(\tau)\|_{L^2}^2 \,\mathrm{d}\tau -2 \int_0^t\int_{\mathbbm{R}}\int_{\mathbbm{R}_+} {\rm{div}}{p}(0,y,\tau)p_1(0,y,\tau) \,\mathrm{d}x \mathrm{d}y \mathrm{d}\tau.
\end{aligned}
\end{equation}
In order to get the estimate of the last term of $\eqref{eq-p-cjfgj}$, we may convert the boundary terms ${\rm{div}}{p}(0,y,\tau)$ and $p_1(0,y,\tau)$ to some boundary terms of $v$, because no information of $p$ is given at the boundary. For $v(0,y,t)=v_y(0,y,t)=0$, the equation $\eqref{eq-vp-rd}$ gives
\begin{align}
    &p_1(0,y,t)={\rm{div}}{p}_x(0,y,t)-v_x(0,y,t),\\
	&{\rm{div}}{p}_x(0,y,t)=-v_{xt}(0,y,t)-v_x^2(0,y,t)-2 \bar{u}_x(0)v_x(0,y,t)-u_-v_{xx}(0,y,t)-g'(u_-)v_{xy}(0,y,t),\label{eq-divpx0}
\end{align}
Recalling $\eqref{eq-vx0-bj}$, it follows from $\eqref{eq-p-cjfgj}$ that
\begin{equation}\label{eq-divpp-cjf}
	\begin{aligned}[b]
	&\int_0^t\int_{\mathbbm{R}}\int_{\mathbbm{R}_+} (|\nabla{\rm{div}}{p}|^2+2({\rm{div}}{p})^2+|p|^2) \,\mathrm{d}x \mathrm{d}y \mathrm{d}\tau \\
	=& \int_0^t \|\nabla{v}(\tau)\|_{L^2}^2 \,\mathrm{d}\tau +2u_- \int_0^t\int_{\mathbbm{R}} v_x(0,y,\tau)p_1(0,y,\tau) \,\mathrm{d}y \mathrm{d}\tau \\
    =& \int_0^t \|\nabla{v}(\tau)\|_{L^2}^2 \,\mathrm{d}\tau +2u_- \int_0^t\int_{\mathbbm{R}} v_x(0,y,\tau)({\rm{div}}{p}_x(0,y,\tau)-v_x(0,y,\tau)) \,\mathrm{d}y \mathrm{d}\tau\\
    =& \int_0^t \|\nabla{v}(\tau)\|_{L^2}^2 \,\mathrm{d}\tau-2u_-\int_0^t\int_{\mathbbm{R}} |\nabla{v}(0,y,\tau)|^2 \,\mathrm{d}y \mathrm{d}\tau +2u_- \int_0^t\int_{\mathbbm{R}} v_x(0,y,\tau){\rm{div}}{p}_x(0,y,\tau) \,\mathrm{d}y \mathrm{d}\tau.
	\end{aligned}
\end{equation}
 By employing $\eqref{eq-divpx0}$, the last term of $\eqref{eq-divpp-cjf}$ can be divided into the following estimates $\eqref{eq-divpx0-1}$-$\eqref{eq-divpx0-4}$:
\begin{equation}\label{eq-divpx0-1}
\begin{aligned}[b]
	\int_0^t\int_{\mathbbm{R}} v_x(0,y,t)v_{xt}(0,y,t) \,\mathrm{d}y \mathrm{d}\tau
	&=\int_0^t\int_{\mathbbm{R}} \{v_x^2(0,y,t)\}_t \,\mathrm{d}y \mathrm{d}\tau\\
	&=\int_{\mathbbm{R}} v_x^2(0,y,t) \,\mathrm{d}y -\int_{\mathbbm{R}} v_x^2(0,y,0) \,\mathrm{d}y\\
	&\le \int_{\mathbbm{R}} \|v_x(\cdot,y,t)\|_{L_x^\infty}^2 \,\mathrm{d}y -\int_{\mathbbm{R}} \|v_x(\cdot,y,0)\|_{L_x^\infty}^2 \,\mathrm{d}y\\
	&\le \|v_x(t)\|_{L^2}^2+\|v_{xx}(t)\|_{L^2}^2+ \| v_0 \|_{H^2}^2,	
\end{aligned}
\end{equation}
\begin{equation}\label{eq-divpx0-2}
	\begin{aligned}
		\int_0^t\int_{\mathbbm{R}} v_x(0,y,t)v_{xx}(0,y,t) \,\mathrm{d}y \mathrm{d}\tau \le \int_0^t\int_{\mathbbm{R}} |\nabla v(0,y,t)|^2 \,\mathrm{d}y \mathrm{d}\tau+\int_0^t\int_{\mathbbm{R}} (\Delta v(0,y,t))^2 \,\mathrm{d}y \mathrm{d}\tau,
	\end{aligned}
\end{equation}
and
\begin{equation}\label{eq-divpx0-3}
	\int_0^t\int_{\mathbbm{R}} v_x(0,y,t)v_{xy}(0,y,t) \,\mathrm{d}y \mathrm{d}\tau=0,
\end{equation}
by using integration by parts.
Since $\sup\limits_{0\le \tau \le t}\| \nabla{v}(\tau) \|_{L^\infty}$ and $\| \bar{u}_x\|_{L^\infty}$ are bounded, we get
\begin{equation}\label{eq-divpx0-4}
	\begin{aligned}
		\int_0^t\int_{\mathbbm{R}} v_x^3(0,y,t)+\bar{u}_x(0,y,t)v_x^2(0,y,t) \,\mathrm{d}y \mathrm{d}\tau
		\le C \int_0^t\int_{\mathbbm{R}} v_x^2(0,y,t) \,\mathrm{d}y \mathrm{d}\tau.
	\end{aligned}
\end{equation}
Thus, the last term of $\eqref{eq-divpp-cjf}$ can be estimated as
\begin{equation}\label{eq-vxdivp-0}
	\begin{aligned}
		\int_0^t\int_{\mathbbm{R}} v_x(0,y,t){\rm{div}}{p}_x(0,y,t) \,\mathrm{d}y \mathrm{d}\tau
		\le C (\| v(t) \|_{H^2}^2+  \| v_0 \|_{H^2}^2)+C \int_0^t\int_{\mathbbm{R}} (|\nabla{v}(0,y,t)|^2+|\Delta{v}(0,y,t)|^2) \,\mathrm{d}y \mathrm{d}\tau.
	\end{aligned}
\end{equation}
 Substituting $\eqref{eq-vxdivp-0}$ into $\eqref{eq-divpp-cjf}$, combining $\eqref{eq-v-jg}$ and $\eqref{eq-dlnabv}$, we get $\eqref{eq-p-h2-gj}$.
  $\hfill\Box$
\\

Next, we consider the higher order estimates of $v$ and $p$. We still need to get the third-order estimation of $v$ first, and then get the higher-order estimation of $p$ through the estimation of $v$. To do this, we may use the structure of $\eqref{eq-vp-rd}$ to convert $p$ to $v$, just like we did before. But in this way, we can only get the $L^2$-estimates of $\nabla\Delta {v}$ instead of $\nabla^3{v}$, see $\eqref{eq-nav3-3}$. From $\eqref{dengjia-nab3}$, we note that $\| \nabla^3{v}(t)\|_{L^2} $ is equivalent to $\| \nabla\Delta {v}(t) \|_{L^2} +\| \nabla {v}_{yy}(t) \|_{L^2} $. Thanks to $v_{yy}(0,y,t)= 0$, the boundary terms generated in estimating $\| {v}_{yy}(t) \|_{H^1}$ is well handled, which promotes us to obtain $H^3$-estimates of $v$. Thus, we firstly give the estimates of $\| {v}_{yy}(t) \|_{L^2} $. Secondly, we give the estimates of $\| \nabla^2{\rm{div}}{p}(t) \|_{L^2} $ which will appear in the estimation of $\|\nabla\Delta {v} (t)\|_{L^2} $. In the end, we show the estimate of $\|\nabla\Delta {v} (t)\|_{L^2}  $, which completes the proof of third-order estimation of $v$.

 \begin{lem}\label{lem-vyy}
 	Under the same assumptions as Proposition $\ref{prop-priori}$, it holds that
 	\begin{equation}\label{eq-vyy}
	\begin{aligned}[b]
		&\| v_{yy}(t) \|_{H^1}^2+\int_0^t (\|\sqrt{\bar{u}_x}v_{yy}(\tau)\|_{L^2}^2+\|\sqrt{\bar{u}_x} \nabla{v}_{yy}(\tau)\|_{L^2}^2) \,\mathrm{d}\tau +\int_0^t \|\nabla{v}_{yy}(\tau)\|_{L^2}^2 \,\mathrm{d}\tau+\int_0^t \int_{\mathbbm{R}} |\nabla{v}_{yy}(0,y,\tau)|^2 \,\mathrm{d}y\mathrm{d}\tau\\
		\le& CM_0^2+C (\varepsilon_0+M_0^4)\int_0^t \| \nabla \Delta{v}(\tau) \|_{L^2}^2 \,\mathrm{d}\tau,
	\end{aligned}
\end{equation}
for $t\in[0,T]$.
 \end{lem}
  {\it\bfseries Proof.}
We can get from $v_{yy}\times\partial_{yy}\eqref{eq-vp-rd}_1+p_{yy}\cdot \partial_{yy}\eqref{eq-vp-rd}_2+\nabla{v}_{yy}\cdot \nabla \partial_{yy}\eqref{eq-vp-rd}_1- \nabla{\rm{div}}{p}_{yy}\cdot \partial_{yy} \eqref{eq-vp-rd}_2$ that
\begin{equation}\label{eq-vyy-nabvyy}
	\begin{aligned}[b]
		&\left\{\frac{1}{2}v_{yy}^2+\frac{1}{2}|\nabla{v}_{yy}|^2\right\}_t+\frac{1}{2} \bar{u}_x(v_{yy}^2+|\nabla{v}_{yy}|^2)+\bar{u}_xv_{xyy}^2+|\nabla{v}_{yy}|^2+{\rm{div}}(p_{yy}v_{yy} )\\
		&\ \ \ \ +\left\{ \frac{1}{2}vv_{yy}^2+\frac{1}{2} \bar{u}v_{yy}^2+\frac{1}{2} v|\nabla{v}_{yy}|^2+\frac{1}{2} \bar{u}|\nabla{v}_{yy}|^2 \right\}_x
		+\left\{\frac{1}{2}  g'(v+\bar{u})v_{yy}^2+\frac{1}{2} g'(v+\bar{u})|\nabla{v}_{yy}|^2  \right\}_y \\
		& =-\frac{1}{2} v_xv_{yy}^2-2v_yv_{xy}v_{yy}- g'''(v+\bar{u})v_y^3v_{yy}-\frac{5}{2}g''(v+\bar{u})v_yv_{yy}^2- \frac{1}{2}  v_x|\nabla{v}_{yy}|^2-v_{yy}\nabla{v}_{yy}\cdot \nabla{v}_x\\
		& \ \ \ -2v_y \nabla{v}_{xy}\cdot\nabla{v}_{yy}-2v_{xy}\nabla{v}_y\cdot \nabla{v}_{yy}-v_{xyy}\nabla{v}\cdot \nabla{v}_{yy}-\bar{u}_{xx}v_{yy}v_{xyy}-g^{(4)}(v+\bar{u})(\nabla{v}\cdot\nabla{v}_{yy}+\bar{u}_xv_{xyy})v_y^3\\
		& \ \ \ -3g'''(v+\bar{u})v_y^2 \nabla{v}_y\cdot \nabla{v}_{yy}-3g'''(v+\bar{u})(\nabla{v}\cdot \nabla{v}_{yy}+\bar{u}_xv_{xyy})v_yv_{yy}\\
		& \ \ \ -3g''(v+\bar{u})\nabla{v}_y\cdot \nabla{v}_{yy}v_{yy} -\frac{5}{2}g''(v+\bar{u})v_y|\nabla{v}_{yy}|^2-g''(v+\bar{u})(\nabla{v}\cdot \nabla{v}_{yy}+\bar{u}_xv_{xyy})v_{yyy}.
	\end{aligned}
\end{equation}
Integrating $\eqref{eq-vyy-nabvyy}$ over $\mathbbm{R}_+\times\mathbbm{R}$, we obtain
\begin{equation}\label{eq-vyyH1}
	\begin{aligned}[b]
		&\frac{\mathrm{d}}{\mathrm{d}t}(\| v_{yy} \|_{L^2}^2 +\| \nabla{v}_{yy} \|_{L^2}^2 )+\|\sqrt{\bar{u}_x} v_{yy} \|_{L^2}^2 +\|\sqrt{\bar{u}_x}  \nabla{v}_{yy} \|_{L^2}^2+\| \nabla{v}_{yy} \|_{L^2}^2+\int_{\mathbbm{R}} |\nabla{v}_{yy}(0,y,\tau)|^2 \,\mathrm{d}y \\
		\le& C \int_{\mathbbm{R}}\int_{\mathbbm{R}_+} (|\nabla{v}||v_{yy}|^2+|v_y||\nabla{v}_y||v_{yy}|+|v_y|^3|v_{yy}|+|\nabla{v}||\nabla{v}_{yy}|^2+|v_{yy}||\nabla v_x||\nabla{v}_{yy}|+|v_y||\nabla v_{xy}||\nabla{v}_{yy}|) \,\mathrm{d}x \mathrm{d}y \\
		&  +C \int_{\mathbbm{R}}\int_{\mathbbm{R}_+} (|\nabla{v}_y|^2|\nabla{v}_{yy}|+|\bar{u}_{xx}||v_{yy}||\nabla{v}_{yy}|+(|\nabla{v}|+\bar{u}_x)|v_y|^3|\nabla{v}_{yy}|+|v_y|^2|\nabla{v}_y||\nabla{v}_{yy}|)\,\mathrm{d}x \mathrm{d}y\\
		& +C \int_{\mathbbm{R}}\int_{\mathbbm{R}_+} (|\nabla{v}|+\bar{u}_x)|v_y||v_{yy}||\nabla{v}_{yy}|+|\nabla{v}_y||v_{yy}||\nabla{v}_{yy}|+|v_y||\nabla{v}_{yy}|^2+(|\nabla{v}|+\bar{u}_x)|v_{yyy}||\nabla{v}_{yy}|) \,\mathrm{d}x \mathrm{d}y.
	\end{aligned}
\end{equation}
Note that we have obtained the first and second derivative estimates of $v$ in $\eqref{eq-v-jg}$ and $\eqref{eq-dlnabv}$. Thus, using Cauchy inequality, the first three terms on the right-hand side of $\eqref{eq-vyyH1}$ are bounded by
\begin{equation*}
	\int_{\mathbbm{R}}\int_{\mathbbm{R}_+} (|\nabla{v}||v_{yy}|^2+|v_y||\nabla{v}_y||v_{yy}|+|v_y|^3|v_{yy}|)  \,\mathrm{d}x \mathrm{d}y
	\le C(\| \nabla{v} \|_{L^2}^2+\| \nabla^2{v} \|_{L^2}^2).
\end{equation*}
Using $\eqref{dengjia-nab3}$, $\eqref{eq-Lyinf}$, $\eqref{eq-Lxinf}$ and Cauchy inequality, the fifth, sixth and seventh terms on the right-hand side of $\eqref{eq-vyyH1}$ can be estimated as
\begin{equation}\label{eq-vyyH1-0}
	\begin{aligned}[b]
		\int_{\mathbbm{R}}\int_{\mathbbm{R}_+} |v_{yy}||\nabla v_{x}||\nabla{v}_{yy}| \,\mathrm{d}x \mathrm{d}y
		&\le \| \nabla{v}_{yy}(t) \|_{L^2}\| v_{yy}(t) \|_{L_x^2(L_y^ \infty)}\| \nabla{v}_x(t) \|_{L_x^\infty(L_y^2)}\\
		&\le  \| \nabla{v}_{yy}(t) \|_{L^2}^{3/2}\| \nabla^2{v}(t) \|_{L^2}\| \nabla{v}_{xx} \|_{L^2}^{1/2} \\
		&\le \frac{1}{4} \|\nabla{v}_{yy}(t)\|_{L^2}^2+C \| \nabla^2{v}(t) \|_{L^2}^4  \| \nabla {v}_{xx}(t) \|_{L^2}^2\\
		&\le \frac{1}{4} \|\nabla{v}_{yy}(t)\|_{L^2}^2+C M_0^4  (\| \nabla \Delta {v}(t) \|_{L^2}^2+\| \nabla {v}_{yy}(t) \|_{L^2}^2),
	\end{aligned}
\end{equation}
\begin{equation}\label{eq-vyyH1-1}
	\begin{aligned}[b]
		\int_{\mathbbm{R}}\int_{\mathbbm{R}_+} |v_y||\nabla v_{xy}||\nabla{v}_{yy}| \,\mathrm{d}x \mathrm{d}y
		&\le \| v_y(t) \|_{L^\infty} \| \nabla v_{xy}(t) \|_{L^2} \| \nabla v_{yy}(t) \|_{L^2} \\
        &\le \| v_y(t) \|_{L^\infty}\| \nabla v_{yy}(t) \|_{L^2}^2+\| v_y(t) \|_{L^\infty}\| \nabla v_{xy}(t) \|_{L^2}^2\\
        &\le \| v_y(t) \|_{L^\infty}\| \nabla v_{yy}(t) \|_{L^2}^2+\| v_y(t) \|_{L^\infty}(\| \nabla v_{yy}(t) \|_{L^2}^2+ \| \nabla \Delta{v}(t) \|_{L^2}^2)\\
		&\le C\varepsilon_0 \|\nabla{v}_{yy}(t)\|_{L^2}^2+C \varepsilon_0 \| \nabla \Delta{v}(t) \|_{L^2}^2,
	\end{aligned}
\end{equation}
and
\begin{equation}\label{eq-vyyH1-2}
	\begin{aligned}[b]
		\int_{\mathbbm{R}}\int_{\mathbbm{R}_+} |\nabla v_{y}|^2|\nabla{v}_{yy}| \,\mathrm{d}x \mathrm{d}y
		&\le \| \nabla{v}_{yy}(t) \|_{L^2}\|\nabla v_{y}(t) \|_{L_x^2(L_y^ \infty)}\| \nabla v_{y}(t) \|_{L_x^\infty(L_y^2)} \\
		&\le  \| \nabla{v}_{yy}(t) \|_{L^2}^{3/2}\| \nabla{v}_y(t) \|_{L^2}\| \nabla{v}_{xy} \|_{L^2}^{1/2} \\
		&\le \frac{1}{4} \|\nabla{v}_{yy}(t)\|_{L^2}^2+C \| \nabla{v}_y \|_{L^2}^4  \| \nabla {v}_{xx}(t) \|_{L^2}^2\\
		&\le \frac{1}{4} \|\nabla{v}_{yy}(t)\|_{L^2}^2+C M_0^4  (\| \nabla \Delta {v}(t) \|_{L^2}^2+\| \nabla {v}_{yy}(t) \|_{L^2}^2).
	\end{aligned}
\end{equation}
By employing Cauchy inequality, the remaining terms on the right-hand side of $\eqref{eq-vyyH1}$ are bounded by
$$C(\varepsilon_0+\delta)\| \nabla{v}_{yy}(t) \|_{L^2}^2  +C (\| \nabla{v}(t) \|_{L^2}^2+\| \nabla^2{v}(t) \|_{L^2}^2 ).$$
Hence, applying $\eqref{eq-vyyH1-1}$, for some small $\varepsilon_0$, $M_0$ and $\delta$, $\eqref{eq-vyyH1}$ can be rewritten as
\begin{equation}\label{eq-vyy-nabvyy-djf}
	\begin{aligned}[b]
		&\frac{\mathrm{d}}{\mathrm{d}t}(\| v_{yy} \|_{L^2}^2 +\| \nabla{v}_{yy} \|_{L^2}^2 )+\|\sqrt{\bar{u}_x} v_{yy} \|_{L^2}^2 +\|\sqrt{\bar{u}_x}  \nabla{v}_{yy} \|_{L^2}^2+\| \nabla{v}_{yy} \|_{L^2}^2+\int_{\mathbbm{R}} |\nabla{v}_{yy}(0,y,t)|^2 \,\mathrm{d}y \\
		\le& C (\varepsilon_0+M_0^4)\| \nabla \Delta{v}(t) \|_{L^2}^2 +C (\| \nabla{v}(t) \|_{L^2}^2+\| \nabla^2{v}(t) \|_{L^2}^2 ).
	\end{aligned}
\end{equation}
Integrating $\eqref{eq-vyy-nabvyy-djf}$ over $[0,t]$, the desired estimate $\eqref{eq-vyy}$ is obtained.
  $\hfill\Box$

\begin{lem}\label{lem-p-2}Under the same assumptions as Proposition $\ref{prop-priori}$, for $t\in[0,T]$, it holds that
	\begin{equation}\label{eq-p-2}
		\begin{aligned}
			\int_0^t (\| \nabla^2{\rm{div}}{p}(\tau) \|_{L^2}^2+\| \nabla p(\tau) \|_{L^2}^2+\| \nabla{\rm{div}}{p}(\tau) \|_{L^2}^2 ) \,\mathrm{d}\tau
			\le  CM_0^2+C (\varepsilon_0+M_0^4) \int_0^t \| \nabla \Delta{v}(\tau) \|_{L^2}^2 \,\mathrm{d}\tau.
		\end{aligned}
	\end{equation}
\end{lem}
 {\it\bfseries Proof.}
Square $\eqref{eq-vp-rd}_2$ in the form $\nabla^2{\rm{div}}{p}-\nabla{p}=\nabla^2 v$, and then integrate the resulting equation in $(x,y)\in \mathbbm{R}_+\times \mathbbm{R}$. Consequently, we have
\begin{equation}\label{eq-nabdivp-1}
	\begin{aligned}[b]
		& \| \nabla^2{\rm{div}}{p}(t)\|_{L^2}^2+2\| \nabla{\rm{div}}{p}(t)\|_{L^2}^2+\| \nabla{p}(t)\|_{L^2}^2  \\
		=& \| \nabla^2 v(t) \|_{L^2}^2   +2 \int_0^t\int_{\mathbbm{R}}\int_{\mathbbm{R}_+} {\rm{div}}\{\nabla{\rm{div}}{p}\nabla{p}\} \,\mathrm{d}x \mathrm{d}y \\
		=&  \| \nabla^2 v(t) \|_{L^2}^2   -2 \int_{\mathbbm{R}} \nabla{\rm{div}}{p}(0,y,t)\cdot\nabla{p_1}(0,y,t) \, \mathrm{d}y. \\
	\end{aligned}
\end{equation}
To convert the boundary terms of $p$ to $v$, we make use of the equation $\eqref{eq-rd2-dengjia}$. Then the last term of $\eqref{eq-nabdivp-1}$ can be rewritten as
\begin{equation}\label{eq-nabdivp-2}
	\begin{aligned}[b]
		&\int_{\mathbbm{R}} \nabla{\rm{div}}{p}(0,y,t)\cdot \nabla{p_1}(0,y,t) \, \mathrm{d}y \\
		=&\int_{\mathbbm{R}} {\rm{div}}{p}_x(0,y,t)p_{1x}(0,y,t)+{\rm{div}}{p}_y(0,y,t)p_{1y}(0,y,t) \, \mathrm{d}y \\
		=&\int_{\mathbbm{R}} {\rm{div}}{p}_x(0,y,t)({\rm{div}}{p}(0,y,t)-p_{2y}(0,y,t)) \, \mathrm{d}y
		-\int_{\mathbbm{R}} {\rm{div}}{p}_{yy}(0,y,t)p_{1}(0,y,t) \, \mathrm{d}y \\
		=&\int_{\mathbbm{R}} {\rm{div}}{p}_x(0,y,t)\left({\rm{div}}{p}(0,y,t)-{\rm{div}}{p}_{yy}(0,y,t)+v_{yy}(0,y,t)\right) \, \mathrm{d}y\\
		&-\int_{\mathbbm{R}} {\rm{div}}{p}_{yy}(0,y,t)({\rm{div}}{p}_x(0,y,t)-v_x(0,y,t)) \, \mathrm{d}y .
	\end{aligned}
\end{equation}
Notice that $v_{yy}(0,y,t)=0$. We divide the estimate on the right side of the inequality $\eqref{eq-nabdivp-2}$ into the following three estimates $\eqref{eq-nabdivp-3}$-$\eqref{eq-nabdivp-4}$. The equation $\eqref{eq-vp-rd}_1$ gives ${\rm{div}}{p}_{yy}(0,y,t)=-u_-v_{xyy}(0,y,t)$. Using Cauchy inequality, the right-hand side of $\eqref{eq-nabdivp-2}$ can be estimated as
\begin{equation}\label{eq-nabdivp-3}
	\int_{\mathbbm{R}} {\rm{div}}{p}_{yy}(0,y,t)v_x(0,y,t) \,\mathrm{d}y
	\le C\int_0^t\int_{\mathbbm{R}} v_{xyy}^2(0,y,t) \,\mathrm{d}y+C\int_0^t\int_{\mathbbm{R}} v_{x}^2(0,y,t) \,\mathrm{d}y ,
\end{equation}
\begin{equation}\label{eq-nabdivp-3-1}
\begin{aligned}[b]
	\int_{\mathbbm{R}} {\rm{div}}{p}_x(0,y,t){\rm{div}}{p}(0,y,t) \, \mathrm{d}y
	&\le \int_{\mathbbm{R}} \|{\rm{div}}{p}_x(\cdot,y,t)\|_{L_x^\infty} \|{\rm{div}}{p}(\cdot,y,t)\|_{L_x^\infty} \, \mathrm{d}y \\
	&\le  \| \nabla^2{\rm{div}}{p}(t) \|_{L^2}^\frac{1}{2} \| \nabla{\rm{div}}{p}(t) \|_{L^2} \| {\rm{div}}{p}(t) \|_{L^2}^\frac{1}{2} \\
	&\le \frac{1}{8}  \| \nabla^2{\rm{div}}{p}(t) \|_{L^2}^2  + C ( \| \nabla{\rm{div}}{p}(t) \|_{L^2}^2+ \| {\rm{div}}{p}(t) \|_{L^2}^2)  ,	
\end{aligned}
\end{equation}
and
\begin{equation}\label{eq-nabdivp-4}
\begin{aligned}[b]
	\int_{\mathbbm{R}} {\rm{div}}{p}_{yy}(0,y,t){\rm{div}}{p}_x(0,y,t) \, \mathrm{d}y
	=&\int_{\mathbbm{R}} u_-v_{xyy}(0,y,t){\rm{div}}{p}_x(0,y,t) \, \mathrm{d}y \\
	\le& \int_{\mathbbm{R}} u_-v_{xyy}(0,y,t)\|{\rm{div}}{p}_x(\cdot,y,t)\|_{L_x^\infty} \, \mathrm{d}y \\
	\le& \frac{1}{8}  \| \nabla^2{\rm{div}}{p}(t) \|_{L^2}^2 + C \| \nabla{\rm{div}}{p}(t) \|_{L^2}^2   +C \int_{\mathbbm{R}} |\nabla {v}_{yy}(0,y,t)|^2\,\mathrm{d}y .
\end{aligned}
\end{equation}
Substituting $\eqref{eq-nabdivp-3}$-$\eqref{eq-nabdivp-4}$ into $\eqref{eq-nabdivp-2}$, we can get from $\eqref{eq-nabdivp-1}$ that
\begin{equation}\label{eq-nab2divp-djf}
\begin{aligned}[b]
	 &\| \nabla^2{\rm{div}}{p}(\tau) \|_{L^2}^2+\| \nabla p(\tau) \|_{L^2}^2+\| \nabla{\rm{div}}{p}(\tau) \|_{L^2}^2 \\
	\le& C\| \nabla^2{v}(t) \|_{L^2}^2+  +C \int_{\mathbbm{R}} |\nabla {v}_{yy}(0,y,t)|^2\,\mathrm{d}y+C \int_{\mathbbm{R}} |\nabla {v}(0,y,t)|^2\,\mathrm{d}y+C \| \nabla{\rm{div}}{p}(t) \|_{L^2}^2+C\| {\rm{div}}{p}(t) \|_{L^2}^2.
\end{aligned}
\end{equation}
Integrating $\eqref{eq-nab2divp-djf}$ over $[0,t]$, combining Lemma $\ref{lem-dlnabv}$ and Lemma $\ref{lem-vyy}$, we have
\begin{equation*}
	\int_0^t (\| \nabla^2{\rm{div}}{p}(\tau) \|_{L^2}^2+\| \nabla p(\tau) \|_{L^2}^2+\| \nabla{\rm{div}}{p}(\tau) \|_{L^2}^2 ) \,\mathrm{d}\tau
	\le CM_0^2+C (\varepsilon_0+M_0^4)\int_0^t \| \nabla \Delta{v}(\tau) \|_{L^2}^2 \,\mathrm{d}\tau,	
\end{equation*}
 which completes the proof of Lemma $\ref{lem-p-2}$.
 $\hfill\Box$

\begin{lem}\label{lem-dlnav}Under the same assumptions as Proposition $\ref{prop-priori}$, the estimate for the third order derivatives is obtained as follows:
\begin{equation}\label{eq-nabv-xxx}
\begin{aligned}[b]
    \| \Delta{v}(t) \|_{L^2}^2+\| \nabla \Delta{v}(t) \|_{L^2}^2 \!+\! \int_0^t (\|\sqrt{\bar{u}_x}\Delta{v}(\tau)\|_{L^2}^2+\|\sqrt{\bar{u}_x}\nabla \Delta{v}(\tau)\|_{L^2}^2) \,\mathrm{d}\tau \!+\! \int_0^t \|\nabla \Delta{v}(\tau)\|_{L^2}^2  \,\mathrm{d}\tau& \\
     + \int_0^t\int_{\mathbbm{R}} |\nabla \Delta{v}(0,y,\tau)|^2 \,\mathrm{d}y \mathrm{d}\tau &\le C M_0^2,	
\end{aligned}
\end{equation}
 	for $t\in[0,T]$.
\end{lem}
 {\it\bfseries Proof.}
We can get from $\Delta v\times \Delta \eqref{eq-vp-rd}_1+\nabla{\rm{div}}{p}\cdot \nabla{\rm{div}}\eqref{eq-vp-rd}_2+\nabla \Delta{v}\cdot \nabla \Delta \eqref{eq-vp-rd}_1- \nabla \Delta {\rm{div}}{p}\cdot \nabla{\rm{div}}\eqref{eq-vp-rd}_2$ that
\begin{equation}\label{eq-nav3-feijifen}
	\begin{aligned}[b]
		&\left\{\frac{1}{2} (\Delta v)^2+\frac{1}{2} |\nabla\Delta{v}|^2\right\}_t+\frac{1}{2} \bar{u}_x( (\Delta v)^2+|\nabla\Delta{v}|^2)+\bar{u}_x(\Delta v_x)^2+{\rm{div}} \{\nabla{\rm{div}}{p}\Delta v \}+|\nabla \Delta{v}|^2\\
		&+\left\{\frac{1}{2} v(\Delta v)^2+\frac{1}{2} \bar{u}(\Delta v)^2+\frac{1}{2}  v|\nabla \Delta{v}|^2+\frac{1}{2} \bar{u}|\nabla \Delta{v}|^2\right\}_x\\
		=&-\frac{1}{2} v_x(\Delta v)^2-2 \nabla{v}\cdot \nabla{v}_x\Delta v- \bar{u}_{xxx}v\Delta v-3 \bar{u}_{xx}v_x\Delta v-2 \bar{u}_xv_{xx}\Delta v-\Delta(g(v+\bar{u})_y)\Delta v\\
		&- \frac{1}{2} v_x|\nabla \Delta{v}|^2-\Delta v \nabla{v}_x\cdot\nabla \Delta{v}-\Delta v_x \nabla{v}\cdot\nabla \Delta{v}-2 \nabla{v}_x\cdot\nabla \Delta{v}v_{xx}-2v_{x}\nabla{v}_{xx}\cdot \nabla \Delta{v}\\
		&-2v_{xy}\nabla{v}_y\cdot\nabla \Delta{v}-2v_y \nabla{v}_{xy}\cdot \nabla \Delta{v}-\bar{u}_{xxxx}v\Delta v_x- \bar{u}_{xxx}\nabla{v}\cdot\nabla \Delta{v}-\bar{u}_{xx}\Delta v\Delta v_x-3 \bar{u}_{xxx}v_x\Delta v_x\\
		&-3 \bar{u}_{xx}v\nabla{v}_x\cdot \Delta v_x-2 \bar{u}_{xx}v_{xx}\Delta v_x-2 \bar{u}_xv\nabla v_{xx}\cdot \Delta v_x- \nabla\Delta(g(v+\bar{u})_y)\cdot \nabla\Delta v,
	\end{aligned}
\end{equation}
where we have used from $\eqref{eq-vp-rd}_2$ that
\begin{equation*}
	 |\nabla \Delta{v}|^2=|\nabla \Delta {\rm{div}}{p}|^2+|\nabla {\rm{div}}{p}|^2-2{\rm{div}}\{\Delta {\rm{div}}{p}\nabla {\rm{div}}{p}\}.
\end{equation*}
By simple calculation, we can obtain
\begin{equation}\label{eq-nabdelg-nabdelv}
	\begin{aligned}[b]
		& \nabla\Delta(g(v+\bar{u})_y)\cdot\nabla\Delta v \\
		\sim&  (\nabla{v}\cdot\nabla\Delta v+\bar{u}_x\Delta v_x)(v_x+ \bar{u}_x)^2v_y  +(v_x+ \bar{u}_x)(\nabla{v}_x\cdot \nabla\Delta{v}+\bar{u}_{xx}\Delta v_x)v_y  +(\nabla{v}\cdot\nabla\Delta v+\bar{u}_x\Delta v_x)\bar{u}_{xx}v_y\\
		&+\bar{u}_{xxx}\Delta v_xv_y  +(\nabla{v}\cdot\nabla\Delta v+\bar{u}_x\Delta v_x)v_y^3  +\nabla{v}_y\cdot\nabla\Delta vv_y^2  +(\nabla{v}\cdot\nabla\Delta v+\bar{u}_x\Delta v_x)\Delta vv_y+|\nabla\Delta{v}|^2v_y\\
		&+\left( (v_x+\bar{u}_x)^2+\bar{u}_{xx}+v_{y}^2+\Delta v \right)\nabla{v}_y\cdot \nabla\Delta{v}+ (\nabla{v}\cdot\nabla\Delta v+\bar{u}_x\Delta v_x)\Delta v_y+\{\frac{1}{2} g'(v+\bar{u})|\nabla\Delta{v}|^2\}_y\\
        & +(\nabla{v}\cdot\nabla\Delta v+\bar{u}_x\Delta v_x)(\nabla{v}\cdot \nabla{v}_y+\bar{u}_xv_{xy})+\nabla{v}_x\cdot \nabla\Delta{v}v_{xy}+v_x \nabla{v}_{xy}\cdot \nabla\Delta{v}+\nabla{v}_y\cdot \nabla\Delta{v}v_{yy}\\
        &+v_y \nabla{v}_{yy}\cdot \nabla\Delta{v}+\bar{u}_{xx}v_{xy}\Delta v_x+\bar{u}_x \nabla{v}_{xy}\cdot \nabla\Delta{v}.
	\end{aligned}
\end{equation}
Notice that the first six terms on right-hand side of $\eqref{eq-nav3-feijifen}$ have been estimated in Lemma $\ref{lem-dlnabv}$ after integration in $(x,y)\in\mathbbm{R}_+\times\mathbbm{R}$. So we don't repeat the calculation here. After integration in $(x,y)\in\mathbbm{R}_+\times\mathbbm{R}$, the remaining terms on right-hand side of $\eqref{eq-nav3-feijifen}$ are bounded by
\begin{equation*}
\begin{aligned}
	&\| \nabla{v} \|_{L^\infty} \| \nabla\Delta{v} \|_{L^2}^2+\| \nabla{v} \|_{L^\infty} (\| \nabla{v} \|_{L^2}+\| \nabla^2{v} \|_{L^2}) \| \nabla\Delta{v} \|_{L^2}+C \delta (\| \nabla{v} \|_{L^2} +\| \nabla^2{v} \|_{L^2} )\| \nabla\Delta{v} \|_{L^2} \\
	&+\int_{\mathbbm{R}}\int_{\mathbbm{R}_+} \left(|\nabla^2{v}||\nabla\Delta{v}|+|\nabla{v}|(|\nabla{v}_{xx}|+|\nabla{v}_{xy}|)|\nabla\Delta{v}|+\bar{u}_x (|\nabla{v}_{xx}|+|\nabla{v}_{xy}|)|\nabla\Delta{v}|\right) \,\mathrm{d}x \mathrm{d}y,	
\end{aligned}
\end{equation*}
where we have used the fact that $\eqref{eq-DND-v}$, $\| \nabla{v} \|_{L^\infty}\le C \varepsilon_0 $ and $\| \partial_x^i\bar{u} \|_{L^\infty}\le C \delta, i=2,3$. Integrating $\eqref{eq-nav3-feijifen}$ over $\mathbbm{R}_+\times\mathbbm{R}$, and using $\| \nabla{v} \|_{L^\infty}\le C \varepsilon_0$ again and applying Cauchy inequality, for some small $\varepsilon_0$ and $\delta$, we have
\begin{equation}\label{eq-nav3-3}
	\begin{aligned}[b]
		&\frac{\mathrm{d}}{\mathrm{d}t}(\| \Delta{v}(t) \|_{L^2}^2+\| \nabla \Delta{v}(t) \|_{L^2}^2)\!+\! (\| \sqrt{\bar{u}_x} \Delta{v}(t) \|_{L^2}^2+\| \sqrt{\bar{u}_x}  \nabla \Delta{v}(t) \|_{L^2}^2)   \!+\!\|   \nabla \Delta{v}(t) \|_{L^2}^2  \\
		&\ \ \ \ \ \ \ \ \ \ \ +\int_{\mathbbm{R}} | \nabla \Delta{v}(0,y,t)|^2 \, \mathrm{d}y  \\
		\le& C \left( \| \nabla{v}(t) \|_{L^2}^2+\| \nabla^2{v}(t) \|_{L^2}^2   +\int_{\mathbbm{R}} |{\rm{div}}{p}_x(0,y,t) \Delta v(0,y,t)| \,\mathrm{d}y   \right.\\
		&\ \  \left.+\int_{\mathbbm{R}}\int_{\mathbbm{R}_+} \left(|\nabla^2{v}||\nabla\Delta{v}|+|\nabla{v}|(|\nabla{v}_{xx}|+|\nabla{v}_{xy}|)|\nabla\Delta{v}|+\bar{u}_x (|\nabla{v}_{xx}|+|\nabla{v}_{xy}|)|\nabla\Delta{v}|\right) \,\mathrm{d}x \mathrm{d}y  \right).
	\end{aligned}
\end{equation}
From $\eqref{eq-v-jf}$ and $\eqref{eq-dlnabv}$, the first two terms on the right-hand side of $\eqref{eq-nav3-3}$ are bounded by $CM_0^2$. By Cauchy inequality, $\eqref{eq-dlnabv}$ and $\eqref{eq-p-2}$, the third term on the right-hand side of $\eqref{eq-nav3-3}$ is bounded by
\begin{equation}\label{eq-divpdlv-0}
\begin{aligned}[b]
	&\int_0^t\int_{\mathbbm{R}} |{\rm{div}}{p}_x(0,y,t) \Delta v(0,y,t)| \,\mathrm{d}y \\
	\le& \frac{1}{2} \int_{\mathbbm{R}} |{\rm{div}}{p}_x(0,y,t)|^2 \,\mathrm{d}y  +\frac{1}{2} \int_0^t\int_{\mathbbm{R}} | \Delta v(0,y,t)|^2 \,\mathrm{d}y \\
	\le&\| \nabla{\rm{div}}{p}(t) \|_{L^2}^2+ \| \nabla^2{\rm{div}}{p}(t) \|_{L^2}^2 +\int_{\mathbbm{R}} | \Delta v(0,y,t)|^2 \,\mathrm{d}y .
\end{aligned}
\end{equation}
  Using Lemma $\ref{lem-nab2v-dengjia}$, \eqref{eq-dlnabv} and the Cauchy inequality, the fourth term can be estimated as
\begin{equation}\label{eq-nabdelv-gj1}
	\begin{aligned}[b]
		C\int_{\mathbbm{R}}\int_{\mathbbm{R}_+} |\nabla^2{v}|^2| \nabla \Delta{v}| \,\mathrm{d}x \mathrm{d}y
		\le& \| \nabla^2 v \|_{L_x^2(L_y^\infty)} \| \nabla^2 v \|_{L_x^\infty(L_y^2)}  \| \nabla\Delta{v} \|_{L^2} \\
		\le& \| \nabla^2{v}(t) \|_{L^2}\| \nabla^3{v}(t) \|_{L^2}^2\\
        \le& CM_0(\| \nabla\Delta{v}(t) \|_{L^2}^2+\| \nabla{v}_{yy}(t) \|_{L^2}^2).
	\end{aligned}
\end{equation}
By applying Lemma $\ref{lem-nab2v-dengjia}$, the last two terms on the right-hand side of $\eqref{eq-nav3-3}$ can estimated as
\begin{equation}\label{eq-nabdelv-gj1-2}
 	\begin{aligned}[b]
        &\int_{\mathbbm{R}}\int_{\mathbbm{R}_+} \left(|\nabla{v}|(|\nabla{v}_{xx}|+|\nabla{v}_{xy}|)|\nabla\Delta{v}|+\bar{u}_x (|\nabla{v}_{xx}|+|\nabla{v}_{xy}|)|\nabla\Delta{v}|\right) \,\mathrm{d}x \mathrm{d}y\\
        \le& (\| \nabla{v} \|_{L^\infty}+\| \bar{u}_x \|_{L^\infty} )\| \nabla^3 {v}(t) \|_{L^2}^2\\
        \le& C(\varepsilon_0+\delta) (\| \nabla \Delta{v}(t) \|_{L^2}^2+\|\nabla{v}_{yy}(t) \|_{L^2}^2).
 	\end{aligned}
 \end{equation}
Substitute $\eqref{eq-divpdlv-0}$-$\eqref{eq-nabdelv-gj1-2}$ into $\eqref{eq-nav3-3}$ and choose small $\varepsilon_0$, $\delta$ and $M_0$ such that $((\varepsilon_0+\delta+M_0)\| \nabla\Delta{v}(t) \|_{L^2}^2\le \frac{1}{2}\| \nabla\Delta{v}(t) \|_{L^2}^2$. Consequently, we have
\begin{equation}\label{eq-vH3-djf}
	\begin{aligned}[b]
		&\frac{\mathrm{d}}{\mathrm{d}t}(\| \Delta{v}(t) \|_{L^2}^2+\| \nabla \Delta{v}(t) \|_{L^2}^2)\!+\! (\| \sqrt{\bar{u}_x} \Delta{v}(t) \|_{L^2}^2+\| \sqrt{\bar{u}_x}  \nabla \Delta{v}(t) \|_{L^2}^2)   \!+\!\|   \nabla \Delta{v}(t) \|_{L^2}^2 +\int_{\mathbbm{R}} | \nabla \Delta{v}(0,y,t)|^2 \, \mathrm{d}y  \\
		\le& C \left( \| \nabla{v}(t) \|_{H^1}^2 +\| \nabla{\rm{div}}{p}(t) \|_{L^2}^2 +\| \nabla^2{\rm{div}}{p}(t) \|_{L^2}^2	+(\varepsilon_0+\delta+M_0)\| \nabla{v}_{yy}(t) \|_{L^2}^2+\int_{\mathbbm{R}} (\Delta v(0,y,t))^2 \,\mathrm{d}y  \right).	
	\end{aligned}
\end{equation}
To treat the terms $\| \nabla^2{\rm{div}}{p}(t) \|_{L^2}^2$ and $\| \nabla{v}_{yy}(t) \|_{L^2}^2$ on the right-hand side of $\eqref{eq-vH3-djf}$, substituting the inequalities $\eqref{eq-vyy}$ and $\eqref{eq-p-2}$ into $\eqref{eq-vH3-djf}$ and integrating the resulting inequality over $[0,t]$, we deduce that
\begin{equation}
	\begin{aligned}[b]
		&\| \Delta{v}(t) \|_{H^1}^2+\int_0^t (\| \sqrt{\bar{u}_x} \Delta{v}(\tau) \|_{L^2}^2+\| \sqrt{\bar{u}_x}  \nabla \Delta{v}(\tau) \|_{L^2}^2+\|\nabla \Delta{v}(t) \|_{L^2}^2) \,\mathrm{d}\tau +\int_0^t   \int_{\mathbbm{R}} | \nabla \Delta{v}(0,y,\tau)|^2 \, \mathrm{d}y \mathrm{d}\tau\\
		&\le CM_0^2 \!+\!C\int_0^t  (\| \nabla{v}(\tau) \|_{H^1}^2 \!+\!\| \nabla{\rm{div}}{p}(t) \|_{L^2}^2) \,\mathrm{d}\tau \!+\!C (\varepsilon_0 \!+\!M_0^4) \int_0^t \| \nabla \Delta{v}(\tau) \|_{L^2}^2 \,\mathrm{d}\tau \!+\! C \int_0^t \int_{\mathbbm{R}} (\Delta v(0,y,t))^2 \,\mathrm{d}y \mathrm{d}\tau.
	\end{aligned}
\end{equation}
In the end, we may choose smallness $M_0$ and $\varepsilon_0$, combining the results of Lemmas $\ref{lem-v}$-$\ref{lem-p-1-4}$, the estimate $\eqref{eq-nabv-xxx}$ holds.
 $\hfill\Box $

\begin{lem}\label{lem-nab3v}Under the same assumptions as Lemma $\ref{lem-dlnabv}$, we have
\begin{equation}\label{eq-vxy-jf}
		\| \nabla^3v(t) \|_{L^2}^2+\int_0^t \| \nabla^3{v}(\tau) \|_{L^2}^2  \,\mathrm{d}\tau   \le C M_0^2,
\end{equation}
	for $t\in[0,T]$.
\end{lem}
 {\it\bfseries Proof.}
Since $\eqref{eq-vyy}$ and $\eqref{eq-nabv-xxx}$ give
\begin{equation*}
	\| \nabla v_{yy}(t) \|_{L^2}^2+\int_0^t \| \nabla v_{yy}(\tau) \|_{L^2}^2  \,\mathrm{d}\tau   \le C M_0^2,
\end{equation*}
it follows that
\begin{equation*}
	\| v_{xxx}(t) \|_{L^2}^2\!+\!\int_0^t \| v_{xxx}(\tau) \|_{L^2}^2  \,\mathrm{d}\tau
	\le C(\| \Delta v_x(t) \|_{L^2}^2\!+\!\|  v_{xyy}(t) \|_{L^2}^2)\!+\!C\int_0^t (\| \Delta v_x(\tau) \|_{L^2}^2\!+\!\| v_{xyy}(\tau) \|_{L^2}^2) \,\mathrm{d}\tau \le CM_0^2,
\end{equation*}
and
\begin{equation*}
	\| v_{xxy}(t) \|_{L^2}^2\!+\!\int_0^t \| v_{xxy}(\tau) \|_{L^2}^2  \,\mathrm{d}\tau
	\le C (\| \Delta v_y(t) \|_{L^2}^2\!+\!\| v_{yyy}(t) \|_{L^2}^2)\!+\!C\int_0^t (\| \Delta v_y(\tau) \|_{L^2}^2\!+\!\|  v_{yyy}(\tau) \|_{L^2}^2) \,\mathrm{d}\tau \le CM_0^2.
\end{equation*}
Therefore, combining Lemma $\ref{lem-vyy}$, the proof of Lemma $\ref{lem-nab3v}$ is completed.
  $\hfill\Box$
\\

We only get the estimates of integration $p$ in $(x,y,t)\in \mathbbm{R}_+\times\mathbbm{R}\times[0,t]$ above. In order to obtain the estimates of $p$ in $(x,y)\in \mathbbm{R}_+\times\mathbbm{R}$, it is not feasible to use the previous method like $\eqref{eq-p-cjfgj}$ and $\eqref{eq-divpp-cjf}$, because we can't treat $v_{xt}(0,y,t)$ which appears in $\eqref{eq-divpx0-1}$. However, combining with the existing $H^3$-estimates of $v$, we can get the desired estimate of ${\rm{div}}{p}$ by using $\eqref{eq-vp-rd}_1$, if the estimation of $v_t$ can be obtained.
\begin{lem}\label{lem-vt}Under the same assumptions as Proposition $\ref{prop-priori}$, the estimates of $v_t$ can be obtained. For any $t\in[0,T]$,
	\begin{align}
			\|v_t(t)\|_{H^2}^2 +\int_0^t \|\nabla v_t(\tau)\|_{H^1}^2 \,\mathrm{d}\tau + &\int_0^t\int_{\mathbbm{R}} (|\nabla v_t(0,y,\tau)|^2+|\Delta v_t(0,y,\tau)|^2) \,\mathrm{d}y \mathrm{d}\tau \le C  M_0^2,\label{eq-vt-H2}\\
			&\int_0^t \left(\|{\rm{div}}p_t(\tau)\|_{H^1}^2+\|p_t(\tau)\|_{L^2}^2\right) \, \mathrm{d}\tau \le C M_0^2.\label{eq-pt-jf}
	\end{align}
\end{lem}
 {\it\bfseries Proof.}
 Differentiating $\eqref{eq-vp-rd}$ with respect to $t$, we have
\begin{equation}\label{eq-vp-rd-t-1}
	\begin{cases}
		v_{tt}+v_tv_x+vv_{xt}+(\bar{u} v_{t})_x+(g'(v+\bar{u})v_t)_y+{\rm{div}}{p}_t=0,\\
		-\nabla{\rm{div}}{p}_t+p_t+\nabla{v}_t=0.
	\end{cases}
\end{equation}
Set $w=v_t$ and $z=p_t$. The equations $\eqref{eq-vp-rd-t-1}$ can be rewritten as
\begin{equation}\label{eq-wz}
	\begin{cases}
		w_t+wv_x+vw_x+(\bar{u}w)_x+(g'(v+\bar{u})w)_y+{\rm{div}}{z}=0,\\
		-\nabla{\rm{div}}{z}+z+\nabla w=0.
	\end{cases}
\end{equation}
Comparing $\eqref{eq-vp-rd}$ with $\eqref{eq-wz}$, we note that only the second, third and fifth terms of $\eqref{eq-wz}_1$ are different from $\eqref{eq-vp-rd}_1$ in format. After the previous estimates of $\eqref{eq-vp-rd}$, it is not difficult to find that the first and fourth terms of $\eqref{eq-wz}_1$ and the first and second terms of the $\eqref{eq-wz}_2$ can produce some good terms. And the boundary condition satisfies $w(0,y,t)=0$. In other words, the good terms of $\eqref{eq-wz}$ and $\eqref{eq-vp-rd}$ are consistent in format. In order to give the estimates of $w$ and $z$, we just need to estimate the terms generated by the the second, third and fifth terms of $\eqref{eq-wz}_1$. In the following estimation, we will use a \emph{priori} assumption
\begin{equation}\label{eq-vt-prioriass}
	\sup_{0\le t\le T}\| w(t) \|_{L^\infty}\le C \varepsilon_0.
\end{equation}
Thus, we need to show that $w\in H^2(\mathbbm{R}_+\times\mathbbm{R})$.\\
\indent We can get from $w\times\eqref{eq-wz}_1+z\cdot \eqref{eq-wz}_2+\nabla{w}\cdot \nabla \eqref{eq-wz}_1- \nabla{\rm{div}}{z}\cdot \eqref{eq-wz}_2$ that
\begin{equation}\label{eq-vt-1}
	\begin{aligned}[b]
		&\frac{\mathrm{d}}{\mathrm{d}t} \left(\|w(t)\|_{L^2}^2+\|\nabla{w}(t)\|_{L^2}^2\right)+\|\nabla{w}(t)\|_{L^2}^2+\int_{\mathbbm{R}} |\nabla w(0,y,t)|^2  \,\mathrm{d}y \\
		\le& C\int_{\mathbbm{R}}\int_{\mathbbm{R}_+} (v_xw^2+vw_xw+(g'(v+\bar{u})w)_yw+ \nabla (wv_x)\cdot \nabla{w}+\nabla (vw_x)\cdot \nabla w+\nabla (g'(v+\bar{u})w)_y \cdot \nabla w) \,\mathrm{d}x \mathrm{d}y\\
		\le& C \int_{\mathbbm{R}}\int_{\mathbbm{R}_+} (|\nabla v|w^2+|\nabla{v}||\nabla{w}|^2+w |\nabla^2 v|| \nabla w|+(|\nabla{v}|| \nabla{w}|+\bar{u}_xw_x)(v_yw+w_y)) \,\mathrm{d}x \mathrm{d}y\\
& \ \ \ \ + C\int_{\mathbbm{R}}\int_{\mathbbm{R}_+} \{\frac{1}{2}vw^2 \}_x \,\mathrm{d}x \mathrm{d}y.
	\end{aligned}
\end{equation}
Integrating $\eqref{eq-vt-1}$ over $[0,t]$ and using Cauchy inequality, $\| \nabla{v}(t) \|_{L^\infty}\le C \varepsilon_0$ and $\| w(t) \|_{L^\infty}\le C \varepsilon_0$, we have
\begin{equation}\label{eq-w-jf}
	\begin{aligned}[b]
		& \|w(t)\|_{L^2}^2+\|\nabla{w}(t)\|_{L^2}^2+\int_0^t \|\nabla{w}(\tau)\|_{L^2}^2 \,\mathrm{d}\tau +\int_0^t \int_{\mathbbm{R}} |\nabla w(0,y,\tau)|^2  \,\mathrm{d}y \mathrm{d}\tau  \\
		\le& CM_0^2+C \| \nabla{v}(t) \|_{L^\infty} \int_0^t \| w(\tau) \|_{L^2}^2 \,\mathrm{d}\tau + C\| \nabla{v}(t) \|_{L^\infty} \int_0^t  \|\nabla w(\tau) \|_{L^2}^2 \,\mathrm{d}\tau\\
         &+C\| w(t) \|_{L^\infty}\int_0^t \| \nabla^2 v(\tau) \|_{L^2}\| \nabla w(\tau) \|_{L^2} \,\mathrm{d}\tau +C(\| \nabla{v}(t) \|_{L^\infty}^2+\| \bar{u}_x \|_{L^\infty} )\int_0^t \| \nabla{w}(\tau) \|_{L^2}\| w(\tau)\|_{L^2} \,\mathrm{d}\tau\\
         &+C(\| \nabla{v}(t) \|_{L^\infty}+\| \bar{u}_x \|_{L^\infty} )\int_0^t \| \nabla{w}(\tau) \|_{L^2}^2 \,\mathrm{d}\tau \\
		\le& C(\varepsilon_0+\delta)\int_0^t \| \nabla{w}(\tau) \|_{L^2}^2 \,\mathrm{d}\tau +C \int_0^t \| w(\tau) \|_{L^2}^2 \,\mathrm{d}\tau +C \int_0^t  \| \nabla^2 v(\tau) \|_{L^2}^2 \,\mathrm{d}\tau. 	
	\end{aligned}
\end{equation}
 Since from $\eqref{eq-vp-rd}_1$ gives
\begin{equation*}
	\int_0^t\int_{\mathbbm{R}}\int_{\mathbbm{R}_+} w^2 \,\mathrm{d}x \mathrm{d}y \mathrm{d}\tau
	= \int_0^t\int_{\mathbbm{R}}\int_{\mathbbm{R}_+} v_t^2 \,\mathrm{d}x \mathrm{d}y \mathrm{d}\tau
	\le C \int_0^t (\| \nabla{v}(\tau) \|_{L^2}^2+\| {\rm{div}}{p}(\tau) \|_{L^2}^2)   \,\mathrm{d}\tau
	\le CM_0^2,
\end{equation*}
it follows from $\eqref{eq-w-jf}$ that
\begin{equation}\label{eq-w-jl}
	\| w(t) \|_{H^1}^2+\int_0^t \| \nabla w(\tau) \|_{L^2}^2  \,\mathrm{d}\tau+\int_0^t\int_{\mathbbm{R}} |\nabla w(0,y,\tau)| \,\mathrm{d}y \mathrm{d}\tau\le CM_0^2.
\end{equation}
Similar to $\eqref{eq-delv-nab2v}$, it holds that $\int_{\mathbbm{R}}\int_{\mathbbm{R}_+} (\Delta w)^2 \,\mathrm{d}x \mathrm{d}y=\int_{\mathbbm{R}}\int_{\mathbbm{R}_+} |\nabla^2 w|^2 \,\mathrm{d}x \mathrm{d}y$ owing to $w_y(0,y,t)=0$. The higher order derivatives of $w$ can get from $\nabla{w}\cdot \nabla \eqref{eq-wz}_1+ {\rm{div}}{z}\times {\rm{div}}\eqref{eq-wz}_2+ \Delta{w}\times \Delta\eqref{eq-wz}_1- \Delta{\rm{div}}{z}\times {\rm{div}}\eqref{eq-wz}_2 $ that
\begin{equation}\label{eq-nab2w}
	\begin{aligned}[b]
		&\|\nabla w(t)\|_{L^2}^2+\|\nabla^2{w}(t)\|_{L^2}^2+\int_0^t \|\nabla^2{w}(\tau)\|_{L^2}^2 \,\mathrm{d}\tau +\int_0^t \int_{\mathbbm{R}} (\Delta w(0,y,t))^2 \,\mathrm{d}y \mathrm{d}\tau \\
		\le&CM_0^2+C \left|\int_0^t \int_{\mathbbm{R}}\int_{\mathbbm{R}_+} \left\{ \Delta(wv_x)\Delta w+\Delta ({v}w_x)\Delta w+\Delta\{(g'(v+\bar{u})w)_y\}\Delta w \right\} \,\mathrm{d}x \mathrm{d}y \mathrm{d}\tau \right|,
	\end{aligned}
\end{equation}
where we have the fact that
$$\int_0^t\int_{\mathbbm{R}}\int_{\mathbbm{R}_+} (\nabla (wv_x)\cdot \nabla{w}+\nabla (vw_x)\cdot \nabla w+\nabla (g'(v+\bar{u})w)_y \cdot \nabla w) \,\mathrm{d}x \mathrm{d}y \mathrm{d}\tau\le CM_0^2$$
owing to $\eqref{eq-vt-1}$-$\eqref{eq-w-jl}$. The right-hand side of $\eqref{eq-nab2w}$ can be estimated as follows. By utilizing $\eqref{eq-Lyinf}$, $\eqref{eq-Lxinf}$, $\eqref{eq-del-fh}$, $\eqref{eq-del-g}$ and Cauchy inequality, we obtain
\begin{equation}\label{eq-nab2w-1}
	\begin{aligned}[b]
		&\int_0^t\int_{\mathbbm{R}}\int_{\mathbbm{R}_+} (\Delta(wv_x)\Delta w+\Delta ({v}w_x)\Delta w) \,\mathrm{d}x \mathrm{d}y \mathrm{d}\tau\\
		\le& C \int_0^t\int_{\mathbbm{R}}\int_{\mathbbm{R}_+} (v_x(\Delta w)^2+w \Delta v_x \Delta w+ \nabla w\cdot\nabla v_x \Delta w+\Delta v w_x \Delta w+v\Delta{w_x}\Delta{w}+\nabla{v}\cdot \nabla w_x\Delta w) \,\mathrm{d}x \mathrm{d}y \mathrm{d}\tau\\
		\le& C \int_0^t   \left(\| \nabla{v}(\tau) \|_{L^\infty}\| \nabla^2 w(\tau) \|_{L^2}^2+\| w(\tau) \|_{L^\infty}\| \nabla^3{v}(\tau) \|_{L^2}\| \Delta w(\tau)\|_{L^2}+\| \Delta v(\tau) \|_{L_x^2(L_y^\infty)} \| w_x(\tau) \|_{L_x^\infty(L_y^2)}\| \Delta w(\tau) \|_{L^2} \right) \,\mathrm{d}\tau\\
           & \ \ \ \ +\|\nabla v\|_{L^\infty}\int_0^t  \int_{\mathbbm{R}} (\Delta{w}(0,y,t))^2 \,\mathrm{d}\tau\\
		\le & \left(C \varepsilon_0+\frac{1}{4}\right)\int_0^t \| \nabla^2 w(\tau) \|_{L^2}^2 \,\mathrm{d}\tau +\int_0^t \| \nabla^3{v}(\tau) \|_{L^2}^2 \,\mathrm{d}\tau+\sup\limits_{0\le t\le T}\left\{\| \Delta v(t) \|_{L^2}^2\| \Delta v_y(t) \|_{L^2}^2\right\} \int_0^t \| \nabla{w}(\tau) \|_{L^2}^2  \,\mathrm{d}\tau\\
           & \ \ \ \ +\|\nabla v\|_{L^\infty}\int_0^t  \int_{\mathbbm{R}} (\Delta{w}(0,y,t))^2 \,\mathrm{d}\tau,
	\end{aligned}
\end{equation}
and
\begin{equation}\label{eq-nab2w-2}
	\begin{aligned}[b]
		&\int_0^t \int_{\mathbbm{R}}\int_{\mathbbm{R}_+}\Delta\{(g'(v+\bar{u})w)_y\}\Delta w \,\mathrm{d}x \mathrm{d}y \mathrm{d}\tau\\
		\le& C \int_0^t\int_{\mathbbm{R}}\int_{\mathbbm{R}_+} \left| \left( (v_x+\bar{u}_x)^2+\bar{u}_{xx}+v_y^2+\Delta v \right)v_yw\Delta w+(\Delta v_yw+v_y \Delta w+\nabla{v}_y\cdot \nabla{w})\Delta w \right| \,\mathrm{d}x \mathrm{d}y \mathrm{d}\tau\\
		&+C \int_0^t\int_{\mathbbm{R}}\int_{\mathbbm{R}_+} \left| \nabla{v}\cdot \nabla{v}_yw \Delta w+\nabla{v}\cdot \nabla{w}v_y \Delta w+\bar{u}_xv_{xy}w \Delta w+\bar{u}_xw_xv_y \Delta w \right| \,\mathrm{d}x \mathrm{d}y \mathrm{d}\tau \\
		&+C \int_0^t\int_{\mathbbm{R}}\int_{\mathbbm{R}_+} \left| \left( (v_x+\bar{u}_x)^2+\bar{u}_{xx}+v_y^2+\Delta v \right)w_y \Delta w+\nabla{v}\cdot \nabla{w}_y\Delta w+\bar{u}_xw_{xy}\Delta w \right| \,\mathrm{d}x \mathrm{d}y \mathrm{d}\tau \\
		\le& \left(C \varepsilon_0+\frac{1}{4}\right)\int_0^t \| \nabla^2 w(\tau) \|_{L^2}^2 \,\mathrm{d}\tau+C \int_0^t (\| w(\tau) \|_{L^2}^2+\|\nabla w(\tau) \|_{L^2}^2+\| \nabla^2{v}(\tau) \|_{L^2}^2+\| \nabla^3{v}(\tau) \|_{L^2}^2 ) \,\mathrm{d}\tau \\
		&+C\int_0^t \| \nabla{v}_y(\tau) \|_{L_x^2(L_y^\infty)}\| \nabla{w}(\tau) \|_{L_x^\infty(L_y^2)}\| \Delta w(\tau) \|_{L^2} \,\mathrm{d}\tau \\
		 \le& \left(C \varepsilon_0+\frac{1}{4}\right)\int_0^t \| \nabla^2 w(\tau) \|_{L^2}^2 \,\mathrm{d}\tau+C \int_0^t (\| w(\tau) \|_{L^2}^2+\|\nabla w(\tau) \|_{L^2}^2+\| \nabla^2{v}(\tau) \|_{L^2}^2+\| \nabla^3{v}(\tau) \|_{L^2}^2 ) \,\mathrm{d}\tau.
	\end{aligned}
\end{equation}
In deriving the last inequality of $\eqref{eq-nab2w-2}$ above, we have used that $\|\nabla{v}_y(t)\|_{L^2} $ and $\|\nabla{v}_{yy}(t)\|_{L^2} $ are bounded, which yields
\begin{equation*}
	\begin{aligned}
		&\int_0^t \| \nabla{v}_y(\tau) \|_{L_x^2(L_y^\infty)}\| \nabla{w}(\tau) \|_{L_x^\infty(L_y^2)}\| \Delta w(\tau) \|_{L^2} \,\mathrm{d}\tau\\
		\le&\int_0^t \| \nabla{v}_y(\tau) \|_{L^2}^{\frac{1}{2} }\| \nabla{v}_{yy}(\tau) \|_{L^2}^{\frac{1}{2} }\| \nabla{w}(\tau) \|_{L^2}^{\frac{1}{2} }\| \nabla^2{w}(\tau) \|_{L^2}^{\frac{3}{2} }\,\mathrm{d}\tau \\
		\le &\frac{1}{4}\int_0^t \| \nabla^2{w}(\tau) \|_{L^2}^2 \,\mathrm{d}\tau+\sup\limits_{0\le t\le T}\left\{\| \nabla v_y(t) \|_{L^2}^2\| \nabla v_{yy}(t) \|_{L^2}^2\right\} \int_0^t \| \nabla{w}(\tau) \|_{L^2}^2  \,\mathrm{d}\tau.
	\end{aligned}
\end{equation*}
Substituting $\eqref{eq-nab2w-1}$ and $\eqref{eq-nab2w-2}$ into $\eqref{eq-nab2w}$, combining $\eqref{eq-w-jf}$, we get the desired inequality $\eqref{eq-vt-H2}$. The proof of $\eqref{eq-pt-jf}$ is similar to Lemma $\ref{lem-p-1-4}$. We omit the details.
 $\hfill\Box$

After we get the $H^2$-estimates of $v_t$ and $H^3$-estimates of $v$, the estimates of the lower derivatives of $p$ and ${\rm{div}}{p}$ can be easily obtained by the equation $\eqref{eq-vp-rd}_1$. Using relation between $p_1$ and $p_2$ subtly, we can get the higher derivatives estimation of $p$ and ${\rm{div}}{p}$.
\begin{lem}\label{lem-p-gj}Under the same assumptions as Proposition $\ref{prop-priori}$, we can obtain the estimates of perturbation $p$.
	\begin{equation}
		\|{\rm{div}}p(t)\|_{H^3}^2+\|p(t)\|_{H^3}^2+\int_0^t (\| \nabla^3{\rm{div}}{p}(\tau)\|_{L^2}^2+\| \nabla^3{p}(\tau)\|_{L^2}^2+\| \nabla^2{p}(\tau)\|_{L^2}^2) \,\mathrm{d}\tau\le C M_0^2,
	\end{equation}
	for $t\in[0,T]$.
\end{lem}
{\it\bfseries Proof.}
From $\eqref{eq-vp-rd}$, we can obtain a series of estimates of ${\rm{div}}{p}$:
\begin{align}
    &\| {\rm{div}}{p}(t)\|_{L^2}^2\le C\| v_t(t)\|_{L^2}^2+C\| v(t)\|_{H^1}^2,\nonumber\\
    &\| \nabla {\rm{div}}{p}(t)\|_{L^2}^2\le C\| v_t(t)\|_{H^1}^2+C\| v(t)\|_{H^2}^2,\nonumber\\
	&\|\nabla^2 {\rm{div}}{p}(t)\|_{L^2}^2\le C\| v_t(t)\|_{H^2}^2+C\| v(t)\|_{H^3}^2,\label{eq-p-xyjf-1}\\
	&\|\nabla \Delta {\rm{div}}{p}(t)\|_{L^2}^2\le \|\nabla {\rm{div}}{p}(t)\|_{L^2}^2+\| \nabla \Delta{v}(t)\|_{L^2}^2.\label{eq-p-xyjf-2}
\end{align}
Consequently, the corresponding estimate of $p$ can be obtained:
\begin{equation*}
	\begin{aligned}
		&\| p(t)\|_{L^2}^2\le \| \nabla{\rm{div}}{p}(t) \|_{L^2}^2+\| \nabla{v}(t) \|_{L^2}^2,\\
		&\| \nabla p(t) \|_{L^2} \le \| \nabla^2{\rm{div}}{p}(t) \|_{L^2}^2+\| \nabla^2{v}(t) \|_{L^2}^2.
	\end{aligned}
\end{equation*}
Next, we consider the higher regularity for ${\rm{div}}{p}$ and $p$. By the definition,
$$\nabla^3{\rm{div}}{p}=({\rm{div}}{p}_{xxx},{\rm{div}}{p}_{xxy},{\rm{div}}{p}_{xyy},{\rm{div}}{p}_{yyy}).$$
This shows that $\nabla^3{\rm{div}}{p}$ consists of $\nabla{\rm{div}}{p}_{xx}$ and $\nabla{\rm{div}}{p}_{yy}$. Once we give the $L^2$-estimates of $\nabla{\rm{div}}{p}_{xx}$, the $L^2$-estimates of $\nabla{\rm{div}}{p}_{yy}$ can be obtained immediately according to $\eqref{eq-p-xyjf-2}$. However, applying $\partial_x^2$ to $\eqref{eq-vp-rd}_2$ will produce good terms ($\| {\rm{div}}{p}_{xx} \|_{H^1}^2, \| p_{xx} \|_{L^2}^2 $) and boundary term $p_{1xx}(0,y,t)$ that we can't solve. To overcome it, we rewrite the equation $\eqref{eq-vp-rd}_2$ as $\eqref{eq-rd2-dengjia}$. We find that $\nabla\partial_x\eqref{eq-rd2-dengjia}_2$ can produce good terms ($\| \nabla{\rm{div}}{p}_{xx} \|_{L^2}^2, \| \nabla p_{1x} \|_{L^2}^2,  \| p_{1xxx} \|_{L^2}^2 $) which can control the boundary term  $p_{1xx}(0,y,t)$. Therefore, we focus on the equations $\eqref{eq-rd2-dengjia}$ in the later estimates.\\
\indent Applying the differential operator $\partial_x\partial_y$ and $\nabla \partial_x$ to $\eqref{eq-rd2-dengjia}_1$, respectively, we get
\begin{equation*}
\begin{aligned}
	{\rm{div}}{p}_{xxy}-p_{1xy}=v_{xxy},\\
	{\rm{div}}{p}_{xxx}-p_{1xx}=v_{xxx}.	
\end{aligned}
\end{equation*}
 Squaring these two equations, it holds that
 \begin{align}
 	({\rm{div}}{p}_{xxy})^2+(p_{1xy})^2+2(p_{1xyy})^2-2\left\{{\rm{div}}{p}_{xx}p_{1xy} \right\}_y=v_{xxy}^2-2p_{1xxx}p_{1xyy},\label{eq-divp-H3-1}\\
 	({\rm{div}}{p}_{xxx})^2+(p_{1xx})^2+2(p_{1xxx})^2-2\left\{{\rm{div}}{p}_{xx}p_{1xx} \right\}_x=v_{xxx}^2-2p_{1xyy}p_{1xxx}.\label{eq-divp-H3-2}
 \end{align}
 Integrating $\eqref{eq-divp-H3-1}$ and $\eqref{eq-divp-H3-2}$ over $\mathbbm{R}_+\times\mathbbm{R}$, and then adding the resulting equations, we obtain
 \begin{equation}\label{eq-divp-H3-3}
 \begin{aligned}[b]
  	&\int_{\mathbbm{R}}\int_{\mathbbm{R}_+} |\nabla {\rm{div}}{p}_{xx}|^2+|\nabla p_{1x}|^2+2(p_{1xyy})^2+2(p_{1xxx})^2 \,\mathrm{d}x \mathrm{d}y\\
 	=&\int_{\mathbbm{R}}\int_{\mathbbm{R}_+} |\nabla{v}_{xx}|^2 \,\mathrm{d}x \mathrm{d}y-4 \int_{\mathbbm{R}}\int_{\mathbbm{R}_+} p_{1xyy}p_{1xxx} \,\mathrm{d}x \mathrm{d}y -2 \int_{\mathbbm{R}} {\rm{div}}{p}_{xx}(0,y,t)p_{1xx}(0,y,t) \,\mathrm{d}y.	
 \end{aligned}
 \end{equation}
 By $\eqref{eq-p-xyjf-1}$ and $\eqref{eq-p1p2}$, the last two terms of $\eqref{eq-divp-H3-3}$ are bounded by
 \begin{equation}\label{eq-divp-H3-4}
 	\begin{aligned}[b]
 		\int_{\mathbbm{R}}\int_{\mathbbm{R}_+} p_{1xyy}p_{1xxx} \,\mathrm{d}x \mathrm{d}y
 		\le \frac{1}{2} \int_{\mathbbm{R}}\int_{\mathbbm{R}_+} (p_{2xxy}+p_{1xxx})^2 \,\mathrm{d}x \mathrm{d}y
 		\le \frac{1}{2} \int_{\mathbbm{R}}\int_{\mathbbm{R}_+} ({\rm{div}}{p}_{xx})^2 \,\mathrm{d}x \mathrm{d}y,
 	\end{aligned}
 \end{equation}
 and
 \begin{equation}\label{eq-divp-H3-5}
 	\begin{aligned}[b]
 		\int_{\mathbbm{R}} {\rm{div}}{p}_{xx}(0,y,t)p_{1xx}(0,y,t) \,\mathrm{d}y
 		&\le C_\varepsilon \int_{\mathbbm{R}} \| {\rm{div}}{p}_{xx} \|_{L_x^\infty}^2 \,\mathrm{d}y+\varepsilon \int_{\mathbbm{R}} \|{p}_{1xx} \|_{L_x^\infty}^2 \,\mathrm{d}y\\
 		&\le \varepsilon(\| {\rm div}p_{xxx}(t) \|_{L^2}^2+\| p_{1xx}(t) \|_{L^2}^2+\|p_{1xxx}(t)\|_{L^2}^2 )+ C_\varepsilon\| {\rm div}p_{xx}(t) \|_{L^2}^2,
 	\end{aligned}
 \end{equation}
 where $\varepsilon$ is a small positive constant and $C_\varepsilon$ is a fixed positive constant depending on $\varepsilon$. Substituting $\eqref{eq-divp-H3-4}$ and $\eqref{eq-divp-H3-5}$ into $\eqref{eq-divp-H3-3}$, combining $\eqref{eq-vxy-jf}$ and $\eqref{eq-p-xyjf-1}$, we obtain
 \begin{equation}\label{eq-nabp3-1}
 	\int_{\mathbbm{R}}\int_{\mathbbm{R}_+} |\nabla {\rm{div}}{p}_{xx}|^2+|\nabla p_{1x}|^2+(p_{1xyy})^2+(p_{1xxx})^2 \,\mathrm{d}x \mathrm{d}y\le CM_0^2.
 \end{equation}
By $\eqref{eq-p-xyjf-2}$, we get further
 \begin{equation*}
 \| {\rm{div}}{p}_{xyy}(t) \|_{L^2}^2 \le C(\| \Delta {\rm{div}}{p}_{x}(t) \|_{L^2}^2+\| {\rm{div}}{p}_{xxx}(t) \|_{L^2}^2)\le CM_0^2,
 \end{equation*}
 and
 \begin{equation*}
 \| {\rm{div}}{p}_{yyy}(t) \|_{L^2}^2 \le C(\| \Delta {\rm{div}}{p}_{y}(t) \|_{L^2}^2+\| {\rm{div}}{p}_{xxy}(t) \|_{L^2}^2)\le CM_0^2,
 \end{equation*}
 which yields
 \begin{align}
 	&\| \nabla^3{\rm{div}}{p}(t) \|_{L^2}^2\le CM_0^2,\nonumber\\
 	&\| \nabla^2 p(t) \|_{L^2}^2\le \| \nabla^3{\rm{div}}{p}(t) \|_{L^2}^2+\| \nabla^3{v}(t)\|_{L^2}^2\le CM_0^2.\nonumber
 \end{align}
 Next, we try to get the $L^2$-estimates of $\nabla^3 p$. By the definition
 \begin{equation*}
 	\nabla^3p=
 	\left(
  \begin{array}{cccc}
    p_{1xxx} & p_{1xxy} & p_{1xyy} & p_{1yyy}\\
    p_{2xxx} & p_{2xxy} & p_{2xyy} & p_{2yyy}\\
  \end{array}
\right).
 \end{equation*}
 The equation $\eqref{eq-p1p2}$ deduce that
	\begin{align}
	    &{\rm{div}}{p}_{xx}-p_{1xxx}=p_{2xxy}=p_{1xyy}={\rm{div}}{p}_{yy}-p_{2yyy},\label{eq-p1xxx}\\
		&p_{2xxx}=p_{1xxy}={\rm{div}}{p}_{xy}-p_{2xyy}={\rm{div}}{p}_{xy}-p_{1yyy}.\label{eq-p2xxx}
	\end{align}
 The $L^2$-estimates of $p_{1xxx}$ and $p_{1xyy}$ have been obtained in $\eqref{eq-nabp3-1}$, which yields $L^2$-estimates of $p_{2xxy}$ and $p_{2yyy}$ in the equation $\eqref{eq-p1xxx}$ can be estimated, combining $\eqref{eq-p-xyjf-1}$. Thus, in order to complete the $L^2$-estimates of $\nabla^3 p$, we just need to estimate any third order derivative of $p$ in $\eqref{eq-p2xxx}$. Once this estimate is obtained, the remaining term can be easily obtained according to the relation $\eqref{eq-p2xxx}$. Now, we give an estimate of $p_{2xyy}$ for a brief explanation.\\
\indent Just like the previous method, differentiating $\eqref{eq-rd2-dengjia}$ respect to $x$ and $y$, and squaring the resulting equation in the form of ${\rm{div}}{p}_{xyy}-p_{2xy}=v_{xyy}$, we have
 \begin{equation}\label{eq-p-H3-1}
 	({\rm{div}}{p}_{xyy})^2+(p_{2xy})^2+(p_{2xyy})^2-2\left\{ {\rm{div}}{p}_{xy}p_{2xy} \right\}_y=v_{xyy}^2-2p_{1xxy}p_{2xyy}.
 \end{equation}
Integrate $\eqref{eq-p-H3-1}$ over $\mathbbm{R}_+\times\mathbbm{R}$. Notice that
\begin{equation*}
	\int_{\mathbbm{R}}\int_{\mathbbm{R}_+} p_{1xxy}p_{2xyy} \,\mathrm{d}x \mathrm{d}y \le \frac{1}{2} \int_{\mathbbm{R}}\int_{\mathbbm{R}_+} ({\rm{div}}{p}_{xy})^2 \,\mathrm{d}x \mathrm{d}y.
\end{equation*}
It follows from $\eqref{eq-p-H3-1}$ that
\begin{equation*}
	\| {\rm{div}}{p}_{xyy}(t)\|_{L^2}^2+\| p_{2xy}(t)\|_{L^2}^2+ \| p_{2xyy}(t) \|_{L^2}^2\le CM_0^2.
\end{equation*}
Consequently, using $\eqref{eq-p2xxx}$, we have
\begin{equation*}
	\begin{aligned}
		&\| p_{1yyy}(t) \|_{L^2}^2=\| p_{2xyy}(t) \|_{L^2}^2\le CM_0^2,\\
		&\| p_{2xxx}(t) \|_{L^2}^2=\| p_{1xxy}(t) \|_{L^2}^2 \le C( \| {\rm{div}}{p}_{xy}(t) \|_{L^2}^2+\|p_{2xyy}(t) \|_{L^2}^2)\le CM_0^2,
	\end{aligned}
\end{equation*}
which yields $\| \nabla^3 p(t) \|_{L^2}^2\le CM_0^2$. Using the similar method, we can get
\begin{equation*}
	\int_0^t (\| \nabla^3{\rm{div}}{p}(\tau)\|_{L^2}^2+\| \nabla^3{p}(\tau)\|_{L^2}^2+\| \nabla^2{p}(\tau)\|_{L^2}^2) \,\mathrm{d}\tau\le CM_0^2.
\end{equation*}
 Therefore, the proof of Lemma $\ref{lem-p-gj}$ is completed.
$\hfill\Box$

\subsection{Large-time Behavior} \label{subsec-4}
Combining the standard theory of the existence and uniqueness of the local solution with the \emph{a priori} estimates, one can extend the local solution for problem \eqref{eq-vp-rd} globally, that is
\begin{equation*}
\begin{cases}
	v\in C^0([0,\infty);H^3), \quad \nabla v \in L^2(0,\infty;H^2),\\[1mm]
	p\in C^0([0,\infty);H^3)\cap L^2(0,\infty;H^3), \quad  {\rm{div}}{p}\in C^0([0,\infty);H^3)\cap L^2(0,\infty;H^3),\\[1mm]
	v_t\in C^0([0,\infty);H^2), \quad \nabla v_t \in L^2(0,\infty;H^1), \quad  p_t\in L^2(0,\infty;L^2),  \quad  {\rm{div}}{p}_t\in L^2(0,\infty;H^1).
\end{cases}
\end{equation*}
Based on this, we can deduce that for $t\ge 0$,
\begin{equation}\label{eq-bigtime}
\begin{aligned}[b]
	&\| v(t) \|_{H^3}^2+\| p(t) \|_{H^3}^2+\| {\rm{div}}{p}(t) \|_{H^3}^2+\int_0^\infty (\| \nabla{v}(t) \|_{H^2}^2+\| p(t) \|_{H^3}^2+\| {\rm{div}}{p}(t) \|_{H^3}^2 )  \,\mathrm{d}t <\infty,\\
	&\| v_t(t) \|_{H^2}^2+\int_0^\infty (\| \nabla{v}_t(t) \|_{H^1}^2+\| p_t(t) \|_{L^2}^2+\| {\rm{div}}{p}_t(t) \|_{H^1}^2 )  \,\mathrm{d}t <\infty.
\end{aligned}
\end{equation}
In other to show the desired large time behavior in Theorem $\ref{thm-main}$, using Gagliardo-Nirenberg inequality
$$\| f \|_{L^\infty}\le C \| f \|_{L^2}^\frac{1}{2} \| \nabla^2 f\|_{L^2}^\frac{1}{2}, \ \ \ \ \| \nabla f \|_{L^\infty}\le C \| f \|_{L^2}^\frac{1}{3} \| \nabla^3 f\|_{L^2}^\frac{2}{3},$$
we just need to proof that
\begin{equation}\label{eq-divpt}
\begin{aligned}[b]
	&\int_0^\infty \left|\frac{\mathrm{d}}{\mathrm{d}t}( \|\nabla{v}(t)\|_{L^2}^2+\|\nabla^2{v}(t)\|_{L^2}^2) \right| \,\mathrm{d}t < \infty,\\
	&\int_0^\infty \left|\frac{\mathrm{d}}{\mathrm{d}t}(\| \nabla{\rm{div}}{p} (t)\|_{L^2}^2+\| {\rm{div}}{p}(t) \|_{L^2}^2  +\|p(t)\|_{L^2}^2 )\right| \,\mathrm{d}t < \infty.	
\end{aligned}
\end{equation}
In fact, we can obtain
\begin{equation*}
	\int_0^\infty \left|\frac{\mathrm{d}}{\mathrm{d}t}( \|\nabla{v}(t)\|_{L^2}^2\!+\!\|\nabla^2{v}(t)\|_{L^2}^2) \right| \,\mathrm{d}t
	\le C \int_0^\infty (\|\nabla{v}(t)\|_{L^2}^2+\|\nabla{v}_t(t)\|_{L^2}^2+\|\nabla^2{v}(t)\|_{L^2}^2+\|\nabla^2{v}_t(t)\|_{L^2}^2)  \, \mathrm{d}t
	 < \infty,
\end{equation*}
and
\begin{equation*}
	\begin{aligned}[b]
		&\int_0^\infty \left|\frac{\mathrm{d}}{\mathrm{d}t}(\| \nabla{\rm{div}}{p} (t)\|_{L^2}^2+\| {\rm{div}}{p}(t) \|_{L^2}^2  +\|p(t)\|_{L^2}^2 )\right| \,\mathrm{d}t\\
		\le&C \int_0^\infty (\| \nabla{\rm{div}}{p}(t) \|_{L^2}^2+\| \nabla{\rm{div}}{p}_t(t) \|_{L^2}^2
		+\| {\rm{div}}{p}(t) \|_{L^2}^2+\| {\rm{div}}p_t (t) \|_{L^2}^2+\|p(t) \|_{L^2}^2+\|p_t(t) \|_{L^2}^2)  \,\mathrm{d}t \\
		<& \infty.
	\end{aligned}
\end{equation*}
Therefore, applying $\eqref{eq-bigtime}$ and $\eqref{eq-divpt}$, we can get
$$\| \nabla{v}(t) \|_{L^2}, \ \ \| \nabla^2{v}(t) \|_{L^2}, \ \ \| \nabla{\rm{div}}{p} (t)\|_{L^2},  \ \ \| {\rm{div}}{p}(t) \|_{L^2},  \ \ \|p(t)\|_{L^2}\rightarrow 0, \quad \text{as } \ t\rightarrow \infty,$$
it follows that as $t\rightarrow \infty$,
\begin{equation*}
	\begin{aligned}[b]
	&\|\nabla^k v(t) \|_{L^\infty}\le C \| \nabla^k v(t) \|_{L^2}^{1/2}\| \nabla^{k+2} v(t) \|_{L^2}^{1/2} \rightarrow 0, \ \ \ \ k=0,1 \\
    &\| p(t) \|_{L^\infty}\le C \| p(t) \|_{L^2}^{1/2}\| \nabla^2{p}(t) \|_{L^2}^{1/2} \rightarrow 0, \\
    &\| \nabla p(t) \|_{L^\infty}\le C \| p(t) \|_{L^2}^{1/3}\| \nabla^3{p}(t) \|_{L^2}^{2/3} \rightarrow 0,\\
    &\| \nabla {\rm{div}}{p}(t) \|_{L^\infty} \le C \| \nabla {\rm{div}}{p}(t) \|_{L^2}^{1/2}\| \nabla^{3} {\rm{div}}{p}(t) \|_{L^2}^{1/2} \rightarrow 0.
    \end{aligned}
\end{equation*}
Hence, we complete the proof of $\eqref{eq-asymptotic-behavior}$ in Theorem $\ref{thm-main}$.	

\section{Convergence Rate for the Planar Stationary Solution}\label{sec-5}
\subsection{Time and Space Weighted Energy Estimates}
\begin{prop}
	Suppose the same assumptions as in Theorem $\ref{thm-main}$ hold. Then, there exists a positive constant $C$ such that, for an arbitrary constant $\varepsilon>0$ it holds that
	\begin{equation}\label{eq-convergence-rate}
	\begin{aligned}[b]
		\sum_{l=0}^{2}(1+t)^{\alpha+l+\varepsilon} \| \partial_y^lv(t)\|_{H^{3-l}}^2 +\sum_{l=0}^{2}\int_0^t (1+\tau)^{\alpha+l+\varepsilon}\|\partial_y^l\nabla v(\tau)\|_{H^{2-l}}^2 \,\mathrm{d}\tau \le C(1+t)^\varepsilon  M_\alpha^2,
	\end{aligned}
	\end{equation}
	and
	\begin{equation}\label{eq-q-convergence-rate}
		\begin{aligned}[b]
			(1+t)^{\alpha+1+\varepsilon}(\| {\rm{div}}{p}_y(t) \|_{H^2}^2+\| p_{1y}(t)\|_{H^2}^2+\| p_{2}(t) \|_{H^3}^2)
			\le C(1+t)^\varepsilon M_\alpha^2.
		\end{aligned}
	\end{equation}
\end{prop}

\indent Applying the inequality $\eqref{eq-v-Linfty}$-$\eqref{eq-nabvy-inf}$ and the estimates $\eqref{eq-convergence-rate}$-$\eqref{eq-q-convergence-rate}$, we can get the decay rate of the solution $(u,q)$ of $\eqref{yfc}$ (see $\eqref{eq-convergence-rate-thm}$, $\eqref{eq-q-sjsjgj}$) and complete the proof of the Theorem $\ref{thm-main}$. We devote ourselves to show the estimates $\eqref{eq-convergence-rate}$ and $\eqref{eq-q-convergence-rate}$ as follows.
\begin{lem}\label{Lem-v-H2-gamma}Suppose that $v(x,y,t)$ is a solution to the problem $\eqref{eq-vp-rd}$-$\eqref{eq-vp-initial}$ satisfying $v\in H^3(\mathbbm{R}_+\times\mathbbm{R}\times[0,T])$ for $T>0$. Assume that the initial data satisfy $v_0 \in L_{\alpha,2}^2(\mathbbm{R}_+\times \mathbbm{R})$, $\alpha\ge0$. For any $\beta,\gamma\in[0,\alpha]$, there are positive constants $\varepsilon_1$ and $C=C(\varepsilon_1)$, which are independent of $T$ and $\gamma$ such that if $\sup_{0\le \tau\le T}\|v(\tau)\|_{H^3}+\delta\le\varepsilon_1$, then it holds that for $t\in[0,T]$,
\begin{equation}\label{eq-vsjl-alpha}
\begin{aligned}[b]
&(1+t)^\gamma|v(t)|_{\beta,2}^2 +\beta \int_0^t (1+\tau)^\gamma|v(\tau)|_{\beta-1,2}^2  \,\mathrm{d}\tau
+\int_0^t (1+\tau)^\gamma (|\nabla{v}(t)|_{\beta,1}^2 +|\nabla{\rm{div}}{p}(t)|_{\beta}^2+ |{\rm{div}}{p}(t)|_{\beta}^2+|p(t)|_{\beta}^2)   \,\mathrm{d}\tau\\
&+\int_0^t (1+\tau)^\gamma \int_{\mathbbm{R}} |\nabla{v}(0,y,\tau)|^2 \,\mathrm{d}y \mathrm{d}\tau  + \int_0^t (1+\tau)^\gamma \int_{\mathbbm{R}} |\nabla^2{v}(0,y,\tau)|^2 \,\mathrm{d}y \mathrm{d}\tau \\
\le&  C|v_0|_{\beta,2}^2+C \gamma \int_0^t (1+\tau)^{\gamma-1} |v(\tau)|_{\beta,2}^2 \,\mathrm{d}\tau
+C\beta \int_0^t (1+\tau)^\gamma (\| p (\tau)\|_{L^2}^2+\| {\rm{div}}{p}(\tau) \|_{L^2}^2+\| \nabla^2{v}(\tau) \|_{L^2}^2+\| \nabla{v}(\tau) \|_{L^2}^2 ) \,\mathrm{d}\tau.
\end{aligned}
\end{equation}
In particular, $\gamma=\beta=0$, for any $t\in[0,T]$,
\begin{equation}\label{eq-vsjl2}
\begin{aligned}[b]
 \|v(t)\|_{H^2}^2+\int_0^t \ (\|\nabla{v}(\tau) \|_{H^1}^2+\|\nabla{\rm{div}}{p}(\tau)  \|_{L^2}^2+\|{\rm{div}}{p}(\tau)  \|_{L^2}^2+\| p(\tau) \|_{L^2}^2)\,d\tau&\\
 +\int_0^t (1+\tau)^\gamma \int_{\mathbbm{R}} (|\nabla{v}(0,y,\tau)|^2 +|\nabla^2{v}(0,y,\tau)|^2) \,\mathrm{d}y \mathrm{d}\tau&\le C M_0^2.
\end{aligned}
\end{equation}
\end{lem}
{\it\bfseries Proof.}
Multiplying $\eqref{eq-vnabv}$ by $2(1+x)^\beta$, $\beta\in[0,\alpha]$ and then integrating it over $\mathbbm{R}_+\times\mathbbm{R}$, we have
\begin{equation}\label{eq-vnabv-beta}
	\begin{aligned}[b]
		 &\frac{\mathrm{d}}{\mathrm{d}t}(|v(t)|_\beta^2+|\nabla{v}(t)|_\beta^2)+|u_+|\beta(|v(t)|_{\beta-1}^2+|\nabla{v}(t)|_{\beta-1}^2)+2|\nabla{v}(t)|_\beta^2+|u_-|\int_{\mathbbm{R}} |\nabla{v}(0,y,t)|^2 \,\mathrm{d}y \\
		 \le&  \beta \int_{\mathbbm{R}}\int_{\mathbbm{R}_+}(1+x)^{\beta-1} (\frac{2}{3}|v|^3+|v||\nabla{v}|^2+2|p_1||v|) \,\mathrm{d}x \mathrm{d}y \\
		 & \ \ \ \ \ \ \ \ \ \ \ \ \ \ \ \ \ \  +C \int_{\mathbbm{R}}\int_{\mathbbm{R}_+} (1+x)^{\beta}(|\nabla{v}|^3+|\bar{u}_{xx}||v||\nabla{v}|+\bar{u}_x|\nabla{v}|^2) \,\mathrm{d}x \mathrm{d}y.
	\end{aligned}
\end{equation}
Here we have used the fact that $u_-\le \bar{u}\le u_+\le0$ and
\begin{equation*}
- \beta (1+x)^{\beta-1}\bar{u} \ge  -u_+\beta(1+x)^{\beta-1}.
\end{equation*}
Under the assumptions of Lemma $\ref{Lem-v-H2-gamma}$, by Sobolev inequality, for some positive constant $C$, it holds
$$|v(t)|_{L^\infty}\le C\varepsilon_1, \ \ \ \ |\nabla{v}(t)|_{L^\infty}\le C\varepsilon_1.$$
Thus, the first, second and fourth terms on the right-hand side of $\eqref{eq-vnabv-beta}$ are bounded by
\begin{equation}\label{eq-v-beta-1}
	 \beta \int_{\mathbbm{R}}\int_{\mathbbm{R}_+}(1+x)^{\beta-1} |v|(|v|^2+|\nabla{v}|^2) \,\mathrm{d}x \mathrm{d}y
	\le \beta \|v(t)\|_{L^\infty} (|v(t)|_{\beta-1}^2+|\nabla{v}(t)|_{\beta-1}^2)
	\le C\varepsilon_1\beta (|v(t)|_{\beta-1}^2+|\nabla{v}(t)|_{\beta-1}^2),
\end{equation}
and
\begin{equation}\label{eq-v-beta-2}
	\int_{\mathbbm{R}}\int_{\mathbbm{R}_+} (1+x)^{\beta}|\nabla{v}|^3 \,\mathrm{d}x \mathrm{d}y
	\le \|\nabla{v}(t)\|_{L^\infty}|\nabla{v}(t)|_{\beta}^2
	\le C\varepsilon_1 |\nabla{v}(t)|_{\beta}^2.
\end{equation}
From Proposition $\ref{prop-fhb}$ and $\eqref{eq-ND-beta}$, the last two terms are estimated as
\begin{equation}\label{eq-v-beta-4}
	\int_{\mathbbm{R}}\int_{\mathbbm{R}_+} (1+x)^{\beta}\bar{u}_x|\nabla{v}|^2 \,\mathrm{d}x \mathrm{d}y
	\le C \delta |\nabla{v}(t)|_\beta^2.
\end{equation}
and
\begin{equation}\label{eq-v-beta-3}
\begin{aligned}[b]
	\int_{\mathbbm{R}}\int_{\mathbbm{R}_+} (1+x)^{\beta}|\bar{u}_{xx}||v||\nabla{v}| \,\mathrm{d}x \mathrm{d}y
	&\le C \beta \delta \| \nabla{v}(t) \|_{L^2}^2+ C\delta \| \nabla{v}(t) \|_{L^2}^2+C \delta |\nabla{v}(t)|_\beta^2.
\end{aligned}
\end{equation}
Notice that the last two terms of $\eqref{eq-v-beta-3}$ can be absorbed in the fifth term on the left of $\eqref{eq-vnabv-beta}$. Finally, the third term on the right-hand side of $\eqref{eq-vnabv-beta}$ is bounded by
\begin{equation}\label{eq-pv-deta}
	\beta \int_{\mathbbm{R}}\int_{\mathbbm{R}_+}(1+x)^{\beta-1} |p_1||v| \,\mathrm{d}x \mathrm{d}y
	\le \frac{|u_+|}{4}\beta|v(t)|_{\beta-1}^2+C \beta \int_{\mathbbm{R}}\int_{\mathbbm{R}_+} (1+x)^{\beta-1}|p_1|^2 \,\mathrm{d}x \mathrm{d}y.
\end{equation}
Substituting $\eqref{eq-v-beta-1}$-$\eqref{eq-pv-deta}$ into $\eqref{eq-vnabv-beta}$, we obtain
\begin{equation}\label{eq-vnabv-beta-jg}
	\begin{aligned}[b]
		 &\frac{\mathrm{d}}{\mathrm{d}t}(|v(t)|_\beta^2+|\nabla{v}(t)|_\beta^2)+\frac{|u_+|}{2} \beta(|v(t)|_{\beta-1}^2+|\nabla{v}(t)|_{\beta-1}^2)+|\nabla{v}(t)|_\beta^2+|u_-|\int_{\mathbbm{R}} |\nabla{v}(0,y,t)|^2 \,\mathrm{d}y\\
		 \le&  C \beta |p_1|_{\beta-1}^2+C \beta \| \nabla{v}(t) \|_{L^2}^2.
	\end{aligned}
\end{equation}
Similarly, multiplying $\eqref{eq-nabdelv-wf}$ by $2(1+x)^\beta$, $\beta\in[0,\alpha]$ and then integrating the resulting equation over $\mathbbm{R}_+\times\mathbbm{R}$, we get
\begin{equation}\label{eq-nabvdelv-beta}
	\begin{aligned}[b]
		 &\frac{\mathrm{d}}{\mathrm{d}t}(|\nabla v(t)|_\beta^2+|\Delta{v}(t)|_\beta^2)+|u_+|\beta(|\nabla v(t)|_{\beta-1}^2+|\Delta{v}(t)|_{\beta-1}^2)+2|\nabla^2{v}(t)|_\beta^2+|u_-|\int_{\mathbbm{R}} (|\nabla{v}(0,y,t)|^2+|\Delta{v}(0,y,t)|^2) \,\mathrm{d}y \\
		 \le& 2 \int_{\mathbbm{R}} {\rm{div}}{p}(0,y,t)v_x(0,y,t) \,\mathrm{d}y+ \beta \int_{\mathbbm{R}}\int_{\mathbbm{R}_+}(1+x)^{\beta-1} (|v||\nabla{v}|+|v||\Delta{v}|^2+4|v_{xy}||v_y|+2|{\rm{div}}{p}||v_x|) \,\mathrm{d}x \mathrm{d}y\\
		 & +C \int_{\mathbbm{R}}\int_{\mathbbm{R}_+} (1+x)^{\beta}(\cdots) \,\mathrm{d}x \mathrm{d}y.
	\end{aligned}
\end{equation}
 From ${\rm{div}}{p}(0,y,t)=-u_-v_x(0,y,t)$, the first term is bounded by
 $$2 \int_{\mathbbm{R}} {\rm{div}}{p}(0,y,t)v_x(0,y,t) \,\mathrm{d}y\le 2|u_-|\int_{\mathbbm{R}} v_x^2(0,y,t) \,\mathrm{d}y.$$
 The second and the third terms on the right-hand side of $\eqref{eq-nabvdelv-beta}$ can be absorbed in the third and the fourth terms on the left-hand side, respectively, owing to $\| v \|_{L^\infty}\le C \varepsilon_1$. The last term on the right-hand side can be estimated by employing the similar method as $\eqref{eq-dlnabv-1}$-$\eqref{eq-dlnabv-4}$. Consequently, we have
\begin{equation}
	C \int_{\mathbbm{R}}\int_{\mathbbm{R}_+} (1+x)^{\beta}(...) \,\mathrm{d}x \mathrm{d}y\le C \beta  \| \nabla{v}(t) \|_{L^2}^2+ C\delta \| \nabla{v}(t) \|_{L^2}^2+ C(\varepsilon_1+\delta)|\nabla v|_\beta^2+C(\varepsilon_1+\delta)|\nabla^2 v|_\beta^2.
\end{equation}
Using Cauchy inequality, the remaining terms are bounded by
\begin{equation}
	\beta \int_{\mathbbm{R}}\int_{\mathbbm{R}_+}(1+x)^{\beta-1} (4|v_{xy}||v_y|+2|{\rm{div}}{p}||v_x|) \,\mathrm{d}x \mathrm{d}y
	\le \frac{|u_+|}{4}\beta|\nabla{v}|_{\beta-1}^2+C \beta(|\nabla^2{v}|_{\beta-1}^2+|{\rm{div}}{p}|_{\beta-1}^2).
\end{equation}
Thus, for some small $\delta$ and $\varepsilon_1$, we can conclude that
\begin{equation}\label{eq-nabvdelv-beta-jg}
	\begin{aligned}[b]
		 &\frac{\mathrm{d}}{\mathrm{d}t}(|\nabla v(t)|_\beta^2+|\Delta{v}(t)|_\beta^2)+\frac{|u_+|}{2} \beta(|\nabla v(t)|_{\beta-1}^2+|\Delta{v}(t)|_{\beta-1}^2)+|\nabla^2{v}(t)|_\beta^2+|u_-|\int_{\mathbbm{R}} |\Delta{v}(0,y,t)|^2 \,\mathrm{d}y \\
		 \le&|u_-|\int_{\mathbbm{R}} |\nabla{v}(0,y,t)|^2 \,\mathrm{d}y+  C \beta  \| \nabla{v}(t) \|_{L^2}^2+ C\delta \| \nabla{v}(t) \|_{L^2}^2+C(\varepsilon_1+\delta)|\nabla v|_\beta^2+C \beta(|\nabla^2{v}|_{\beta-1}^2+|{\rm{div}}{p}|_{\beta-1}^2)\\
		 \le &|u_-|\int_{\mathbbm{R}} |\nabla{v}(0,y,t)|^2 \,\mathrm{d}y+ C \beta  \| \nabla{v}(t) \|_{L^2}^2+C(\varepsilon_1+\delta)|\nabla v|_\beta^2+C \beta(|\nabla^2{v}|_{\beta-1}^2+|{\rm{div}}{p}|_{\beta-1}^2).
	\end{aligned}
\end{equation}
We can get from $v_y\times\partial_y\eqref{eq-vp-rd}_1+p_y\cdot \partial_y\eqref{eq-vp-rd}_2+\nabla{v}_y\cdot \nabla \partial_y\eqref{eq-vp-rd}_1- \nabla{\rm{div}}{p}_y\cdot \partial_y\eqref{eq-vp-rd}_2$ that
\begin{equation}\label{eq-vnabv-y-1}
	\begin{aligned}[b]
		&\left\{\frac{1}{2}v_y^2+\frac{1}{2}|\nabla{v}_y|^2\right\}_t+\frac{1}{2} \bar{u}_x(v_y^2+|\nabla{v}_y|^2)+\bar{u}_xv_{xy}^2+|\nabla{v}_y|^2+{\rm{div}}\{p_yv_y \}\\
		&\ \ \ \ +\left\{ \frac{1}{2}vv_y^2+\frac{1}{2} \bar{u}v_y^2+\frac{1}{2} v|\nabla{v}_y|^2+\frac{1}{2} \bar{u}|\nabla{v}_y|^2 \right\}_x
		+\left\{\frac{1}{2}  g'(v+\bar{u})v_y^2+\frac{1}{2} g'(v+\bar{u})|\nabla{v}_y|^2  \right\}_y \\
		&\ \ \ \   =-\frac{1}{2} v_xv_y^2- \frac{1}{2} g''(v+\bar{u})v_y^3-\frac{1}{2} v_x|\nabla{v}_y|^2-v_y \nabla{v}_x\cdot \nabla{v}_y -v_{xy}\nabla{v}\cdot \nabla{v}_y - \bar{u}_{xx}v_yv_{xy}- \frac{3}{2} g''(v+\bar{u})v_y|\nabla{v}_y|^2\\
		& \ \ \ \ \ \ \ \ -g''(v+\bar{u})\bar{u}_xv_{xy}v_{yy}-g''(v+\bar{u})v_{yy}\nabla{v}\cdot \nabla{v}_y-g'''(v+\bar{u})v_y^2 \nabla{v}\cdot \nabla{v}_y-g'''(v+\bar{u})\bar{u}_x v_y^2v_{xy}.
	\end{aligned}
\end{equation}
The terms on right-hand side of $\eqref{eq-vnabv-y-1}$ are bounded by
\begin{equation}
	\begin{aligned}[b]
		&\| \nabla{v}(t) \|_{L^\infty} \| \nabla{v}(t) \|_{L^2}^2+ (\| \nabla{v}(t) \|_{L^\infty}+\| \bar{u}_x \|_{L^\infty}  )\| \nabla^2{v}(t) \|_{L^2}^2 +(\| \bar{u}_{xx}  \|_{L^\infty}+\| \nabla{v}(t) \|_{L^\infty} )\| \nabla{v}(t) \|_{L^2}\| \nabla^2{v}(t) \|_{L^2}\\
		\le& C(\varepsilon_1+\delta)\| \nabla^2{v}(t) \|_{L^2}^2+C(\varepsilon_1+\delta)\| \nabla{v}(t) \|_{L^2}^2.
	\end{aligned}
\end{equation}
Multiplying $\eqref{eq-vnabv-y-1}$ by $2(1+x)^\beta$, $\beta\in[0,\alpha]$ and then integrating the resulting equation over $\mathbbm{R}_+\times\mathbbm{R}$, we also obtain that
\begin{equation}\label{eq-vynabvy-beta-0}
	\begin{aligned}[b]
		 &\frac{\mathrm{d}}{\mathrm{d}t}(|v_y(t)|_\beta^2+|\nabla{v}_y(t)|_\beta^2)+|u_+|\beta(|v_y(t)|_{\beta-1}^2+|\nabla{v}_y(t)|_{\beta-1}^2)+2|\nabla{v}_y(t)|_\beta^2+|u_-|\int_{\mathbbm{R}} |\nabla{v}_y(0,y,t)|^2 \,\mathrm{d}y \\
		 \le&  \beta \left|\int_{\mathbbm{R}}\int_{\mathbbm{R}_+}(1+x)^{\beta-1} (v{v}_y+v|\nabla{v}_y|^2+2{p}_{1y}v_y) \,\mathrm{d}x \mathrm{d}y\right| +C \int_{\mathbbm{R}}\int_{\mathbbm{R}_+} (1+x)^{\beta}(\cdots) \,\mathrm{d}x \mathrm{d}y\\
		 \le & C \varepsilon_1|v_y|_{\beta-1}^2+C\varepsilon_1|\nabla v_y|_{\beta-1}^2+2\beta \left|\int_{\mathbbm{R}}\int_{\mathbbm{R}_+} (1+x)^{\beta-1}p_{1}y_{yy} \,\mathrm{d}x \mathrm{d}y\right|+C(\varepsilon_1+\delta)(|\nabla{v}|_{\beta}^2+|\nabla^2{v}|_{\beta}^2)\\
		 \le& C \varepsilon_1|v_y|_{\beta-1}^2+(C\varepsilon_1+\frac{|u_+|}{4}\beta)|\nabla v_y|_{\beta-1}^2+C \beta|p_1|_{\beta-1}^2+C(\varepsilon_1+\delta)(|\nabla{v}|_{\beta}^2+|\nabla^2{v}|_{\beta}^2).
	\end{aligned}
\end{equation}
Therefore, we may choose some small $\varepsilon_1$ and $\delta$ such that $\eqref{eq-vynabvy-beta-0}$ can be rewritten as
\begin{equation}\label{eq-vynabvy-beta-jg}
\begin{aligned}[b]
	&\frac{\mathrm{d}}{\mathrm{d}t}(|v_y(t)|_\beta^2+|\nabla{v}_y(t)|_\beta^2)+\frac{|u_+|}{2} \beta(|v_y(t)|_{\beta-1}^2+|\nabla{v}_y(t)|_{\beta-1}^2)+|\nabla{v}_y(t)|_\beta^2+|u_-|\int_{\mathbbm{R}} |\nabla{v}_y(0,y,t)|^2 \,\mathrm{d}y \\
	\le& C \beta|p_1|_{\beta-1}^2+C(\varepsilon_1+\delta)(|\nabla{v}|_{\beta}^2+|\nabla^2{v}|_{\beta}^2).	
\end{aligned}
\end{equation}
By utilizing Lemma $\ref{lem-nab2v-dengjia}$, we can obtain from $2\times\eqref{eq-vnabv-beta-jg}+\eqref{eq-nabvdelv-beta-jg}+\eqref{eq-vynabvy-beta-jg}$ that
\begin{equation}\label{eq-vH2-sjsj}
	\begin{aligned}[b]
	&\frac{\mathrm{d}}{\mathrm{d}t}(|v(t)|_\beta^2+|\nabla{v}(t)|_\beta^2+|\nabla^2{v}(t)|_\beta^2)+|u_+| \beta(|v(t)|_{\beta-1}^2+|\nabla{v}(t)|_{\beta-1}^2+|\nabla^2{v}(t)|_{\beta-1}^2)+|\nabla{v}(t)|_\beta^2+|\nabla^2{v}(t)|_\beta^2\\
	& \ \ \ \ \ \ \ \ \ \ \ \ \ \ \ \ +|u_-|\int_{\mathbbm{R}} (|\nabla{v}(0,y,t)|^2+|\nabla^2{v}(0,y,t)|^2) \,\mathrm{d}y \\
	\le& C \beta\left(|p_1(t)|_{\beta-1}^2+|{\rm{div}}{p}(t)|_{\beta-1}^2+|\nabla^2 v(t)|_{\beta-1}^2+\| \nabla{v}(t) \|_{L^2}^2 \right).
	\end{aligned}
\end{equation}
On the other hand, $\eqref{eq-vp-rd}_2$ gives
\begin{equation}\label{eq-p-beta-jg}
	\begin{aligned}[b]
		& \int_{\mathbbm{R}}\int_{\mathbbm{R}_+} (1+x)^{\beta}(|\nabla{\rm{div}}{p}|^2+2({\rm{div}}{p})^2+|p|^2) \,\mathrm{d}x \mathrm{d}y\\
		=& \int_{\mathbbm{R}}\int_{\mathbbm{R}_+} (1+x)^{\beta}|\nabla{v}|^2 \,\mathrm{d}x \mathrm{d}y
		+2 \int_{\mathbbm{R}}\int_{\mathbbm{R}_+} \{(1+x)^{\beta}({\rm{div}}p)p_1\}_x \,\mathrm{d}x \mathrm{d}y
		-2 \beta\int_{\mathbbm{R}}\int_{\mathbbm{R}_+}  (1+x)^{\beta-1}({\rm{div}}p)p_1 \,\mathrm{d}x \mathrm{d}y \\
		=&|\nabla{v}(t)|_{\beta}^2-2 \int_{\mathbbm{R}}({\rm{div}}p(0,y,t))p_1(0,y,t) \,\mathrm{d}y -2\beta\int_{\mathbbm{R}}\int_{\mathbbm{R}_+}  (1+x)^{\beta-1}({\rm{div}}p)p_1 \,\mathrm{d}x \mathrm{d}y\\
		\le& |\nabla{v}(t)|_\beta^2-2u_-\int_{\mathbbm{R}}  v_x(0,y,t)v_{xt}(0,y,t)\,\mathrm{d}y+C_{u_-}^1|u_-|\int_{\mathbbm{R}} |\nabla{v}(0,y,t)|^2 \,\mathrm{d}y +C_{u_-}^2|u_-|\int_{\mathbbm{R}} (\Delta{v}(0,y,t))^2 \,\mathrm{d}y\\
		&\ \  +C\beta(|{\rm{div}}{p}(t)|_{\beta-1}^2+|p_1(t)|_{\beta-1}^2),
	\end{aligned}
\end{equation}
where $C_{u_-}^1=C(\| v_x \|_{L^\infty}+\| \bar{u}_x \|_{L^\infty}+|u_-|+1)$ and $C_{u_-}^2=C|u_-|$ are fixed constants. Here, we have used the equations $\eqref{eq-vx0-bj}$ and $\eqref{eq-divpx0}$ which yields
\begin{equation}
	\begin{aligned}[b]
		-2\int_{\mathbbm{R}}({\rm{div}}p(0,y,t))p_1(0,y,t) \,\mathrm{d}y
		=&2u_-\int_{\mathbbm{R}}v_x(0,y,t)({\rm{div}}{p}_x(0,y,t)-v_x(0,y,t)) \,\mathrm{d}y\\
		\le & -2u_-\int_{\mathbbm{R}} v_x(0,y,t)v_{xt}(0,y,t)\,\mathrm{d}y+Cu_-^2\int_{\mathbbm{R}} (\Delta{v}(0,y,t))^2 \,\mathrm{d}y\\
		&\ \ +C|u_-|(\| v_x \|_{L^\infty}+\| \bar{u}_x \|_{L^\infty}+|u_-|+1)\int_{\mathbbm{R}} |\nabla{v}(0,y,t)|^2 \,\mathrm{d}y .
	\end{aligned}
\end{equation}
Finally, let $\lambda_1>\max\left\{C_{u_-}^1,C_{u_-}^2,1\right\}$. We can get from $\lambda_1\times\eqref{eq-vH2-sjsj}+\eqref{eq-p-beta-jg}$ that
\begin{equation}\label{eq-lam1lam2}
	\begin{aligned}[b]
		&\frac{\mathrm{d}}{\mathrm{d}t}\left(\lambda_1|v(t)|_\beta^2+\lambda_1|\nabla{v}(t)|_\beta^2+\lambda_1|\nabla^2{v}(t)|_\beta^2\right) +\lambda_1|u_+|\beta(|v(t)|_{\beta-1}^2+|\nabla{v}(t)|_{\beta-1}^2+|\nabla^2{v}(t)|_{\beta-1}^2)\\
		&\ \ +(\lambda_1-1)|\nabla{v}(t)|_\beta^2+\lambda_1|\nabla^2{v}(t)|_\beta^2+|\nabla{\rm{div}}{p}(t)|_{\beta}^2+2|{\rm{div}}{p}(t)|_{\beta}^2+|p(t)|_{\beta}^2\\
		&\ \ +(\lambda_1-C_{u_-}^1)|u_-|\int_{\mathbbm{R}} |\nabla{v}(0,y,t)|^2 \,\mathrm{d}y+(\lambda_1-C_{u_-}^2)|u_-|\int_{\mathbbm{R}} |\nabla^2{v}(0,y,t)|^2 \,\mathrm{d}y\\
		\le&  C \beta(|\nabla^2{v}(t)|_{\beta-1}^2+|p(t)|_{\beta-1}^2+|{\rm{div}}{p}(t)|_{\beta-1}^2+\| \nabla{v}(t) \|_{L^2}^2 )
		-2u_-\int_{\mathbbm{R}}  v_x(0,y,t)v_{xt}(0,y,t)\,\mathrm{d}y.
	\end{aligned}
\end{equation}
We note that there exist sufficiently large constants $R_1$, $R_2$ and $R_3$ satisfying $\frac{C\beta}{1+R_i}\le \frac{1}{4}, i=1,2,3$, such that
\begin{equation}\label{eq-beta-1-1}
\begin{aligned}[b]
	& C \beta\int_{\mathbbm{R}}\int_{\mathbbm{R}_+} (1+x)^{\beta-1}|{\rm{div}}p|^2 \,\mathrm{d}x \mathrm{d}y  \\
	\le& C_{ R_1}\beta \int_{\mathbbm{R}}\int_{0}^{R_1} |{\rm{div}}p|^2 \,\mathrm{d}x \mathrm{d}y +\frac{C\beta}{1+R_1}\int_{\mathbbm{R}}\int_{R_1}^{\infty} (1+x)^{\beta} |{\rm{div}}p|^2 \,\mathrm{d}x \mathrm{d}y \,\mathrm{d}x \mathrm{d}y\\
    \le&  C_{ R_1}\beta \int_{\mathbbm{R}}\int_{0}^{\infty} |{\rm{div}}p|^2 \,\mathrm{d}x \mathrm{d}y+\frac{1}{4}\int_{\mathbbm{R}}\int_{0}^{\infty} (1+x)^{\beta} |{\rm{div}}p|^2 \,\mathrm{d}x \mathrm{d}y \,\mathrm{d}x \mathrm{d}y,
\end{aligned}
\end{equation}
\begin{equation}\label{eq-beta-1-2}
	 C \beta\int_{\mathbbm{R}}\int_{\mathbbm{R}_+} (1+x)^{\beta-1}|p_1|^2 \,\mathrm{d}x \mathrm{d}y \le  C_{R_2}\beta \| p(t) \|_{L^2}^2  +\frac{1}{4}|p(t)|_{\beta}^2,
\end{equation}
and
\begin{equation}\label{eq-beta-1-3}
	C\beta \int_{\mathbbm{R}}\int_{\mathbbm{R}_+} (1+x)^{\beta-1}|\nabla^2{v}|^2 \,\mathrm{d}x \mathrm{d}y
	\le C_{R_3}\beta \|\nabla^2{v}(t)\|_{L^2}^2   +\frac{1}{4}|\nabla^2{v}(t)|_{\beta}^2.
\end{equation}
Substituting $\eqref{eq-beta-1-1}$-$\eqref{eq-beta-1-3}$ into $\eqref{eq-lam1lam2}$, and then multiplying $\eqref{eq-lam1lam2}$ by $(1+t)^\gamma$, $\gamma\in [0,\alpha]$, and integrating the resulting inequality over $[0,t]$, we get
\begin{equation}\label{eq-lam1lam2-jit}
	\begin{aligned}[b]
        &(1+t)^\gamma\left(\lambda_1|v(t)|_\beta^2+\lambda_1|\nabla{v}(t)|_\beta^2+\lambda_1|\nabla^2{v}(t)|_\beta^2\right) +\lambda_1|u_+|\beta \int_0^t (1+\tau)^\gamma(|v(t)|_{\beta-1}^2+|\nabla{v}(\tau)|_{\beta-1}^2+|\nabla^2{v}(\tau)|_{\beta-1}^2) \,\mathrm{d}\tau \\
		&\ \ +\int_0^t (1+\tau)^\gamma \left\{(\lambda_1-1)|\nabla{v}(\tau)|_\beta^2+(\lambda_1- \frac{1}{4})|\nabla^2{v}(\tau)|_\beta^2+\frac{3}{4}|\nabla{\rm{div}}{p}(\tau)|_{\beta}^2+\frac{7}{4}|{\rm{div}}{p}(\tau)|_{\beta}^2+|p(\tau)|_{\beta}^2\right\} \,\mathrm{d}\tau\\
		&\ \ +\int_0^t (1+\tau)^\gamma \left\{(\lambda_1-C_{u_-}^1)|u_-|\int_{\mathbbm{R}} |\nabla{v}(0,y,\tau)|^2 \,\mathrm{d}y+(\lambda_1-C_{u_-}^2)|u_-|\int_{\mathbbm{R}} |\nabla^2{v}(0,y,\tau)|^2 \,\mathrm{d}y \right\} \,\mathrm{d}\tau\\
		\le&  C|v_0|_{\beta,2}^2+C \gamma \int_0^t (1+\tau)^{\gamma-1} (|v(\tau)|_\beta^2+|\nabla{v}(\tau)|_{\beta}^2+|\nabla^2{v}(\tau)|_\beta^2) \,\mathrm{d}\tau+C\beta \int_0^t (1+\tau)^\gamma (\| p (\tau)\|_{L^2}^2+\| {\rm{div}}{p}(\tau) \|_{L^2}^2 ) \,\mathrm{d}\tau\\
		&+C\beta \int_0^t (1+\tau)^\gamma (\| \nabla^2{v}(\tau) \|_{L^2}^2+\| \nabla{v}(\tau) \|_{L^2}^2 ) \,\mathrm{d}\tau-2u_-\int_0^t (1+\tau)^\gamma \int_{\mathbbm{R}}  v_x(0,y,t)v_{xt}(0,y,t)\,\mathrm{d}y \mathrm{d}\tau.
	\end{aligned}
\end{equation}
Using Cauchy inequality and Sobolev inequality $\eqref{eq-h-Rzheng}$, the last term of $\eqref{eq-lam1lam2-jit}$ is bounded by
\begin{equation}\label{eq-lam1lam2-jit-1}
	\begin{aligned}[b]
		&\int_0^t (1+\tau)^\gamma \int_{\mathbbm{R}}  v_x(0,y,\tau)v_{xt}(0,y,\tau)\,\mathrm{d}y \mathrm{d}\tau\\
		=&\frac{1}{2} (1+t)^\gamma\int_{\mathbbm{R}} v_x^2(0,y,t) \,\mathrm{d}y -\frac{1}{2} \int_{\mathbbm{R}} v_x^2(0,y,0) \,\mathrm{d}y
		-\frac{1}{2} \gamma \int_0^t (1+\tau)^{\gamma-1}\int_{\mathbbm{R}}  v_x^2(0,y,\tau) \,\mathrm{d}y  \,\mathrm{d}\tau \\
		\le& \frac{1}{2}(1+t)^\gamma \int_{\mathbbm{R}} \| v_x(\cdot,y,t) \|_{L^\infty}^2  \,\mathrm{d}y+\frac{1}{2}\int_{\mathbbm{R}} \| v_x(\cdot,y,0) \|_{L^\infty}^2  \,\mathrm{d}y+\gamma \int_0^t (1+\tau)^{\gamma-1}\int_{\mathbbm{R}} \| v_x^2(\cdot,y,\tau) \|_{L^\infty}^2 \,\mathrm{d}y \mathrm{d}\tau \\
		\le&  \frac{1}{2} (1+t)^\gamma\| v_x(t) \|_{L^2}^2 +\frac{1}{2} (1+t)^\gamma\| v_{xx}(t) \|_{L^2}^2+\| v_0 \|_{H^2}^2+C\gamma \int_0^t (1+\tau)^{\gamma-1}( \| v_x(\tau) \|_{L^2}^2 +\| v_{xx}(\tau) \|_{L^2}^2) \,\mathrm{d}\tau.
	\end{aligned}
\end{equation}
Notice that if $\lambda_1$ is sufficiently large such that $\lambda_1>\max\left\{C_{u_-}^1,C_{u_-}^2,1\right\}$, the first two terms on the right-hand side of $\eqref{eq-lam1lam2-jit-1}$ can be absorbed in $(1+t)^\gamma|\nabla{v}(t)|_\beta^2$ and $(1+t)^\gamma|\nabla^2{v}(t)|_\beta^2$ on the left-hand side of $\eqref{eq-lam1lam2-jit}$, respectively. And the last term of $\eqref{eq-lam1lam2-jit-1}$ can be absorbed in $C \gamma \int_0^t (1+\tau)^{\gamma-1} (|\nabla{v}(\tau)|_{\beta}^2+|\nabla^2{v}(\tau)|_\beta^2) \,\mathrm{d}\tau$, which appears on the right-hand side of $\eqref{eq-lam1lam2-jit}$. Consequently, the desired inequality $\eqref{eq-vsjl-alpha}$ can be obtained.
$\hfill\Box$

Applying the arguments of Kawashima and Matsumura in \cite{Kawashima1985}, we can obtain the following Lemma \ref{lem-alpha-gamma}. For the sake of the completeness of the present paper, we give the main idea of proving the lemma \ref{lem-alpha-gamma}. The details can be referred to \cite{Kawashima1985}.
\begin{lem}\label{lem-alpha-gamma}
Let $\gamma\in[0,\alpha]\cap\mathbbm{Z}$. There are positive constants $\varepsilon_2(\le\varepsilon_1)$ and $C=C(\varepsilon_2)$, which are independent of $T$ and $\gamma$ such that if $ \sup_{0\le \tau\le T}\|v(\tau)\|_{H^3}+\delta\le\varepsilon_2$, then
\begin{equation}\label{eq-vsjl-bate}
\begin{aligned}[b]
&(1+t)^\gamma |v(t)|_{\alpha-\gamma,2}^2+(\alpha-\gamma)\int_0^t \ (1+\tau)^\gamma |v(\tau)|_{\alpha-\gamma-1,2}^2 \,\mathrm{d}\tau
 +\int_0^t (1+\tau)^\gamma \int_{\mathbbm{R}} (|\nabla{v}(0,y,\tau)|^2+|\nabla^2{v}(0,y,\tau)|^2) \,\mathrm{d}y  \,\mathrm{d}\tau\\
 & \ \ +\int_0^t (1+\tau)^\gamma (|\nabla v(\tau)|_{\alpha-\gamma,1}^2+|\nabla{\rm{div}}{p}(\tau)|_{\alpha-\gamma}^2+|{\rm{div}}{p}(\tau)|_{\alpha-\gamma}^2+|p(\tau)|_{\alpha-\gamma}^2)\,\mathrm{d}\tau \le C M_\alpha^2,
\end{aligned}
\end{equation}
holds for $t\in[0,T]$. Moreover, for any $0\le\gamma\le[\alpha]$, the following estimate holds:
\begin{equation}\label{eq-vsjl-gamma}
\begin{aligned}[b]
(1+t)^\gamma \|v(t)\|_{H^2}^2+\int_0^t (1+\tau)^\gamma (\|\nabla{v} (\tau)\|_{H^1}^2+\|\nabla{\rm{div}}{p}(\tau) \|_{L^2}^2+\|{\rm{div}}{p}(\tau) \|_{L^2}^2+\| p(\tau)\|_{L^2}^2) \,\mathrm{d}\tau &\\
+\int_0^t (1+\tau)^\gamma \int_{\mathbbm{R}} (|\nabla{v}(0,y,\tau)|^2+|\nabla^2{v}(0,y,\tau)|^2) \,\mathrm{d}y  \,\mathrm{d}\tau&\le C M_\alpha^2.
\end{aligned}
\end{equation}
\end{lem}
{\it\bfseries Proof.}
Letting $\beta=0$ in \eqref{eq-vsjl-alpha}, we have
\begin{equation}\label{eq-vbeta1}
\begin{aligned}[b]
&(1+t)^\gamma \|v(t)\|_{H^2}^2+\int_0^t \ (1+\tau)^\gamma (\|\nabla{v} (\tau)\|_{H^1}^2+\|\nabla{\rm{div}}{p}(\tau) \|_{L^2}^2+\|{\rm{div}}{p}(\tau) \|_{L^2}^2+\| p(\tau)\|_{L^2}^2) \,d{\tau}\\
&\ \ \ \ \ \ \ \ \ \ \ \ \ \ \ \  +\int_0^t (1+\tau)^\gamma \int_{\mathbbm{R}} (|\nabla{v}(0,y,\tau)|^2+|\nabla^2{v}(0,y,\tau)|^2) \,\mathrm{d}y  \,\mathrm{d}\tau\\
\le& C \left(\|v_0\|_{H^2}^2+\gamma \int_0^t \ (1+\tau)^{\gamma-1}\|v(\tau)\|_{H^2}^2 \,d{\tau} \right).
\end{aligned}
\end{equation}

$\mathrm{(i)}$ Firstly, we consider the case of $\gamma=0$, $\beta=\alpha$ in $\eqref{eq-vsjl-alpha}$. Combining with \eqref{eq-vsjl2}, for any $\alpha\ge0$, we get
\begin{equation}\label{eq-vbeta2}
\begin{aligned}[b]
&|v(t)|_{\alpha,2}^2+ \alpha\int_0^t \  |v(\tau)|_{\alpha-1,2}^2  \,d{\tau} +\int_0^t \  ( | \nabla{v}(\tau) |_{\alpha,1}^2 +|\nabla{\rm{div}}{p}(\tau)  |_{\alpha}^2+| {\rm{div}}{p}(\tau) |_{\alpha}^2+| p(\tau) |_{\alpha}^2)  \,d{\tau}\\
&\ \ \ \ \ \ \ \ \ +\int_0^t  \int_{\mathbbm{R}} (|\nabla{v}(0,y,\tau)|^2+|\nabla^2{v}(0,y,\tau)|^2) \,\mathrm{d}y  \,\mathrm{d}\tau\\
\le& C \left(|v_0|_{\alpha,2}^2+\alpha\int_0^t \ (\|{\rm{div}}{p}(\tau) \|_{L^2}^2+\| p(\tau)\|_{L^2}^2+\|\nabla{v} (\tau)\|_{H^1}^2) \,d{\tau} \right)\\
\le& CM_\alpha^2,
\end{aligned}
\end{equation}
where $M_\alpha^2:=\|v_0\|_{H^3}^2+|v_0|_{\alpha,2}^2$ is defined in $\eqref{eq-M0-Mal}$. Therefore we finished the proof for $0\le\alpha<1$.

$\mathrm{(ii)}$ Secondly, letting $\gamma=1$ in \eqref{eq-vbeta1}, and using \eqref{eq-vbeta2}, we have
\begin{equation}\label{eq-vbeta3}
\begin{aligned}[b]
&(1+t) \|v(t)\|_{H^2}^2+\int_0^t \ (1+\tau)(\|\nabla{v} (\tau)\|_{H^1}^2+\|\nabla{\rm{div}}{p}(\tau) \|_{L^2}^2+\|{\rm{div}}{p}(\tau) \|_{L^2}^2+\| p(\tau)\|_{L^2}^2) \,d{\tau}\\
&\ \ \ \ \ \ \ \ \ \ \ \ \ \ \ \ \ \ +\int_0^t (1+\tau)\int_{\mathbbm{R}} (|\nabla{v}(0,y,\tau)|^2+|\nabla^2{v}(0,y,\tau)|^2) \,\mathrm{d}y  \,\mathrm{d}\tau\\
\le& C \left(\|v_0\|_{H^2}^2+ \int_0^t \ \|v(\tau)\|_{H^2}^2 \,d{\tau} \right)\\
\le& CM_\alpha^2,
\end{aligned}
\end{equation}
where $\alpha\ge1$. That is, $\eqref{eq-vsjl-gamma}$ with $\gamma=1$ is proved. Consider the case of $\gamma=1$, $\beta=\alpha-1$ in $\eqref{eq-vsjl-alpha}$. According to the results of $\eqref{eq-vbeta2}$ and $\eqref{eq-vbeta3}$, for any $\alpha\ge1$, we have
\begin{equation}\label{eq-vbeta4}
\begin{aligned}[b]
&(1+t) |v(t)|_{\alpha-1,2}^2+ (\alpha-1)\int_0^t \ (1+\tau)|v(\tau)|_{\alpha-2,2}^2 \,d{\tau} +\int_0^t (1+\tau) \int_{\mathbbm{R}} (|\nabla{v}(0,y,\tau)|^2+|\nabla^2{v}(0,y,\tau)|^2) \,\mathrm{d}y  \,\mathrm{d}\tau\\
&+\int_0^t \ (1 +\tau) (| \nabla{v}(\tau)  |_{\alpha-1,1}^2 +| \nabla{\rm{div}}{p}(\tau)  |_{\alpha-1}^2+| {\rm{div}}{p}(\tau)  |_{\alpha-1}^2+|p (\tau)  |_{\alpha-1}^2 )\,d{\tau}\\
\le& C \left(|v_0|_{\alpha-1,2}^2+ \int_0^t \ |v(\tau)|_{\alpha-1,2}^2 \,d{\tau} +(\alpha-1)\int_0^t \ (1+\tau)(\|{\rm{div}}{p}(\tau) \|_{L^2}^2+\| p(\tau)\|_{L^2}^2+\|\nabla{v} (\tau)\|_{H^1}^2)\,d{\tau} \right)\\
\le& CM_\alpha^2.
\end{aligned}
\end{equation}
Consequently, we finished the proof for $\alpha<2$.

$\mathrm{(iii)}$ Just as Step $2$, letting $\gamma=2$ in \eqref{eq-vbeta1}, and using \eqref{eq-vbeta4}, we have
\begin{equation}\label{eq-vbeta5}
\begin{aligned}[b]
&(1+t)^2 \|v(t)\|_{H^2}^2+\int_0^t \ (1+\tau)^2 (\|\nabla{v} (\tau)\|_{H^1}^2+\|\nabla{\rm{div}}{p}(\tau) \|_{L^2}^2+\|{\rm{div}}{p}(\tau) \|_{L^2}^2+\| p(\tau)\|_{L^2}^2) \,d{\tau}\\
&\ \ \ \ \ \ \ \ \ \ \ \ \ \ \ \ \ \ +\int_0^t (1+\tau)^2\int_{\mathbbm{R}} (|\nabla{v}(0,y,\tau)|^2+|\nabla^2{v}(0,y,\tau)|^2) \,\mathrm{d}y  \,\mathrm{d}\tau\\
\le& C \left(\|v_0\|_{H^2}^2+2 \int_0^t \ (1+\tau)\|v(\tau)\|_{H^2}^2 \,d{\tau} \right)\\
\le& CM_\alpha^2,
\end{aligned}
\end{equation}
where $\alpha\ge2$ is assumed.  The inequality $\eqref{eq-vsjl-gamma}$ with $\gamma=2$ is obtained. We reflected on the case of $\gamma=2$, $\beta=\alpha-2$ in \eqref{eq-vsjl-alpha}, combining with \eqref{eq-vbeta4} and \eqref{eq-vbeta5}, for any $\alpha\ge2$, we obtain
\begin{equation}\label{eq-vbeta6}
\begin{aligned}
 &(1+t)^2|v(t)|_{\alpha-2,2}^2+ (\alpha-2)\int_0^t \ (1+\tau)^2 |v(\tau)|_{\alpha-3,2}^2  \,d{\tau} +\int_0^t (1\!+\!\tau)^2 \int_{\mathbbm{R}} (|\nabla{v}(0,y,\tau)|^2+|\nabla^2{v}(0,y,\tau)|^2) \,\mathrm{d}y  \,\mathrm{d}\tau \\
 &+\int_0^t \ (1+\tau)^2 (| \nabla{v}(\tau)  |_{\alpha-2,1}^2 +| \nabla{\rm{div}}{p}(\tau)  |_{\alpha-2}^2+| {\rm{div}}{p}(\tau)  |_{\alpha-2}^2+|p (\tau)  |_{\alpha-2}^2 )\,d{\tau}\\
\le& C \left(|v_0|_{\alpha-2,2}^2+2 \int_0^t \ (1+\tau)|v(\tau)|_{\alpha-2,2}^2 \,d{\tau} +(\alpha-2)\int_0^t \ (1+\tau)^2 (\|{\rm{div}}{p}(\tau) \|_{L^2}^2+\| p(\tau)\|_{L^2}^2+\|\nabla{v} (\tau)\|_{H^1}^2) \,d{\tau} \right)\\
\le& C M_\alpha^2.
\end{aligned}
\end{equation}
Thus, the lemma is proved for $\alpha<3$. \\
\indent Repeating the above procedure, we can proof that the equation \eqref{eq-vsjl-bate} holds for any $\alpha\ge0$.
$\hfill\Box$

Referring to \cite{Nishikawa1998} and \cite{Nishikawa2007}, we get the following Lemma $\ref{lem-v-al-integer}$.
\begin{lem}\label{lem-v-al-integer}
	Under the same assumptions as in Lemma $\ref{lem-alpha-gamma}$. then the following estimate holds. For any $\varepsilon>0$
	\begin{equation}\label{eq-integer-al}
		\begin{aligned}[b]
			(1+t)^{\alpha+\varepsilon} &\|v(t)\|_{H^2}^2 +\int_0^t (1+\tau)^{\alpha+\varepsilon} (\|\nabla{v}(\tau) \|_{H^1}^2 \!+\! \|\nabla{\rm{div}}{p}(\tau) \|_{L^2}^2 \!+\! \|{\rm{div}}{p}(\tau) \|_{L^2}^2 \!+\! \| p(\tau)\|_{L^2}^2) \,\mathrm{d}\tau \\
			&+\int_0^t (1+\tau)^{\alpha+\varepsilon}\int_{\mathbbm{R}} (|\nabla{v}(0,y,\tau)|^2+|\nabla^2{v}(0,y,\tau)|^2) \,\mathrm{d}y  \,\mathrm{d}\tau\le C(1+t)^\varepsilon M_\alpha^2.
		\end{aligned}
	\end{equation}
\end{lem}
{\it\bfseries Proof.}
If $\alpha$ is an integer, we take $\gamma=[\alpha]=\alpha$ in $\eqref{eq-vsjl-gamma}$, which yields $\eqref{eq-integer-al}$ in the case of integer $\alpha$. Next, we consider the case of non-integer $\alpha$. Recalling $\eqref{eq-vbeta1}$, we just need to show the boundedness of the last term
$$\gamma \int_0^t \ (1+\tau)^{\gamma-1}\|v(\tau)\|_{H^2}^2 \,d{\tau}. $$
 By $\mathrm{H\ddot{o}lder}$ inequality, it holds that
\begin{equation*}
    \|\nabla^k{v}(t)\|_{L^2}^2=\int_{\mathbbm{R}}\int_{\mathbbm{R}_+} \left\{(1+x)^{1/p}\nabla^k{v}^2\right\}^{1/q} \left\{(1+x)^{-1/q}\nabla^k{v}^2\right\}^{1/p}  \,\mathrm{d}x \mathrm{d}y
	\le |\nabla^k{v}(t)|_{1/p}^{2/q}|\nabla^k{v}(t)|_{-1/q}^{2/p},\ \ k=0,1,2,
\end{equation*}
where $1/p+1/q=1$, $p:=1/(\alpha-[\alpha])$, $q:=1/(1- \alpha+[\alpha])$. Taking $\gamma=[\alpha]$ in $\eqref{eq-vsjl-bate}$, we obtain
\begin{equation*}
	\begin{aligned}[b]
		&(1+t)^{[\alpha]} |v(t)|_{\alpha-[\alpha],2}^2+(\alpha-[\alpha])\int_0^t \ (1+\tau)^{[\alpha]} |v(\tau)|_{\alpha-[\alpha]-1,2}^2 \,\mathrm{d}\tau\\
        & \ \ \ \ \ +\int_0^t (1+\tau)^{[\alpha]} (|\nabla v(\tau)|_{\alpha-[\alpha]}^2+|\nabla{\rm{div}}{p}(\tau)|_{\alpha-[\alpha]}^2+|{\rm{div}}{p}(\tau)|_{\alpha-[\alpha]}^2+|p(\tau)|_{\alpha-[\alpha]}^2)\,\mathrm{d}\tau\\
        &\ \ \ \ \ +\int_0^t (1+\tau)^{[\alpha]} \int_{\mathbbm{R}} (|\nabla{v}(0,y,\tau)|^2+|\nabla^2{v}(0,y,\tau)|^2) \,\mathrm{d}y  \,\mathrm{d}\tau \le C M_\alpha^2,
	\end{aligned}
\end{equation*}
 Using $\mathrm{H\ddot{o}lder}$ inequality again, we get
\begin{equation}\label{eq-v-pqgamma}
	\begin{aligned}[b]
		\int_0^t \ (1+\tau)^{\gamma-1} \|v(\tau)\|_{L^2}^2 \,d{\tau}
		\le& \int_0^t \ (1+\tau)^{\gamma-1}|v(\tau)|_{1/p}^{2/q}|v(\tau)|_{-1/q}^{2/p}  \,d{\tau}\\
		\le& \int_0^t (1+\tau)^{\gamma-1-[\alpha]}\left\{(1+\tau)^{[\alpha]}|v(\tau)|_{1/p}^2\right\}^{1/q}\left\{(1+\tau)^{[\alpha]}|v(\tau)|_{-1/q}^2\right\}^{1/p} \,\mathrm{d}\tau\\
		\le& CM_\alpha^2 \int_0^t (1+\tau)^{\gamma-1- [\alpha]}\left\{(1+\tau)^{[\alpha]}|v(\tau)|_{-1/q}^2\right\}^{1/p} \,\mathrm{d}\tau  \\
		\le& CM_\alpha^{2/q}\left(\int_0^t (1+\tau)^{(\gamma-1- [\alpha])q} \,\mathrm{d}\tau\right)^{1/q}\left(\int_0^t(1+\tau)^{[\alpha]}|v(\tau)|_{-1/q}^{2/q}\,\mathrm{d}\tau\right)^{1/p}\\
		\le& CM_\alpha^2\left(\int_0^t (1+\tau)^{(\gamma-1- [\alpha])q} \,\mathrm{d}\tau\right)^{1/q}.
	\end{aligned}
\end{equation}
Similarly,
\begin{equation}\label{eq-nabv-pqgamma}
	\int_0^t \ (1+\tau)^{\gamma-1} \|\nabla{v}(\tau)\|_{H^1}^2 \,d{\tau}
	\le CM_\alpha^2\left(\int_0^t (1+\tau)^{(\gamma-1- [\alpha])q} \,\mathrm{d}\tau\right)^{1/q}.
\end{equation}
Taking $\gamma=\alpha+ \varepsilon (\varepsilon>0)$ in $\eqref{eq-v-pqgamma}$ and $\eqref{eq-nabv-pqgamma}$, we obtain
\begin{equation}\label{eq-vnabv-gamma}
	\int_0^t \ (1+\tau)^{\alpha+\varepsilon-1} \|v(\tau)\|_{H^2}^2 \,d{\tau}\le C_\varepsilon(1+t)^\varepsilon M_\alpha^2.
\end{equation}
Therefore, substituting $\eqref{eq-vnabv-gamma}$ into $\eqref{eq-vbeta1}$, we get
\begin{equation}\label{eq-vp-al-eps}
\begin{aligned}[b]
	(1+t)^{\alpha+\varepsilon}&\|v(t)\|_{H^2}^2+\int_0^t (1+\tau)^{\alpha+\varepsilon}(\|\nabla{v}(\tau) \|_{H^1}^2 + \|\nabla{\rm{div}}{p}(\tau) \|_{L^2}^2 + \|{\rm{div}}{p}(\tau) \|_{L^2}^2 +\| p(\tau)\|_{L^2}^2) \,\mathrm{d}\tau\\
	&+\int_0^t (1+\tau)^{\alpha+\varepsilon}\int_{\mathbbm{R}} (|\nabla{v}(0,y,\tau)|^2+|\nabla^2{v}(0,y,\tau)|^2) \,\mathrm{d}y  \,\mathrm{d}\tau \le C_\varepsilon(1+t)^\varepsilon M_\alpha^2.	
\end{aligned}
\end{equation}
Consequently, the estimate $\eqref{eq-integer-al}$ holds for the non-integer $\alpha$. The proof of Lemma $\ref{lem-v-al-integer}$ is completed.
$\hfill\Box$

Next, we obtain the estimate for the third order derivative of $v$ by employing the time weighted energy method.
\subsection{Time Weighted Energy Estimates}
\begin{lem}\label{lem-dlv-al-eps}
	Under the same assumptions as in Lemma $\ref{lem-alpha-gamma}$, the estimate holds that
	\begin{equation}\label{eq-dlv-al-eps}
		\begin{aligned}[b]
			(1+t)^{\alpha+\varepsilon} &\|\nabla^2{v}(t)\|_{H^1}^2 +\int_0^t (1+\tau)^{\alpha+\varepsilon} (\|\nabla^3{v}(\tau) \|_{L^2}^2 \!+\! \|\nabla^2{\rm{div}}{p}(\tau) \|_{L^2}^2 \!+\! \|\nabla{\rm{div}}{p}(\tau) \|_{L^2}^2 \!+\! \|\nabla p(\tau)\|_{L^2}^2) \,\mathrm{d}\tau \\
			&+\int_0^t (1+\tau)^{\alpha+\varepsilon}\int_{\mathbbm{R}} |\nabla^3{v}(0,y,\tau)|^2 \,\mathrm{d}y  \,\mathrm{d}\tau\le C(1+t)^\varepsilon M_\alpha^2.
		\end{aligned}
	\end{equation}	
\end{lem}
 {\it\bfseries Proof.}
Multiplying $\eqref{eq-vyy-nabvyy-djf}$ by $(1+t)^\gamma$ and then integrating the resulting inequality over $[0,t]$, we get
\begin{equation}\label{eq-vyy-H1-cjf}
	\begin{aligned}[b]
		&(1+t)^\gamma\| v_{yy}(t) \|_{H^1}^2+\int_0^t (1+\tau)^\gamma (\| \sqrt{\bar{u}_x}v_{yy}(\tau) \|_{L^2}^2+\| \sqrt{\bar{u}_x}\nabla v_{yy}(\tau) \|_{L^2}^2+\| \nabla{v}_{yy}(t) \|_{L^2}^2  ) \,\mathrm{d}\tau\\
		&\ \ \ \ \ \ \ +\int_0^t (1+\tau)^\gamma \int_{\mathbbm{R}} |\nabla{v}_{yy}(0,y,\tau)|^2 \,\mathrm{d}y  \,\mathrm{d}\tau\\
		&\le CM_0^2+C\gamma \int_0^t (1+\tau)^{\gamma-1}\| v_{yy}(\tau) \|_{H^1}^2 \,\mathrm{d}\tau+C(\varepsilon_2+M_0^4)\int_0^t (1+\tau)^\gamma\| \nabla\Delta v(\tau) \|_{L^2}^2  \,\mathrm{d}\tau\\
		&\ \ +C \int_0^t (1+\tau)^\gamma(\| \nabla{v}(\tau) \|_{L^2}^2+\| \nabla^2{v}(\tau) \|_{L^2}^2 ) \,\mathrm{d}\tau.
	\end{aligned}
\end{equation}
At the same time, multiplying $\eqref{eq-vH3-djf}$ by $(1+t)^\gamma$ and integrating the resulting inequality over $[0,t]$, we get
\begin{equation}\label{eq-delv-H1-cjf}
	\begin{aligned}[b]
		&(1+t)^\gamma\| \Delta v(t) \|_{H^1}^2+\int_0^t (1+\tau)^\gamma (\| \sqrt{\bar{u}_x}\Delta v(\tau) \|_{L^2}^2+\| \sqrt{\bar{u}_x}\nabla\Delta v(\tau) \|_{L^2}^2+\| \nabla\Delta {v}(t) \|_{L^2}^2  ) \,\mathrm{d}\tau\\
		&\ \ \ \ \ \ \ +\int_0^t (1+\tau)^\gamma \int_{\mathbbm{R}} |\nabla\Delta {v}(0,y,\tau)|^2 \,\mathrm{d}y  \,\mathrm{d}\tau\\
		&\le CM_0^2 +C\gamma \int_0^t (1+\tau)^{\gamma-1}\| \Delta v(\tau) \|_{H^1}^2 \,\mathrm{d}\tau+C\int_0^t (1+\tau)^\gamma(\| \nabla{v}(\tau) \|_{L^2}^2+\| \nabla^2{v}(\tau) \|_{L^2}^2 ) \,\mathrm{d}\tau \\
		&\ \ +C\int_0^t (1+\tau)^\gamma(\| \nabla{\rm{div}}{p}(\tau) \|_{L^2}^2+\| \nabla^2{\rm{div}}{p}(\tau) \|_{L^2}^2+(\varepsilon_0+\delta+M_0)\| \nabla{v}_{yy}(\tau) \|_{L^2}^2+\int_{\mathbbm{R}} (\Delta v(0,y,\tau))^2 \,\mathrm{d}y ) \,\mathrm{d}\tau.
	\end{aligned}
\end{equation}
In order to treat the term $\| \nabla^2{\rm{div}}{p}(t) \|_{L^2}^2$ on the right-hand side of $\eqref{eq-delv-H1-cjf}$, we
multiply $\eqref{eq-nab2divp-djf}$ by $(1+t)^\gamma$ and then integrate the resulting inequality over $[0,t]$. Combining $\eqref{eq-vyy-H1-cjf}$, we can obtain
\begin{equation}\label{eq-nab2divp-L2-djf}
\begin{aligned}[b]
	 &\int_0^t (1+\tau)^\gamma (\| \nabla^2{\rm{div}}{p}(\tau) \|_{L^2}^2+\| \nabla p(\tau) \|_{L^2}^2+\| \nabla{\rm{div}}{p}(\tau) \|_{L^2}^2) \,\mathrm{d}\tau\\
	\le& C \int_0^t (1+\tau)^\gamma( \| \nabla^2{v}(\tau) \|_{L^2}^2+\int_{\mathbbm{R}} (|\nabla {v}_{yy}(0,y,\tau)|^2+|\nabla {v}(0,y,\tau)|^2) \,\mathrm{d}y +\| \nabla{\rm{div}}{p}(\tau) \|_{L^2}^2+\| {\rm{div}}{p}(\tau) \|_{L^2}^2) \,\mathrm{d}\tau\\
	\le& CM_0^2+C\gamma \int_0^t (1+\tau)^{\gamma-1}\| v_{yy}(\tau) \|_{H^1}^2 \,\mathrm{d}\tau+C(\varepsilon_2+M_0^4)\int_0^t (1+\tau)^\gamma\| \nabla\Delta v(\tau) \|_{L^2}^2  \,\mathrm{d}\tau\\
	&  +  C \int_0^t (1+\tau)^\gamma( \| \nabla{v}(\tau) \|_{L^2}^2+\| \nabla^2{v}(\tau) \|_{L^2}^2+\int_{\mathbbm{R}} |\nabla {v}(0,y,\tau)|^2 \,\mathrm{d}y +\| \nabla{\rm{div}}{p}(\tau) \|_{L^2}^2+\| {\rm{div}}{p}(\tau) \|_{L^2}^2) \,\mathrm{d}\tau.
\end{aligned}
\end{equation}
Substituting $\eqref{eq-nab2divp-L2-djf}$ into $\eqref{eq-delv-H1-cjf}$, the inequality $\eqref{eq-delv-H1-cjf}$ can be rewritten  as
\begin{equation}\label{eq-delv-H1-cjf-1}
	\begin{aligned}[b]
		&(1+t)^\gamma\| \Delta v(t) \|_{H^1}^2+\int_0^t (1+\tau)^\gamma (\| \sqrt{\bar{u}_x}\Delta v(\tau) \|_{L^2}^2+\| \sqrt{\bar{u}_x}\nabla\Delta v(\tau) \|_{L^2}^2+\| \nabla\Delta {v}(\tau) \|_{L^2}^2  ) \,\mathrm{d}\tau\\
		&\ \ \ \ \ \ \ +\int_0^t (1+\tau)^\gamma \int_{\mathbbm{R}} |\nabla\Delta {v}(0,y,\tau)|^2 \,\mathrm{d}y  \,\mathrm{d}\tau\\
		&\le CM_0^2+C\gamma \int_0^t (1+\tau)^{\gamma-1}\| \nabla^2{v}(\tau) \|_{H^1}^2 \,\mathrm{d}\tau\\
		& + C \int_0^t (1+\tau)^\gamma( \| \nabla{v}(\tau) \|_{L^2}^2+\| \nabla^2{v}(\tau) \|_{L^2}^2+\int_{\mathbbm{R}} |\nabla {v}(0,y,\tau)|^2 \,\mathrm{d}y +\| \nabla{\rm{div}}{p}(\tau) \|_{L^2}^2+\| {\rm{div}}{p}(t) \|_{L^2}^2) \,\mathrm{d}\tau\\
		&+C(\varepsilon_2+\delta+M_0^2)\int_0^t (1+\tau)^\gamma \| \nabla{v}_{yy}(\tau) \|_{L^2}^2 \,\mathrm{d}\tau+C\int_0^t (1+\tau)^\gamma\int_{\mathbbm{R}} (\Delta v(0,y,\tau))^2 \,\mathrm{d}y  \,\mathrm{d}\tau.
	\end{aligned}
\end{equation}
 By employing Lemma $\ref{lem-nab2v-dengjia}$, for some small $\varepsilon_2$, $M_0$ and $\delta$, we can get from $\eqref{eq-vyy-H1-cjf}+\eqref{eq-delv-H1-cjf-1}+\eqref{eq-nab2divp-L2-djf}$ that
\begin{equation}\label{eq-nab3v-L2-djf}
	\begin{aligned}[b]
		&(1+t)^\gamma\| \nabla^2 v(t) \|_{H^1}^2+\int_0^t (1+\tau)^\gamma (\| \sqrt{\bar{u}_x}\nabla^2 v(\tau) \|_{L^2}^2+\| \sqrt{\bar{u}_x}\nabla^3 v(\tau) \|_{L^2}^2+\| \nabla^3 {v}(\tau) \|_{L^2}^2  ) \,\mathrm{d}\tau\\
		&\ \ +\int_0^t (1+\tau)^\gamma \int_{\mathbbm{R}} |\nabla^3 {v}(0,y,\tau)|^2 \,\mathrm{d}y  \,\mathrm{d}\tau+\int_0^t (1+\tau)^\gamma (\| \nabla^2{\rm{div}}{p}(\tau) \|_{L^2}^2+\| \nabla p(\tau) \|_{L^2}^2+\| \nabla{\rm{div}}{p}(\tau) \|_{L^2}^2) \,\mathrm{d}\tau\\
		&\le C\left(M_0^2 +\gamma \int_0^t (1+\tau)^{\gamma-1}(\| \nabla^2 v(\tau) \|_{L^2}^2+\| \nabla^3 v(\tau) \|_{L^2}^2) \,\mathrm{d}\tau+\int_0^t (1+\tau)^\gamma(\| \nabla{v}(\tau) \|_{L^2}^2+\| \nabla^2{v}(\tau) \|_{L^2}^2 ) \,\mathrm{d}\tau \right.\\
		&\left. \ \ +\int_0^t (1+\tau)^\gamma(\| \nabla{\rm{div}}{p}(\tau) \|_{L^2}^2+\| {\rm{div}}{p}(\tau) \|_{L^2}^2)\,\mathrm{d}\tau+\int_0^t (1+\tau)^\gamma\int_{\mathbbm{R}} ((\Delta v(0,y,\tau))^2+|\nabla {v}(0,y,\tau)|^2) \,\mathrm{d}y \,\mathrm{d}\tau\right).
	\end{aligned}
\end{equation}
When $\alpha$ is an integer, for any $\gamma\in[0,\alpha]$, the right-hand side of $\eqref{eq-nab3v-L2-djf}$ are bounded by
\begin{equation*}
	\begin{aligned}
	    &C\left(M_0^2 +\gamma \int_0^t (1+\tau)^{\gamma-1} \| \nabla^3 v(\tau) \|_{L^2}^2 \,\mathrm{d}\tau+\int_0^t (1+\tau)^\gamma(\| \nabla{v}(\tau) \|_{L^2}^2+\| \nabla^2{v}(\tau) \|_{L^2}^2 ) \,\mathrm{d}\tau \right.\\
		&\left. \ \ +\int_0^t (1+\tau)^\gamma(\| \nabla{\rm{div}}{p}(\tau) \|_{L^2}^2+\| {\rm{div}}{p}(\tau) \|_{L^2}^2)\,\mathrm{d}\tau+\int_0^t (1+\tau)^\gamma\int_{\mathbbm{R}} ((\Delta v(0,y,\tau))^2+|\nabla {v}(0,y,\tau)|^2) \,\mathrm{d}y \,\mathrm{d}\tau\right)\\
		&\le C M_\alpha^2,
	\end{aligned}
\end{equation*}
owing to $\eqref{eq-vsjl-gamma}$. Then, the inequality $\eqref{eq-nab3v-L2-djf}$ can be rewritten as
\begin{equation}\label{eq-nab3v-L2-1}
	\begin{aligned}[b]
		&(1+t)^\gamma\| \nabla^2 v(t) \|_{H^1}^2+\int_0^t (1+\tau)^\gamma (\| \sqrt{\bar{u}_x}\nabla^2 v(\tau) \|_{L^2}^2+\| \sqrt{\bar{u}_x}\nabla^3 v(\tau) \|_{L^2}^2+\| \nabla^3 {v}(\tau) \|_{L^2}^2  ) \,\mathrm{d}\tau\\
		\le& CM_\alpha^2 +C\gamma \int_0^t (1+\tau)^{\gamma-1}\| \nabla^3 v(\tau) \|_{L^2}^2 \,\mathrm{d}\tau,
	\end{aligned}
\end{equation}
for $\gamma\in[0,\alpha]$. Taking the value of $\gamma$ as $0,1,2,3,\cdots$, we can conclude that
\begin{equation*}
	\int_0^t \ (1+\tau)^{\gamma-1}\|\nabla^3 v (\tau)\|_{L^2}^2  \,d{\tau}\le C M_\alpha^2,\quad \forall \gamma\in[0,\alpha].
\end{equation*}
 Thus, we get from $\eqref{eq-nab3v-L2-1}$ that
\begin{equation}\label{eq-nab3v-L2-2}
		(1+t)^\gamma\| \nabla^2 v(t) \|_{H^1}^2+\int_0^t (1+\tau)^\gamma (\| \sqrt{\bar{u}_x}\nabla^2 v(\tau) \|_{L^2}^2+\| \sqrt{\bar{u}_x}\nabla^3 v(\tau) \|_{L^2}^2+\| \nabla^3 {v}(\tau) \|_{L^2}^2  ) \,\mathrm{d}\tau
		\le CM_\alpha^2,
\end{equation}
for $\gamma\in[0,\alpha]$. For the case of non-integer $\alpha$, taking $\gamma=\alpha+\varepsilon$, $\varepsilon>0$, the last six terms of $\eqref{eq-nab3v-L2-djf}$ is bounded by $C_\varepsilon(1+t)^\varepsilon M_\alpha^2$ owing to $\eqref{eq-vp-al-eps}$. The second term of $\eqref{eq-nab3v-L2-djf}$ can be estimated as
\begin{equation}
	\int_0^t (1+\tau)^{\alpha+\varepsilon-1}\|\nabla^3 v(\tau)\|_{L^2}^2 \,\mathrm{d}\tau \le (1+t)^\varepsilon\int_0^t (1+\tau)^{[\alpha]}\|\nabla^3 v(\tau)\|_{L^2}^2 \,\mathrm{d}\tau.
\end{equation}
Combining $\eqref{eq-nab3v-L2-2}$ with $\gamma=[\alpha]$, we get from $\eqref{eq-nab3v-L2-djf}$ with $\gamma=\alpha+\varepsilon$ that the inequality $\eqref{eq-dlv-al-eps}$ holds for non-integer $\alpha$, which completes the proof of Lemma $\ref{lem-dlv-al-eps}$.
$\hfill\Box$

Therefore, by utilizing weighted energy estimate in $H^3$, we have proved that if $\alpha$ is an integer
\begin{equation}\label{eq-vH3-integer}
(1+t)^{\gamma} \| v(t)\|_{H^{3}}^2 +\int_0^t (1+\tau)^{\gamma}\|\nabla v(\tau)\|_{H^{2}}^2 \,\mathrm{d}\tau \le C M_\alpha^2, \quad \forall \gamma\in[0,\alpha].
\end{equation}
If $\alpha$ is not an integer, it holds that
\begin{equation}\label{eq-vH3-noninteger}
		(1+t)^{\alpha+\varepsilon} \| v(t)\|_{H^{3}}^2 +\int_0^t (1+\tau)^{\alpha+\varepsilon}\|\nabla v(\tau)\|_{H^{2}}^2 \,\mathrm{d}\tau \le C(1+t)^\varepsilon  M_\alpha^2, \quad \forall \varepsilon>0.
\end{equation}
Using the time weighted energy method similar to the Lemma $\ref{lem-dlv-al-eps}$, we can obtain
\begin{equation}\label{eq-vt-alpha}
		\begin{aligned}[b]
			(1+t)^{\alpha+\varepsilon} \|{v}_t(t)\|_{H^2}^2 +\int_0^t (1+\tau)^{\alpha+\varepsilon} \|\nabla{v}_t(\tau) \|_{H^1}^2  \,\mathrm{d}\tau \le C(1+t)^\varepsilon M_\alpha^2, \quad \forall \varepsilon>0.
		\end{aligned}
\end{equation}
Finally, combining the relation between $p$ and $v$ that we used in Lemma $\ref{lem-p-gj}$, we can also get that for any $\varepsilon>0$,
\begin{equation*}
	(1+t)^{\alpha+\varepsilon}(\| {\rm{div}}{p}(t) \|_{H^3}^2+\| p(t) \|_{H^3}^2 )+\int_0^t (1+\tau)^{\alpha+\varepsilon}(\| \nabla^3{\rm{div}}{p}(\tau) \|_{L^2}^2+\| \nabla^3 p(\tau)\|_{L^2}^2+\| \nabla^2 p(\tau) \|_{L^2}^2)  \,\mathrm{d}\tau
	\le C(1+t)^\varepsilon M_\alpha^2.
\end{equation*}
Thus, the desired decay rate of $q$ can be obtained. We complete the proof of $\eqref{eq-q-sjsjgj}$.\\

In order to completes the proof of $\eqref{eq-convergence-rate}$, we need to show the convergence rate of $\partial_y^j v ~ (j=1,2,3)$ and $\partial_y^j \nabla v ~ (j=1,2)$. To explain the convergence rate of $\partial_y^kv(x,y,t)$ and $\partial_y^k \nabla{v}(x,y,t)$ simply, we introduce the time weighted energy norm $E(t)$ and the corresponding dissipation norm $D(t)$ by
\begin{equation}
	\begin{aligned}
		&E(t)^2:=\sum_{k=0}^2(1+t)^k\| \partial_y^k v(t) \|_{{H^{3-k}}}^2,\\
		&D(t)^2:=\sum_{k=0}^2(1+t)^k\| \partial_y^k \nabla{v}(t) \|_{{H^{2-k}}}^2.
	\end{aligned}
\end{equation}
\begin{lem}\label{lem-vy-sjl}
	Under the same assumptions as Lemma $\ref{lem-dlv-al-eps}$, put $N(t):=\sup_{0\le t\le T}E(t)$. There are positive constants $\varepsilon_3~(\le \varepsilon_2)$ and $C=C(\varepsilon_3)$, which are independent of $T$ such that if $N(t)+\delta\le\varepsilon_3$, then the solution $v(x,y,t)$ satisfies
\begin{equation}\label{eq-vy-gammal-jl}
\begin{aligned}
	\sum_{l=0}^2(1+t)^{\alpha+ l+\varepsilon} \| \partial_y^lv\|_{H^{3-l}}^2 +\sum_{l=0}^2\int_0^t (1+\tau)^{\alpha+ l+\varepsilon}\|\partial_y\nabla v\|_{H^{2-l}}^2 \,\mathrm{d}\tau \le C (1+t)^\varepsilon M_\alpha^2,
\end{aligned}		
\end{equation}
and
\begin{equation}\label{eq-vy-bj-sjl}
\begin{aligned}[b]
		&\int_0^t (1+\tau)^{\alpha+1+\varepsilon}\int_{\mathbbm{R}} (|\nabla{v}_y(0,y,\tau)|^2+|\Delta{v}_y(0,y,\tau)|^2) \,\mathrm{d}y  \,\mathrm{d}\tau +\int_0^t (1+\tau)^{\alpha+2+\varepsilon}\int_{\mathbbm{R}} |\nabla{v}_{yy}(0,y,\tau)|^2 \,\mathrm{d}y  \,\mathrm{d}\tau \\
		 \le& C (1+t)^\varepsilon M_\alpha^2,
\end{aligned}
\end{equation}
for $t\in[0,T]$ and any $\varepsilon>0$.
\end{lem}
{\it\bfseries Proof.}
The case of $l=0$ is proved by the combination of Lemmas $\ref{lem-alpha-gamma}$-$\ref{lem-dlv-al-eps}$. Next, we consider $l=1,2$ as follows. Firstly, we consider the case of $l=1$. We note that the case of $l=1$ in $\eqref{eq-vy-gammal-jl}$ contains the decay estimates for $v_y$, $v_{xy}$ and $v_{xxy}$. To do this, we take two steps to deal with it. As a first step, we give the decay rates for $v_y$ and $\nabla v_{y}$. The second step is that we give the convergence rates for $\nabla v_y$ and $\nabla^2 v_{y}$.

\textbf{Step 1.} Integrating $\eqref{eq-vnabv-y-1}$ over $\mathbbm{R}_+\times\mathbbm{R}$, we have
\begin{equation}\label{eq-vy-nabvy}
	\begin{aligned}[b]
		&\frac{\mathrm{d}}{\mathrm{d}t}(\|v_y(t)\|_{L^2}^2+\|\nabla{v}_y(t)\|_{L^2}^2)+(\|\sqrt{\bar{u}_x}v_y(t)\|_{L^2}^2+\|\sqrt{\bar{u}_x} \nabla{v}_y(t)\|_{L^2}^2)+\|\nabla{v}_y(t)\|_{L^2}^2+\int_{\mathbbm{R}} |\nabla{v}_y(0,y,t)|^2 \,\mathrm{d}y \\
		\le& C \left|\int_{\mathbbm{R}}\int_{\mathbbm{R}_+} \Big( \nabla{v}v_y^2+(\nabla^2 v+\bar{u}_{xx})v_y\nabla{v}_y+(\nabla{v}+\bar{u}_x)|\nabla{v}_y|^2+(\nabla{v}+\bar{u}_x)v_y^2\nabla{v}_y \Big) \,\mathrm{d}x \mathrm{d}y \right|.
	\end{aligned}
\end{equation}
The right-hand side of $\eqref{eq-vy-nabvy}$ can be estimated as follows. By the definition of $E(t)$ and $D(t)$,
\begin{equation*}
	\begin{cases}
        &\| \partial_y^j v(\tau) \|_{{H^{3-j}}} \le E(\tau)(1+\tau)^{-\frac{j}{2} }, \quad j=0,1,2,\\[1mm]
        &\| \partial_y^j \nabla{v}(\tau) \|_{{H^{2-j}}}\le D(\tau)(1+\tau)^{-\frac{j}{2} }, \quad j=0,1,2,
	\end{cases}
\end{equation*}
using $\eqref{eq-nabvy-inf}$, the following estimates of $v_y$ in $L^\infty$ norm can be obtained by simple computation:
	\begin{align}
		\| v_y \|_{L^\infty} \le CD(\tau)(1+\tau)^{-\frac{1}{2} },\label{eq-vy-Dt}\\
		\| v_y \|_{L^\infty} \le CE(\tau)(1+\tau)^{-\frac{3}{4} }.\label{eq-vy-Et}
	\end{align}
From $\eqref{eq-vy-Dt}$, the first term and the last term on the right-hand side of $\eqref{eq-vy-nabvy}$ can be estimated as
\begin{equation}\label{eq-nabvvy2-1}
	\int_{\mathbbm{R}}\int_{\mathbbm{R}_+} |\nabla{v}v_y^2| \,\mathrm{d}x \mathrm{d}y
	\le \| \nabla{v} \|_{L^2} \| v_y \|_{L^\infty} \| v_y \|_{L^2}
	\le CE(\tau)D(\tau)^2(1+\tau)^{-1},
\end{equation}
and
\begin{equation*}
	\int_{\mathbbm{R}}\int_{\mathbbm{R}_+}(|\nabla{v}|+\bar{u}_x)v_y^2|\nabla{v}_y|  \,\mathrm{d}x \mathrm{d}y
	\le (\|\nabla{v}\|_{L^\infty}+\|\bar{u}_x\|_{L^\infty})\| v_y \|_{L^\infty} \| v_y \|_{L^2} \|\nabla{v}_y\|_{L^2}\le C E(\tau)D(\tau)^2(1+\tau)^{ -\frac{3}{2}  }.
\end{equation*}
Using $\eqref{eq-vy-Et}$ and $\eqref{eq-DND-vy}$, the remaining terms on the right-hand side of $\eqref{eq-vy-nabvy}$ are bounded by
\begin{equation}\label{eq-nabvvy2-2}
	\int_{\mathbbm{R}}\int_{\mathbbm{R}_+}|\nabla^2 v||v_y||\nabla{v}_y|  \,\mathrm{d}x \mathrm{d}y\le \| \nabla^2 v \|_{L^2}\| v_y \|_{L^\infty} \| \nabla{v}_y \|_{L^2}  \le C E(\tau)D(\tau)^2(1+\tau)^{ -\frac{5}{4} },\\ 	
\end{equation}
and
\begin{equation}\label{eq-nabvvy2-3}
\begin{aligned}[b]
	&\int_{\mathbbm{R}}\int_{\mathbbm{R}_+}|\bar{u}_{xx}||v_y||\nabla{v}_y|  \,\mathrm{d}x \mathrm{d}y
	\le \| \nabla{v}_y \|_{L^2} \| \bar{u}_{xx}v_y \|_{L^2}\le C \delta \| \nabla{v}_y \|_{L^2}^2 ,\\
	&\int_{\mathbbm{R}}\int_{\mathbbm{R}_+}(|\nabla{v}|+\bar{u}_x)|\nabla{v}_y|^2  \,\mathrm{d}x \mathrm{d}y
	\le (\|\nabla{v}\|_{L^\infty}+\|\bar{u}_x\|_{L^\infty})\| \nabla{v}_y\|_{L^2}^2\le C(\varepsilon_3+\delta)\| \nabla{v}_y\|_{L^2}^2.
\end{aligned}
\end{equation}
Substituting $\eqref{eq-nabvvy2-1}$-$\eqref{eq-nabvvy2-3}$ into $\eqref{eq-vy-nabvy}$, we have
\begin{equation}\label{eq-vy-sjl}
	\begin{aligned}[b]
		&\frac{\mathrm{d}}{\mathrm{d}t}(\|v_y(t)\|_{L^2}^2+\|\nabla{v}_y(t)\|_{L^2}^2)+(\|\sqrt{\bar{u}_x}v_y(t)\|_{L^2}^2+\|\sqrt{\bar{u}_x} \nabla{v}_y(t)\|_{L^2}^2)+\|\nabla{v}_y(t)\|_{L^2}^2+\int_{\mathbbm{R}} |\nabla{v}_y(0,y,t)|^2 \,\mathrm{d}y \\
		\le& C E(\tau)D(\tau)^2(1+\tau)^{-1 }.
	\end{aligned}
\end{equation}
 Multiplying $\eqref{eq-vy-sjl}$ by $(1+t)^{\gamma+1}$ $(\gamma\ge 0)$ and integrating the resulting inequality over $\mathbbm{R}_+\times\mathbbm{R}\times[0,t]$, for some small $\delta$ and $\varepsilon_3$, we obtain
\begin{equation}\label{eq-vy-gj}
	\begin{aligned}[b]
        &(1+t)^{\gamma+1}(\| v_y(t) \|_{L^2}^2+\| \nabla{v}_y(t) \|_{L^2}^2)+\int_0^t (1+\tau)^{\gamma+1}(\|\sqrt{\bar{u}}_x v_y(\tau) \|_{L^2}^2+\|\sqrt{\bar{u}}_x \nabla{v}_y(\tau) \|_{L^2}^2) \,\mathrm{d}\tau \\
		&\ \ \ \ \ \ \ \ \ \ \ \ \ \ \ \ \ \ \ \ \ \ \ \ \ \ \ \ \ +\int_0^t (1+\tau)^{\gamma+1}\| \nabla{v}_y(\tau) \|_{L^2}^2 \,\mathrm{d}\tau+\int_0^t (1+\tau)^{\gamma+1}\int_{\mathbbm{R}} |\nabla{v}_y(0,y,\tau)|^2 \,\mathrm{d}y  \mathrm{d}\tau \\
		\le& C\left(M_0^2+(\gamma+1)\int_0^t (1+\tau)^{\gamma}(\| v_y(\tau) \|_{L^2}^2+\| \nabla{v}_y(\tau) \|_{L^2}^2) \,\mathrm{d}\tau+N(t)\int_0^t (1+\tau)^\gamma D(\tau)^2 \,\mathrm{d}\tau \right).
	\end{aligned}
\end{equation}
To complete the case of $l=1$ in $\eqref{eq-vy-gammal-jl}$, we still need to give the decay estimate for $v_{xxy}$.

\textbf{Step 2.} By employing the similar method, we can get from $\nabla v_y\cdot\nabla\partial_y \eqref{eq-vp-rd}_1+{\rm{div}}{p}_y\times {\rm{div}} \partial_y \eqref{eq-vp-rd}_2+\Delta{v}_y\times \Delta \partial_y \eqref{eq-vp-rd}_1- \Delta{\rm{div}}{p}_y\times {\rm{div}} \partial_y \eqref{eq-vp-rd}_2$ that
\begin{equation}\label{eq-v-xxy-sj}
	\begin{aligned}[b]
		&\left\{ \frac{1}{2} |\nabla{v}_y|^2+\frac{1}{2} (\Delta v_y)^2\right\}_t+\frac{1}{2} \bar{u}_x(|\nabla{v}_y|^2+ (\Delta v_y)^2)+\bar{u}_xv_{xy}^2+ (\Delta v_y)^2 \\
		&\ \ \ \ \ \ \ \ \ \ \ \ \ \ \ \ \ +{\rm{div}}\{{\rm{div}}{p}_y \nabla{v}_y\}+\frac{1}{2}\left\{( v+\bar{u})|\nabla{v}_y|^2+ (v +\bar{u}) (\Delta v_y)^2\right\}_x\\
		&=-\frac{1}{2} v_x|\nabla{v}_y|^2-v_y \nabla{v}_x\cdot \nabla{v}_y-v_{xy}\nabla{v}\cdot \nabla{v}_y- \bar{u}_{xx}v_yv_{xy}-\nabla\left(g'(v+\bar{u})v_y\right)_y\cdot \nabla{v}_y\\
		&\ \ -\frac{1}{2} v_x(\Delta v_y)^2-v_y \Delta v_x \Delta v_y-2 \nabla{v}_y\cdot\nabla{v}_x\Delta v_y-\Delta vv_{xy}\Delta v_y-2\nabla{v}\cdot\nabla{v}_{xy}\Delta v_y- \bar{u}_{xxx}v_y\Delta v_y\\
		&\ \ -3 \bar{u}_{xx}v_{xy}\Delta v_y-2 \bar{u}_xv_{xxy}\Delta v_y-\Delta\left(g'(v+\bar{u})v_y\right)_y\Delta v_y:=\sum_{i=1}^{14} \RNum{1}_i.
	\end{aligned}
\end{equation}
 By simple computation,
\begin{equation*}
		 \nabla\left(g'(v+\bar{u})v_y\right)_y\cdot \nabla{v}_y
		\sim (\nabla{v}\cdot \nabla{v}_y+\bar{u}_xv_{xy})v_y^2+v_y|\nabla{v}_y|^2+(\nabla{v}\cdot \nabla{v}_y+\bar{u}_xv_{xy})v_{yy}+\left\{g'(v+\bar{u})|\nabla v_y|^2\right\}_y,
\end{equation*}
and using $\eqref{eq-del-g}$
\begin{equation*}
    \begin{aligned}
		\Delta\left(g'(v+\bar{u})v_y\right)_y\Delta v_y
		\sim& (v_x+\bar{u}_x)^2v_y^2\Delta v_y+\bar{u}_{xx}v_y^2\Delta v_y+v_y^3\Delta v_y+v_y^2\Delta v \Delta v_y+v_y(\Delta v_y)^2+|\nabla v_y|^2\Delta v_y\\
		&+(\nabla{v}\cdot \nabla{v}_y+\bar{u}_xv_{xy})\Delta v_y+(v_x+\bar{u}_x)^2v_{yy}\Delta v_y+\bar{u}_{xx}v_{yy}\Delta v_y+v_y^2v_{yy}\Delta v_y+v_{yy}\Delta v\Delta v_y\\
		&+(\nabla{v}\cdot \nabla{v}_{yy}+\bar{u}_xv_{xyy})\Delta v_y+\left\{g'(v+\bar{u})(\Delta v_y)^2\right\}_y.
	\end{aligned}
\end{equation*}
Notice that the equation $\eqref{eq-vp-rd}_1$ gives
\begin{equation}\label{eq-divpy-0}
	{\rm{div}}{p}_y(0,y,t)=-u_-v_{xy}(0,y,t).
\end{equation}
 After integration over $\mathbbm{R}_+\times\mathbbm{R}$, the boundary term is bounded by
$$\int_{\mathbbm{R}} {\rm{div}}{p}_y(0,y,t)v_{xy}(0,y,t) \, \mathrm{d}y\le C \int_{\mathbbm{R}} (v_{xy}(0,y,t))^2 \,\mathrm{d}x.$$
Then, the terms on the right-hand side of $\eqref{eq-v-xxy-sj}$ can be estimated as
\begin{equation*}
\begin{aligned}
		&\RNum{1}_1+\RNum{1}_3\le \| \nabla{v} \|_{L^\infty} \| \nabla{v}_y \|_{L^2}^2,\\
		&\RNum{1}_2\le \| v_y \|_{L^\infty} \| \nabla{v}_x \|_{L^2}\| \nabla{v}_y \|_{L^2}  \le C E(\tau)D(\tau)^2(1+\tau)^{-\frac{5}{4}},\\
		&\RNum{1}_4\le \| \bar{u}_{xx}v_y \|_{L^2} \| \nabla{v}_y \|_{L^2} \le C\| \nabla{v}_y \|_{L^2}^2,\\
		&\RNum{1}_5\le C\| v_y \|_{L^\infty} \| v_y \|_{L^2} \| \nabla{v}_y \|_{L^2}+C (\| \nabla{v} \|_{L^\infty}+\| \bar{u}_x \|_{L^\infty}  )  \| \nabla{v}_y \|_{L^2}^2\le CE(\tau)D(\tau)^2(1+\tau)^{-\frac{3}{2}}+C \| \nabla{v}_y \|_{L^2}^2,\\
		&\RNum{1}_6\le \| \nabla{v} \|_{L^\infty} \| \Delta v_y \|_{L^2}^2\le C \varepsilon_3 \| \Delta v_y \|_{L^2}^2,\\
		&\RNum{1}_7\le \| v_y \|_{L^\infty} \| \Delta v_x \|_{L^2}\| \Delta v_y \|_{L^2} \le C  E(\tau)D(\tau)^2(1+\tau)^{-\frac{5}{4}},\\
		&\RNum{1}_8\le \| \nabla{v}_y \|_{L_x^2(L_y^\infty)} \| \nabla{v}_x \|_{L_x^\infty(L_y^2)} \| \Delta v_y \|_{L^2}\le C E(\tau)D(\tau)^2(1+\tau)^{-\frac{5}{4}},\\
		&\RNum{1}_9\le \| v_{xy} \|_{L_x^2(L_y^\infty)} \| \Delta v \|_{L_x^\infty(L_y^2)} \| \Delta v_y \|_{L^2}\le C E(\tau)D(\tau)^2(1+\tau)^{-\frac{5}{4}},\\
		&\RNum{1}_{10}\le \| \nabla{v} \|_{L^\infty} \| \nabla^2{v}_y \|_{L^2}^2,\\
		&\RNum{1}_{11}\le \| \bar{u}_{xxx}v_y \|_{L^2}\| \Delta v_y \|_{L^2} \le  \delta\| \Delta v_y \|_{L^2}^2+C \| \nabla{v}_y \|_{L^2}^2,\\
		&\RNum{1}_{12}\le \| \bar{u}_{xx} \|_{L^\infty} \| \nabla{v}_y \|_{L^2} \| \Delta v_y \|_{L^2}\le \delta \| \Delta v_y \|_{L^2}^2+ C \| \nabla{v}_y \|_{L^2}^2,\\
		&\RNum{1}_{13}\le \| \bar{u}_x \|_{L^\infty} \| \nabla^2{v}_y \|_{L^2}^2\le C \delta\| \nabla^2{v}_y \|_{L^2}^2,\\
\end{aligned}
\end{equation*}
and
\begin{equation*}
	\begin{aligned}
		\RNum{1}_{14}\le& C \| v_y \|_{L^\infty} \| v_y \|_{L^2} \| \Delta {v}_y \|_{L^2}+C \| v_y \|_{L^\infty}\| \Delta v \|_{L^2}\| \Delta v_y \|_{L^2}  +C\| v_y \|_{L^\infty} \| \Delta v_y \|_{L^2}^2\\
		&+C\| \nabla{v}_y \|_{L_x^\infty(L_y^2)}\| \nabla{v}_y \|_{L_x^2(L_y^\infty)} \| \Delta v_y \|_{L^2} +C(\| \nabla{v} \|_{L^\infty}+\| \bar{u}_x \|_{L^\infty}+\| \bar{u}_{xx} \|_{L^\infty} )\| \nabla{v}_y \|_{L^2}\| \Delta{v}_y \|_{L^2}\\
		&+C \| \Delta v \|_{L_x^\infty(L_y^2)}\| v_{yy} \|_{L_x^2(L_y^\infty)}\| \Delta v_y \|_{L^2}  +C (\| \nabla{v} \|_{L^\infty}+\| \bar{u}_x \|_{L^\infty})\| \nabla{v}_{yy} \|_{L^2} \| \Delta v_y \|_{L^2} \\
		\le& C E(\tau)D(\tau)^2(1+\tau)^{-\frac{5}{4}} +C(\varepsilon_3+\delta)\| \nabla^2{v}_{yy} \|_{L^2}^2+C \| \nabla{v}_y \|_{L^2}^2.
	\end{aligned}
\end{equation*}
Therefore, we choose smallness $\delta$ and $\varepsilon_3$ such that $C(\delta+\varepsilon_3)\le \frac{1}{4}
$, then the right-hand side of $\eqref{eq-v-xxy-sj}$ is bounded by
\begin{equation*}
	\frac{1}{4}\| \nabla^2 v_y(t) \|_{L^2}^2+ C\int_{\mathbbm{R}}\int_{\mathbbm{R}_+} |\nabla v_{y}|^2 \,\mathrm{d}x \mathrm{d}y+CE(t)D(t)^2(1+t)^{-\frac{5}{4}},
\end{equation*}
which follows $\eqref{eq-v-xxy-sj}$ that
\begin{equation}\label{eq-delvy-1}
	\begin{aligned}[b]
		&\frac{\mathrm{d}}{\mathrm{d}t}(\|\nabla{v}_y\|_{L^2}^2+\|\nabla^2 v_y\|_{L^2}^2 )+(\|\sqrt{\bar{u}_x}\nabla{v}_y\|_{L^2}^2+ \|\sqrt{\bar{u}_x}\Delta v_y\|_{L^2}^2)+\|\nabla^2 v_y(t)\|_{L^2}^2+\int_{\mathbbm{R}} |\Delta v_y(0,y,t)|^2 \, \mathrm{d}y \\
		\le& C\left(\int_{\mathbbm{R}} |\nabla{v}_y(0,y,t)|^2 \,\mathrm{d}y +\int_{\mathbbm{R}}\int_{\mathbbm{R}_+} |\nabla v_{y}|^2 \,\mathrm{d}x \mathrm{d}y+E(t)D(t)^2(1+t)^{-\frac{5}{4}}  \right).
	\end{aligned}
\end{equation}
Here we have used the fact that $v_{yy}(0,y,t)=0$ and
\begin{equation*}
	\begin{aligned}
		&(\Delta v_y)^2=|\nabla^2{v}_{y}|^2+2\{v_{xxy}v_{yy}\}_y-2\{v_{xyy}v_{yy}\}_x,\\
		&\int_{\mathbbm{R}}\int_{\mathbbm{R}_+} (\Delta v_y) \,\mathrm{d}x \mathrm{d}y =\int_{\mathbbm{R}}\int_{\mathbbm{R}_+} |\nabla^2{v}_{y}|^2 \,\mathrm{d}x \mathrm{d}y.
	\end{aligned}
\end{equation*}
 Multiply $\eqref{eq-delvy-1}$ by $(1+t)^{\gamma+1}$ $(\gamma\ge 0)$ and integrate the resulting estimate over $[0,t]$. By employing $\eqref{eq-vy-gj}$, we have
\begin{equation}\label{eq-delvy-2}
	\begin{aligned}[b]
		&(1+t)^{\gamma+1}(\|\nabla{v}_y(t)\|_{L^2}^2+\|\nabla^2 v_y(t) \|_{L^2}^2) +\int_0^t (1+\tau)^{\gamma+1}(\|\sqrt{\bar{u}_x}\nabla{v}_y(\tau)\|_{L^2}^2+ \|\sqrt{\bar{u}_x}\Delta v_y(\tau)\|_{L^2}^2) \,\mathrm{d}\tau\\
		& +\int_0^t (1+\tau)^{\gamma+1}\|\nabla^2 v_y(\tau)\|_{L^2}^2 \,\mathrm{d}\tau+\int_0^t (1+\tau)^{\gamma+1} \int_{\mathbbm{R}} |\Delta v_y(0,y,\tau)|^2 \, \mathrm{d}y \mathrm{d}\tau \\
		\le& C\left(M_0^2+(\gamma+1)\int_0^t (1+\tau)^{\gamma}(\|{v}_y(\tau)\|_{L^2}^2+\|\nabla{v}_y(\tau)\|_{L^2}^2+\|\nabla^2 v_y(\tau) \|_{L^2}^2) \,\mathrm{d}\tau +N(t)\int_0^t (1+\tau)^{\gamma}D(\tau)^2 \,\mathrm{d}\tau \right).
	\end{aligned}
\end{equation}
Adding $\eqref{eq-vy-gj}$ and $\eqref{eq-delvy-2}$, for any $\gamma\ge 0$, we have
\begin{equation}\label{eq-vy-H2-sjgj}
	\begin{aligned}[b]
		&(1+t)^{\gamma+1}\|{v}_y(t)\|_{H^2}^2 +\int_0^t (1+\tau)^{\gamma+1}(\|\sqrt{\bar{u}_x}{v}_y(\tau)\|_{L^2}^2+\|\sqrt{\bar{u}_x}\nabla{v}_y(\tau)\|_{L^2}^2+ \|\sqrt{\bar{u}_x}\Delta v_y(\tau)\|_{L^2}^2) \,\mathrm{d}\tau\\
		& +\int_0^t (1+\tau)^{\gamma+1}\|\nabla v_y(\tau)\|_{H^1}^2 \,\mathrm{d}\tau+\int_0^t (1+\tau)^{\gamma+1} \int_{\mathbbm{R}} (|\nabla v_y(0,y,\tau)|^2+(\Delta v_y(0,y,\tau))^2) \, \mathrm{d}y \mathrm{d}\tau \\
		\le& C\left(M_0^2+(\gamma+1)\int_0^t (1+\tau)^{\gamma}\|{v}_y(\tau)\|_{H^2}^2 \,\mathrm{d}\tau +N(t)\int_0^t (1+\tau)^{\gamma}D(\tau)^2 \,\mathrm{d}\tau \right).
	\end{aligned}
\end{equation}
Taking $\gamma=\alpha+\varepsilon$ in $\eqref{eq-vy-H2-sjgj}$, for any $\varepsilon>0$, combining $\eqref{eq-vH3-noninteger}$, we can obtain
\begin{equation}\label{eq-vy-H2-aleps}
	\begin{aligned}[b]
		&(1+t)^{\alpha+1+\varepsilon}\|{v}_y(t)\|_{H^2}^2 +\int_0^t (1+\tau)^{\alpha+1+\varepsilon}(\|\sqrt{\bar{u}_x}{v}_y(\tau)\|_{L^2}^2+\|\sqrt{\bar{u}_x}\nabla{v}_y(\tau)\|_{L^2}^2+ \|\sqrt{\bar{u}_x}\Delta v_y(\tau)\|_{L^2}^2) \,\mathrm{d}\tau\\
		& +\int_0^t (1+\tau)^{\alpha+1+\varepsilon}\|\nabla v_y(\tau)\|_{H^1}^2 \,\mathrm{d}\tau+\int_0^t (1+\tau)^{\alpha+1+\varepsilon} \int_{\mathbbm{R}} (|\nabla v_y(0,y,\tau)|^2+(\Delta v_y(0,y,\tau))^2) \, \mathrm{d}y \mathrm{d}\tau \\
		\le& C\left(M_0^2+(\alpha+1+\varepsilon)\int_0^t (1+\tau)^{\alpha+\varepsilon}\|{v}_y(\tau)\|_{H^2}^2 \,\mathrm{d}\tau +N(t)\int_0^t (1+\tau)^{\gamma}D(\tau)^2 \,\mathrm{d}\tau \right)\\
		\le& C(1+t)^\varepsilon M_\alpha^2+C N(t)\int_0^t (1+\tau)^{\gamma}D(\tau)^2 \,\mathrm{d}\tau.
	\end{aligned}
\end{equation}

\textbf{Step 3.} Now we consider the case of $l=2$. Recalling $\eqref{eq-vyyH1}$, we have
\begin{equation}\label{eq-vnabv-yy}
	\begin{aligned}[b]
		&\frac{\mathrm{d}}{\mathrm{d}t}(\| v_{yy} (t)\|_{L^2}^2\!+\!\| \nabla{v}_{yy}(t) \|_{L^2}^2 )\!+\!\|\sqrt{\bar{u}_x} v_{yy}(\tau)  \|_{L^2}^2 \!+\!\|\sqrt{\bar{u}_x}  \nabla{v}_{yy}(\tau)  \|_{L^2}^2\!+\!\| \nabla{v}_{yy}(\tau)  \|_{L^2}^2\!+\!\int_{\mathbbm{R}} |\nabla{v}_{yy}(0,y,\tau)|^2 \,\mathrm{d}y \\
		&\le C \int_{\mathbbm{R}}\int_{\mathbbm{R}_+} (|\nabla{v}||v_{yy}|^2\!+\!|v_y||\nabla{v}_y||v_{yy}|\!+\!|v_y|^3|v_{yy}|\!+\!|\nabla{v}||\nabla{v}_{yy}|^2\!+\!|\nabla{v}_y|^2|\nabla{v}_{yy}|\!+\!|v_{yy}||\nabla^2 v||\nabla{v}_{yy}|) \,\mathrm{d}x \mathrm{d}y \\
		&+C \int_{\mathbbm{R}}\int_{\mathbbm{R}_+} (|v_y||\nabla^2 v_y||\nabla{v}_{yy}|+|\bar{u}_{xx}||v_{yy}||\nabla{v}_{yy}|+(|\nabla{v}|+\bar{u}_x)|v_y|^3|\nabla{v}_{yy}|+|v_y|^2|\nabla{v}_y||\nabla{v}_{yy}|\,\mathrm{d}x \mathrm{d}y\\
		&+C \int_{\mathbbm{R}}\int_{\mathbbm{R}_+} (|\nabla{v}|+\bar{u}_x)|v_y||v_{yy}||\nabla{v}_{yy}|+|\nabla{v}_y||v_{yy}||\nabla{v}_{yy}|+|v_y||\nabla{v}_{yy}|^2+(|\nabla{v}|+\bar{u}_x)|v_{yyy}||\nabla{v}_{yy}|) \,\mathrm{d}x \mathrm{d}y\\
		&:=\sum_{i=1}^{14} \RNum{2}_i.
	\end{aligned}
\end{equation}
Using $\eqref{eq-DND-vy}$, $\eqref{eq-Lxinf}$, $\eqref{eq-Lyinf}$, $\eqref{eq-vy-Dt}$, $\eqref{eq-vy-Et}$ and the $\mathrm{H\ddot{o}lder}$ inequality and the right-hand side can be estimated as
\begin{equation*}
	\begin{aligned}
		&\RNum{2}_1\le C\| \nabla{v} \|_{L_x^2(L_y^\infty)}\| v_{yy} \|_{L_x^\infty(L_y^2)}\| v_{yy} \|_{L^2}\le C E(\tau)D(\tau)^2(1+\tau)^{-2 },\\
		&\RNum{2}_2\le C \| v_y \|_{L^\infty}\| \nabla{v}_y \|_{L^2}\| v_{yy} \|_{L^2}\le  C E(\tau)D(\tau)^2(1+\tau)^{-2 },\\
		&\RNum{2}_3\le C \| v_y \|_{L^\infty}^2 \| v_y \|_{L^2}\| v_{yy} \|_{L^2} \le C \| v_y \|_{L^\infty}^2 \| v_{yy} \|_{L^2}\le  C E(\tau)D(\tau)^2(1+\tau)^{-2 },\\
		&\RNum{2}_4+\RNum{2}_8\le C\| \nabla{v} \|_{L^\infty} \| \nabla{v}_{yy} \|_{L^2}^2+C\| \bar{u}_{xx}v_{yy} \|_{L^2}\| \nabla{v}_{yy} \|_{L^2}\le C(\delta+\varepsilon_3)\| \nabla{v}_{yy} \|_{L^2}^2,\\
		&\RNum{2}_5+ \RNum{2}_{12}\le C \| \nabla{v}_y \|_{L_x^2(L_y^\infty)} \| \nabla{v}_y \|_{L_x^\infty(L_y^2)}\| \nabla{v}_{yy} \|_{L^2}\le C E(\tau)D(\tau)^2(1+\tau)^{-\frac{9}{4}},\\
		&\RNum{2}_6\le C \| v_{yy} \|_{L_x^2(L_y^\infty)} \| \nabla^2{v} \|_{L_x^\infty(L_y^2)} \| \nabla{v}_{yy} \|_{L^2}\le C E(\tau)D(\tau)^2(1+\tau)^{-2},\\
		&\RNum{2}_7\le C \| v_y \|_{L^\infty}\| \nabla^2{v}_{y} \|_{L^2}\| \nabla{v}_{yy} \|_{L^2} \le C E(\tau)D(\tau)^2(1+\tau)^{-\frac{9}{4}},\\
		&\RNum{2}_9\le C \| {v}_y \|_{L^\infty} \| v_y \|_{L^2}\| \nabla{v}_{yy} \|_{L^2}\le C E(\tau)D(\tau)^2(1+\tau)^{-2},\\
		& \RNum{2}_{10}+ \RNum{2}_{11}\le C \| v_y \|_{L^\infty} \| \nabla{v}_y \|_{L^2} \| \nabla{v}_{yy} \|_{L^2}\le C E(\tau)D(\tau)^2(1+\tau)^{-\frac{9}{4}},\\
        & \RNum{2}_{13}+ \RNum{2}_{14}\le C (\| \nabla{v} \|_{L^\infty}+\| \bar{u}_x \|_{L^\infty} )\| \nabla{v}_{yy} \|_{L^2}\le C(\delta+\varepsilon_3)\| \nabla{v}_{yy} \|_{L^2}^2.
	\end{aligned}
\end{equation*}
Thus, the terms on right-hand side of $\eqref{eq-vnabv-yy}$ are bounded by
\begin{equation}\label{eq-nabvyy-rightbound}
	\begin{aligned}
		 C(\delta+\varepsilon_3)\| \nabla{v}_{yy} \|_{L^2}^2+C E(\tau)D(\tau)^2(1+\tau)^{-2}.
	\end{aligned}
\end{equation}
Notice that the first term of $\eqref{eq-nabvyy-rightbound}$ can absorbed in the fifth term on the left-hand side of $\eqref{eq-vnabv-yy}$ if $\varepsilon_3$ and $\delta$ are small enough. Consequently, we can get from $\eqref{eq-vnabv-yy}$ that
\begin{equation}\label{eq-nabvyy-jfgj}
	\begin{aligned}[b]
		&\frac{\mathrm{d}}{\mathrm{d}t}(\| v_{yy} (t)\|_{L^2}^2\!+\!\| \nabla{v}_{yy}(t) \|_{L^2}^2 )\!+\!\|\sqrt{\bar{u}_x} v_{yy}(\tau)  \|_{L^2}^2 \!+\!\|\sqrt{\bar{u}_x}  \nabla{v}_{yy}(\tau)  \|_{L^2}^2\!+\!\| \nabla{v}_{yy}(\tau)  \|_{L^2}^2\!+\!\int_{\mathbbm{R}} |\nabla{v}_{yy}(0,y,\tau)|^2 \,\mathrm{d}y \\
		&\le C E(\tau)D(\tau)^2(1+\tau)^{-2}.
	\end{aligned}
\end{equation}
Multiplying $\eqref{eq-nabvyy-jfgj}$ by $(1+t)^{\gamma+2}$ and integrating the resulting inequality over $\mathbbm{R}_+\times\mathbbm{R}\times[0,t]$, for any $\gamma\ge0$, we obtain
\begin{equation}\label{eq-v-yy-gj}
	\begin{aligned}[b]
		&(1+t)^{\gamma+2}(\| v_{yy}(t) \|_{L^2}^2+\| \nabla{v}_{yy}(t) \|_{L^2}^2)+\int_0^t (1+\tau)^{\gamma+2}(\|\sqrt{\bar{u}}_x v_{yy}(\tau) \|_{L^2}^2+\|\sqrt{\bar{u}}_x \nabla{v}_{yy}(\tau) \|_{L^2}^2) \,\mathrm{d}\tau\\
		&\ \ \ \ \ \ \ \ \ \ \ \ \ \ \ \ \ \ \ \ \ \ \ \ \  +\int_0^t (1+\tau)^{\gamma+2}\| \nabla{v}_{yy}(\tau) \|_{L^2}^2 \,\mathrm{d}\tau+\int_0^t (1+\tau)^{\gamma+2}\int_{\mathbbm{R}} |\nabla{v}_{yy}(0,y,\tau)|^2 \,\mathrm{d}y  \mathrm{d}\tau\\
		\le& C\left(M_0^2+(\gamma+2)\int_0^t (1+\tau)^{\gamma+1}(\| v_{yy}(\tau) \|_{L^2}^2+\| \nabla{v}_{yy}(\tau) \|_{L^2}^2) \,\mathrm{d}\tau+N(t)\int_0^t (1+\tau)^{\gamma} D(\tau)^2 \,\mathrm{d}\tau \right).
	\end{aligned}
\end{equation}
By utilizing $\eqref{eq-vy-H2-aleps}$, taking $\gamma=\alpha+\varepsilon$ in $\eqref{eq-v-yy-gj}$, we can deduce that
\begin{equation}\label{eq-vyy-H1-sjgj}
	\begin{aligned}[b]
		&(1+t)^{\alpha+2+\varepsilon}\| v_{yy}(t) \|_{H^1}^2+\int_0^t (1+\tau)^{\alpha+2+\varepsilon}(\|\sqrt{\bar{u}}_x v_{yy}(\tau) \|_{L^2}^2+\|\sqrt{\bar{u}}_x \nabla{v}_{yy}(\tau) \|_{L^2}^2) \,\mathrm{d}\tau\\
		&\ \ \ \ \ \ \ \ \ \ \ \ \ \ \ \ \ \ \ \ \ \ \ \ \  +\int_0^t (1+\tau)^{\alpha+2+\varepsilon}\| \nabla{v}_{yy}(\tau) \|_{L^2}^2 \,\mathrm{d}\tau+\int_0^t (1+\tau)^{\alpha+2+\varepsilon}\int_{\mathbbm{R}} |\nabla{v}_{yy}(0,y,t)|^2 \,\mathrm{d}y  \mathrm{d}\tau\\
		\le& C\left(M_0^2+(\alpha+2+\varepsilon)\int_0^t (1+\tau)^{\alpha+1+\varepsilon}\| v_{yy}(\tau) \|_{H^1}^2 \,\mathrm{d}\tau+N(t)\int_0^t (1+\tau)^{\alpha+\varepsilon} D(\tau)^2 \,\mathrm{d}\tau \right)\\
		\le& C(1+t)^{\varepsilon}M_\alpha^2+CN(t)\int_0^t (1+\tau)^{\alpha+\varepsilon} D(\tau)^2 \,\mathrm{d}\tau.
	\end{aligned}
\end{equation}
Finally, we need to treat the last term in $\eqref{eq-vy-H2-aleps}$ and $\eqref{eq-vyy-H1-sjgj}$. We add $\eqref{eq-vH3-noninteger}$, $\eqref{eq-vy-H2-aleps}$ and $\eqref{eq-vyy-H1-sjgj}$ to obtain
\begin{equation}\label{eq-vy-gamma-l}
\begin{aligned}[b]
	&\sum_{l=0}^2(1+t)^{\alpha+l+\varepsilon} \| \partial_y^lv\|_{H^{3-l}}^2 +\sum_{l=0}^2\int_0^t (1+\tau)^{\alpha+l+\varepsilon}\|\partial_y^l\nabla v\|_{H^{2-l}}^2 \,\mathrm{d}\tau \\
	\le& C(1+t)^\varepsilon M_\alpha^2+CN(t)\int_0^t (1+\tau)^{\alpha+\varepsilon} D(\tau)^2 \,\mathrm{d}\tau.
\end{aligned}		
\end{equation}
  The last term on the right-hand side of $\eqref{eq-vy-gamma-l}$ can be absorbed in the second term on the left-hand side if $N(t)\le \varepsilon_3$ is sufficiently small. Therefore, the desired inequality $\eqref{eq-vy-gammal-jl}$ can be obtained.
$\hfill\Box$

 \begin{lem}\label{lem-py}
 	Under the same assumptions as Lemma $\ref{lem-vy-sjl}$, it holds that
 	\begin{equation}\label{eq-py-sjl-jl}
 	   \begin{aligned}[b]
 	   		&(1+t)^{\alpha+1+\varepsilon}(\| {\rm{div}}{p}_y(t) \|_{H^2}^2+\| p_{1y}(t)\|_{H^2}^2+\| p_{2}(t)\|_{H^3}^2) \le C(1+t)^\varepsilon M_\alpha^2,
 	   	\end{aligned}	
 	\end{equation}
 	for $t\in[0,T]$ and any $\varepsilon>0$.
 \end{lem}
 {\it\bfseries Proof.}
Rewrite ${\rm{div}}\partial_y\eqref{eq-vp-rd}_{2}$ in the form $\Delta {\rm{div}}{p}_y- {\rm{div}}{p}_y=\Delta v_y$, and then square this equation. Consequently, integrating the resulting equation over $\mathbbm{R}_+\times\mathbbm{R}$, using $\eqref{eq-divpy-0}$, we have
\begin{equation}\label{eq-py-sjv-1}
	\begin{aligned}[b]
		&\int_{\mathbbm{R}}\int_{\mathbbm{R}_+} ((\Delta {\rm{div}}{p}_y)^2+2|\nabla{\rm{div}}{p}_y|^2+({\rm{div}}{p}_y)^2) \,\mathrm{d}x \mathrm{d}y \\
		=&\int_{\mathbbm{R}}\int_{\mathbbm{R}_+} (\Delta v_y)^2 \,\mathrm{d}x \mathrm{d}y-2 \int_{\mathbbm{R}}  {\rm{div}}{p}_{xy}(0,y,t){\rm{div}}{p}_y(0,y,t) \, \mathrm{d}y \\
		=&\int_{\mathbbm{R}}\int_{\mathbbm{R}_+} (\Delta v_y)^2 \,\mathrm{d}x \mathrm{d}y+2u_- \int_{\mathbbm{R}}  {\rm{div}}{p}_{xy}(0,y,t)v_{xy}(0,y,t) \, \mathrm{d}y \\
		\le&\int_{\mathbbm{R}}\int_{\mathbbm{R}_+} (\Delta v_y)^2 \,\mathrm{d}x \mathrm{d}y+\delta \int_{\mathbbm{R}} {\rm{div}}{p}_{xy}(0,y,t) \,\mathrm{d}y +C \delta^{-1}\int_{\mathbbm{R}} (v_{xy}(0,y,t))^2 \,\mathrm{d}y\\
		\le& \int_{\mathbbm{R}}\int_{\mathbbm{R}_+} (\Delta v_y)^2 \,\mathrm{d}x \mathrm{d}y+\delta \int_{\mathbbm{R}}\int_{\mathbbm{R}_+} \left(({\rm{div}}{p}_{xy})^2+({\rm{div}}{p}_{xxy})^2\right) \,\mathrm{d}x \mathrm{d}y  +C \int_{\mathbbm{R}} (v_{xy}^2+v_{xxy}^2) \,\mathrm{d}y.\\
	\end{aligned}
\end{equation}
Notice that the term ${\rm{div}}{p}_{xy}$ on the right-hand side of $\eqref{eq-py-sjv-1}$ can be absorbed in the second term on the left if $\delta$ is sufficiently small. Combining the results of Lemma $\ref{lem-vy-sjl}$, we just need to get the estimate of $\delta\| {\rm{div}}{p}_{xxy} \|_{L^2}^2$. We know that
$$\delta\| {\rm{div}}{p}_{xxy} \|_{L^2}^2\le C\delta \| \Delta {\rm{div}}{p}_y \|_{L^2}^2+C\delta\| {\rm{div}}{p}_{yyy} \|_{L^2}^2,$$
and the former term can be absorbed in the first term of $\eqref{eq-py-sjv-1}$. Therefore, for some small but fixed $\delta$, the inequality $\eqref{eq-py-sjv-1}$ can be rewritten as
\begin{equation}\label{eq-py-sjl-1-1}
	\int_{\mathbbm{R}}\int_{\mathbbm{R}_+} ((\Delta {\rm{div}}{p}_y)^2+2|\nabla{\rm{div}}{p}_y|^2+({\rm{div}}{p}_y)^2) \,\mathrm{d}x \mathrm{d}y\le C \| \nabla{v}_y \|_{H^1}^2+C\delta\| {\rm{div}}{p}_{yyy} \|_{L^2}^2.
\end{equation}
 Now we estimate the last term $\| {\rm{div}}{p}_{yyy} \|_{L^2} $. Rewriting $\nabla\eqref{eq-rd2-dengjia}_{2y}$ in the form $\nabla{\rm{div}}{p}_{yy}-\nabla p_{2y}=\nabla v_{yy}$ and squaring it, we obtain
\begin{equation}\label{eq-py-zk}
	|\nabla{\rm{div}}{p}_{yy}|^2+|\nabla p_{2yy}|^2+|\nabla p_{2y}|^2-2\{\nabla{\rm{div}}{p}_{y}\cdot\nabla p_{2y}\}_y=|\nabla v_{yy}|^2-2 \nabla p_{1xy}\cdot \nabla p_{2yy}.
\end{equation}
Integrating the above equation over $\mathbbm{R}_+\times\mathbbm{R}$, we can get
\begin{equation}\label{eq-py-sjl-2}
	\begin{aligned}[b]
		&\int_{\mathbbm{R}}\int_{\mathbbm{R}_+} (|\nabla{\rm{div}}{p}_{yy}|^2+|\nabla p_{2yy}|^2+|\nabla p_{2y}|^2) \,\mathrm{d}x \mathrm{d}y \\
        =& \int_{\mathbbm{R}}\int_{\mathbbm{R}_+} |\nabla v_{yy}|^2 \,\mathrm{d}x \mathrm{d}y -2 \int_{\mathbbm{R}}\int_{\mathbbm{R}_+} (p_{1xxy}p_{2xyy}+p_{1xyy}p_{2yyy}) \,\mathrm{d}x \mathrm{d}y \\
        \le& \int_{\mathbbm{R}}\int_{\mathbbm{R}_+} |\nabla v_{yy}|^2 \,\mathrm{d}x \mathrm{d}y + \int_{\mathbbm{R}}\int_{\mathbbm{R}_+} |\nabla{\rm{div}}{p}_{y}|^2 \,\mathrm{d}x \mathrm{d}y,
	\end{aligned}
\end{equation}
where we have used the fact that
$$2p_{1xxy}p_{2xyy}+2p_{1xyy}p_{2yyy}\le (p_{1xxy}+p_{2xyy})^2+(p_{1xyy}+p_{2yyy})^2=|\nabla {\rm{div}}{p}_{y}|^2.$$
 Adding $\eqref{eq-py-sjl-2}$ and $\eqref{eq-py-sjl-1-1}$, for some small $\delta$, we can obtain
\begin{equation}\label{eq-py-sjl-2-1}
	\begin{aligned}
		\int_{\mathbbm{R}}\int_{\mathbbm{R}_+} (|\nabla^2{\rm{div}}{p}_{y}|^2+|\nabla{\rm{div}}{p}_y|^2+({\rm{div}}{p}_y)^2+|\nabla p_{2yy}|^2+|\nabla p_{2y}|^2) \,\mathrm{d}x \mathrm{d}y
		\le C \| \nabla{v}_y \|_{H^1}^2.
	\end{aligned}
\end{equation}
Here, we have used the fact that $|\nabla^2 {\rm{div}}{p}_y|^2\le C(\Delta {\rm{div}}{p}_{y})^2+C|\nabla{\rm{div}}{p}_{yy}|^2$. According to $\eqref{eq-rd2-dengjia}$, we can deduce that
\begin{equation}\label{eq-py-sjl-2-2}
	\begin{aligned}[b]
		&\| \nabla^k p_{1y}\|_{L^2}^2\le C(\| \nabla^k {\rm{div}}{p}_{xy} \|_{L^2}^2+\| \nabla^k v_{xy} \|_{L^2}^2)\le C \| \nabla{v}_y \|_{H^1}^2,\quad k=0,1,\\
		&\| \nabla^k p_{2}\|_{L^2}^2\le C(\| \nabla^k {\rm{div}}{p}_{y} \|_{L^2}^2+\| \nabla^k v_{y} \|_{L^2}^2)\le  C \| \nabla{v}_y \|_{H^1}^2+C \| \nabla^k{v}_y \|_{L^2}^2,\quad k=0,1,2.
	\end{aligned}
\end{equation}
Furthermore, using the relation between $p_1$ and $p_2$, that is $\eqref{eq-p1p2}$ and ${\rm{div}}{p}=p_{1x}+p_{2y}$, we have
\begin{equation}
	\begin{aligned}
		&p_{2xxx}=p_{1xxy}={\rm{div}}{p}_{xy}-p_{2xyy}={\rm{div}}{p}_{xy}-p_{1yyy},\\
		&p_{2xxy}=p_{1xyy}={\rm{div}}{p}_{yy}-p_{2yyy}.
	\end{aligned}
\end{equation}
It follows from $\eqref{eq-py-sjl-2-1}$ that
\begin{equation}\label{eq-py-sjl-2-4}
	\| \nabla^3 p_2 \|_{L^2}^2+\| \nabla^2 p_{1y} \|_{L^2}^2\le C \| \nabla{v}_y \|_{H^1}^2.
\end{equation}
Thus, combining $\eqref{eq-vy-gammal-jl}$ and $\eqref{eq-vy-bj-sjl}$, we can get from $(1+t)^{\alpha+1+\varepsilon}[\eqref{eq-py-sjl-2-1}+\eqref{eq-py-sjl-2-2}+\eqref{eq-py-sjl-2-4}]$ that
\begin{equation*}
	\begin{aligned}[b]
		(1+t)^{\alpha+1+\varepsilon}(\| {\rm{div}}{p}_y \|_{H^2}^2+\| p_{1y}\|_{H^2}^2+\| p_{2}\|_{H^3}^2)
		\le C(1+t)^\varepsilon M_\alpha^2,
	\end{aligned}
\end{equation*}
which completes the proof of Lemma $\ref{lem-py}$.
 $ \hfill\Box$

 \vspace{6mm}

\noindent {\bf Acknowledgements:}
The research was supported by the National Natural Science Foundation of China $\#$11771150, $\#$11831003,
Guangdong Basic and Applied Basic Research Foundation $\#$2020B1515310015.

\end{document}